\documentclass[10pt]{article}
\usepackage{amsthm,amsfonts,amsmath,amscd,amssymb,latexsym,epsfig}
\usepackage[all]{xy}
\title{A construction of the Deligne--Mumford orbifold}

\author{Joel W. Robbin\\
        University of Wisconsin
        \and
        Dietmar~A.~Salamon\\
        ETH-Z\"urich}

\date{12 October 2005}

\newtheorem{PARA}{}[section]
\newtheorem{theorem}[PARA]{Theorem}
\newtheorem{corollary}[PARA]{Corollary}
\newtheorem{lemma}[PARA]{Lemma}
\newtheorem{proposition}[PARA]{Proposition}
\newtheorem{definition}[PARA]{Definition}
\newtheorem{remark}[PARA]{Remark}
\newtheorem{example}[PARA]{Example}
\newcommand{\para}{\begin{PARA}\rm}
\newcommand{\arap}{\end{PARA}\rm}
\newcommand{\dfn}{\begin{definition}\rm}
\newcommand{\nfd}{\end{definition}\rm}
\newcommand{\rmk}{\begin{remark}\rm}
\newcommand{\kmr}{\end{remark}\rm}
\newcommand{\xmpl}{\begin{example}\rm}
\newcommand{\lpmx}{\end{example}\rm}
\newcommand{\jdef}[1]{{\bf #1}}
\newcommand{\ds}{\displaystyle}

\newcommand{\Times}[2]{{}_{#1}\!\!\times_{#2}}

\newcommand{\cB}{\mathcal{B}}

\newcommand{\cD}{\mathcal{D}}
\newcommand{\cF}{\mathcal{F}}
\newcommand{\cE}{\mathcal{E}}
\newcommand{\cG}{\mathcal{G}}

\newcommand{\cJ}{\mathcal{J}}

\newcommand{\cL}{\mathcal{L}}
\newcommand{\cM}{\mathcal{M}}

\newcommand{\cN}{\mathcal{N}}
\newcommand{\cP}{\mathcal{P}}
\newcommand{\cQ}{\mathcal{Q}}
\newcommand{\cR}{\mathcal{R}}
\newcommand{\cS}{\mathcal{S}}
\newcommand{\cT}{\mathcal{T}}
\newcommand{\cU}{\mathcal{U}}
\newcommand{\cV}{\mathcal{V}}
\newcommand{\cW}{\mathcal{W}}
\newcommand{\cX}{\mathcal{X}}
\newcommand{\cY}{\mathcal{Y}}
\newcommand{\cZ}{\mathcal{Z}}

\newcommand{\Rectangle}[8]{
\begin{array}{ccc}
{#1} & \mapright{#2} & {#3}\\
\mapdown{#4} && \mapdown{#5} \\
{#6} & \mapright{#7} & {#8}
\end{array}
}
\newcommand{\mapdown}[1]{\Big\downarrow
\rlap{$\vcenter{\hbox{$\scriptstyle#1$}}$}}
\newcommand{\mapright}[1]{\smash{
\mathop{\longrightarrow}\limits^{#1}}}
\newcommand{\medskipnoindent}{\par\medskip\par\noindent}
\newcommand{\dbar}{{\bar\partial}}
\newcommand{\one}{{{\mathchoice \mathrm{ 1\mskip-4mu l} \mathrm{ 1\mskip-4mu l}
\mathrm{ 1\mskip-4.5mu l} \mathrm{ 1\mskip-5mu l}}}}
\newcommand{\A}{{\mathbb{A}}}

\newcommand{\C}{{\mathbb{C}}}
\newcommand{\D}{{\mathbb{D}}}

\renewcommand{\H}{{\mathbb{H}}}

\newcommand{\Q}{{\mathbb{Q}}}
\newcommand{\R}{{\mathbb{R}}}

\newcommand{\T}{{\mathbb{T}}}

\newcommand{\Z}{{\mathbb{Z}}}
\newcommand{\coker}{\mathrm{ coker }}  
\newcommand{\im}{\mathrm{ im }}        
\newcommand{\id}{\mathrm{ id}}         

\newcommand{\codim}{\mathrm{ codim}}   
\newcommand{\INT}{\mathrm{ int}}       
\newcommand{\INDEX}{\mathrm{ index}}   
\newcommand{\Aut}{\mathrm{ Aut}}          
\newcommand{\Diff}{\mathrm{ Diff}}        
\newcommand{\Vect}{\mathrm{ Vect}}        
\newcommand{\Hol}{\mathrm{ Hol}}          
\newcommand{\End}{\mathrm{ End}}          
\newcommand{\eps}{{\varepsilon}}
\newcommand{\Cinf}{C^{\infty}}

\def\NABLA#1{{\mathop{\nabla\kern-.5ex\lower1ex\hbox{$#1$}}}}
\def\Nabla#1{\nabla\kern-.5ex{}_{#1}}
\def\Tabla#1{\Tilde\nabla\kern-.5ex{}_{#1}}
\renewcommand{\Tilde}{\widetilde}

\newcommand{\half}{\mbox{$\frac12$}}
\newcommand{\p}{{\partial}}

\begin{document}

\maketitle



\tableofcontents


\section{Introduction}

According to  Grothendieck~\cite{GROTHENDIECK}, a moduli space is a
space whose elements may be viewed as orbits of a groupoid.\footnote{
By  the term {\em groupoid}, we understand
a category all of whose morphisms are isomorphisms.
}
In this paper,  the main focus is on the
{\em Riemann moduli space} $\cM_g$ of (closed) Riemann surfaces of genus $g$
and various related moduli spaces.
We characterize
the Deligne--Mumford compactification $\bar\cM_g$
by a universal mapping property thus showing that it
is (canonically) an orbifold.
We also treat the related moduli spaces  $\cM_{g,n}$ and  $\bar\cM_{g,n}$.

The points in the moduli space $\cM_g$
are in bijective correspondence with equivalence classes of Riemann surfaces
where two Riemann surfaces are equivalent iff there is an isomorphism
(holomorphic diffeomorphism)\footnote{In the sequel,
when no confusion can result, we will use the term
{\em isomorphism} to signify any bijection between
sets which preserve the appropriate structures.
}
between them; i.e. the Riemann
moduli space is the orbit space of the groupoid whose objects are Riemann
surfaces and whose morphisms are these isomorphisms.
For applications it is important to refine
these groupoids by considering Riemann surfaces with {\em marked points}.
An object is now a {\em marked Riemann surface} of {\em type} $(g,n)$,
i.e. a Riemann surface  of genus $g$  equipped with a sequence
of $n$ distinct points in that surface.
An isomorphism is an isomorphism of Riemann surfaces which
carries the sequence of marked points in the
source to the sequence in the target preserving the indexing.
The corresponding moduli space is denoted $\cM_{g,n}$
and of course $\cM_{g,0}=\cM_g$.

A Riemann surface is a smooth surface $\Sigma$ equipped with a
complex structure $j$. Since any two smooth surfaces of the same
genus are diffeomorphic we may define the Riemann moduli space
as the orbit space under the action of the diffeomorphism group
$\Diff(\Sigma)$ of the space $\cJ(\Sigma)$ of complex structures
$j$ on $\Sigma$:
$$
         \cM_g:= \cJ(\Sigma)/\Diff(\Sigma).
$$
The result is independent of the choice of the substrate $\Sigma$
in the sense that any diffeomorphism $f:\Sigma\to\Sigma'$
induces a bijection $\cJ(\Sigma)\to\cJ(\Sigma')$
and a group isomorphism $\Diff(\Sigma)\to\Diff(\Sigma')$ intertwining
the group actions.
Similarly a marked\footnote{
The reader is cautioned
that the term {\em marked Riemann surface} is often
used with another meaning in the literature.
}
Riemann surface is a triple $(\Sigma,s_*,j)$
where $s_*$ is a finite sequence of $n$ distinct points of $\Sigma$
(i.e. $s_*\in\Sigma^n\setminus\Delta$ where $\Delta$ is the ``fat'' diagonal)
so the corresponding moduli space is
$$
\cM_{g,n}:= \bigl(\cJ(\Sigma)\times(\Sigma^n\setminus
\Delta)\bigr)/\Diff(\Sigma).
$$
This can also be written as
$$
\cM_{g,n}= \cJ(\Sigma)/\Diff(\Sigma,s_*)
$$
where $\Diff(\Sigma,s_*)$ is the subgroup of diffeomorphisms which
fix the points of some particular sequence $s_*$.
Thus in these cases we can replace the groupoid by a group action;
the objects are the points of $\cJ(\Sigma)$.

An object in a groupoid is called {\em stable} iff
its  automorphism group is finite.
A marked Riemann surface  of type  $(g,n)$ is stable if and only if
 $n>\chi(\Sigma)$ where $\chi(\Sigma)=2-2g$  is the Euler characteristic.
In this case each  automorphism group is finite,
but (in the case $g\ge 1$) may be nontrivial.
However, the only automorphism isotopic to
the identity is the identity itself so the identity component
$\Diff_0(\Sigma)$ of $\Diff(\Sigma)$ acts freely on
$\cJ(\Sigma)\times(\Sigma^n\setminus\Delta)$.
The corresponding orbit space
$$
\cT_{g,n}:= \bigl(\cJ(\Sigma)\times
(\Sigma^n\setminus\Delta)\bigr)/\Diff_0(\Sigma)
$$
is called {\em Teichm\"uller space}.
In~\cite{EE} Earle  and  Eells showed
that the projection $\cJ(\Sigma)\to\cT_g$ is a principal fiber bundle
with structure group $\Diff_0(\Sigma)$ and that the base $\cT_g$ is a finite
dimensional smooth manifold of real dimension $6g-6$.
In other words, through each $j\in\cJ(\Sigma)$ there is a smooth slice
for the action of $\Diff_0(\Sigma)$.
(Similar statements hold for $\cT_{g,n}$.)
The total space $\cJ(\Sigma)$ is a complex manifold; the tangent space
at a point $j\in\cJ(\Sigma)$ is the space
$$
T_j\cJ(\Sigma)=\Omega^{0,1}_j(\Sigma,TM)
:=\{\hat j\in\Omega^0(\Sigma,\End(T\Sigma)): j\hat j + \hat j j =0\}
$$
of $(0,1)$ forms on $(\Sigma,j)$ with values in the tangent bundle.
This tangent space is  clearly a complex vector space (the complex
structure is $\hat j\mapsto j\hat j$) and it is not hard to
show (see e.g.~\cite{TROMBA} or Section~\ref{sec:cx}) that this
almost complex structure on $\cJ(\Sigma)$ is integrable and that the
action admits a holomorphic\footnote{
At this point in the discussion this means that the slice is a
complex submanifold of $\cJ(\Sigma)$. After we define the complex
structure on the base a holomorphic slice will be the same thing
as the image of a holomorphic section.
}
slice  through  every point. Since the action of $\Diff_0(\Sigma)$
is (tautologically) by holomorphic diffeomorphisms  of $\cJ(\Sigma)$,
this defines a complex structure on the base $\cT_g$ which is
independent of the choice of the local slice used to define it.
Thus $\cT_g$ is a complex manifold of dimension $3g-3$. Again,
similar results hold for $\cT_{g,n}$. Earle and Eells also showed
that all three spaces in the fibration
$$
\Diff_0(\Sigma)\to\cJ(\Sigma)\to\cT_g\eqno(EE)
$$
are contractible so that the fibration is smoothly trivial and has a
(globally defined) smooth  section. In~\cite{EARLE} Earle showed that
there is no global holomorphic section of $\cJ(\Sigma)\to\cT_g$.
The monograph of Tromba~\cite{TROMBA} contains a nice exposition
of this point of view (and more) and the anthology~\cite{DISCRETE}
is very helpful for understanding the history of the subject and
other points of view.

Now we take a different point of view.
An {\em unfolding} is the germ of a pair $(\pi_A,a_0)$ where
$\pi_A:P\to A$ is a Riemann family and $a_0$ is a point of $A$.
(The term  {\em Riemann family} means that
$\pi_A$ is a proper holomorphic map
and $\dim_C(P)=\dim_C(A)+1$.
The term {\em germ} means that we do not distinguish
between $(\pi_A,a_0)$ and the unfolding which results
by replacing $A$ by a neighborhood of $a_0$ in $A$.)
The fibers $P_a:=\pi^{-1}(a)$ are then complex curves.
The fiber $P_{a_0}$ is called the {\em central fiber}.
A {\em morphism} of unfoldings is a commutative diagram
$$
 \Rectangle{P}{\Phi}{Q}{\pi_A}{\pi_B}{A}{\phi}{B}
$$
where $\Phi$ and $\phi$ are holomorphic, $\phi(a_0)=b_0$
and, for each $a\in A$,
the restriction of $\Phi$ to the fiber $P_a$ is  an isomorphism.
Again,  this is to be understood in the sense of germs:
$\phi$  need only be defined on a neighborhood of $a_0$
and two morphisms are the same iff they
agree on a smaller neighborhood of $a_0$.
An unfolding $(\pi_B:Q\to B,b_0)$ is called {\em universal}
iff for every other unfolding $(\pi_A,a_0)$  every isomorphism
$f:P_{a_0}\to Q_{b_0}$ extends uniquely to a morphism $(\phi,\Phi)$
from $(\pi_A,a_0)$ to $(\pi_B,b_0)$. From  the uniqueness of
the extension it follows that any two universal unfoldings
with the same central fiber are isomorphic in the obvious sense.

Now assume that $\pi_A$ is a submersion so that the fibers are
Riemann surfaces. Using the holomorphic slices for the principal
fiber bundle $\cJ(\Sigma)\to\cT_g$ it is not hard to construct
a universal unfolding of any Riemann surface of genus $\ge 2$
(see Section~\ref{sec:teichmuller}).
Again similar results hold for
$\cT_{g,n}$ (see Section~\ref{sec:teichmuller-n}).

The spaces $\cM_{g,n}$ are not compact.
The {\em Deligne--Mumford} moduli space $\bar\cM_{g,n}$ defined in
Section~\ref{sec:dm-orbifold} is a compactification of $\cM_{g,n}$.
The objects in the corresponding groupoid are commonly called
{\em stable curves of type $(g,n)$}.
Two such curves  need not be homeomorphic.
This moduli space is still the orbit space
of a groupoid but not (in any obvious way) the orbit space of a group action.
 We will characterize $\bar\cM_{g,n}$ by the universal mapping property,
but we will word the definitions so as to avoid the complexities of
algebraic geometry and singularity theory.

It is a well known theorem of algebraic geometry
that a complex curve  $C$ admits a desingularization $u:\Sigma\to C$.
This means that $\Sigma$ is a Riemann surface and that
the restriction of $u$ to the set of  regular points of $u$
is a holomorphic diffeomorphism onto the set of smooth points of the curve $C$.
The desingularization is unique in the sense that if
$u':\Sigma'\to C$ is another desingularization, the
holomorphic diffeomorphism $u^{-1}\circ u'$
extends to a holomorphic diffeomorphism $\Sigma'\to\Sigma$.
A marked complex curve is one which is equipped with a finite sequence
of distinct smooth points.
 A desingularization pulls pack the marking to a marking of $\Sigma$.
That  a marked complex curve $C$ is of {\em type $(g,n)$} means
that the arithmetic genus (see Definition~\ref{def:arithGenus})
of $C$ is $g$ and the number of marked points is $n$.
A {\em nodal curve} is a complex curve with at worst nodal singularities.
For a nodal curve the desingularization $u$ is an immersion and the
critical points occur in pairs. This equips $\Sigma$ with what we call
a {\em nodal structure}. In Section~\ref{sec:RiemannSurface} we use the
term {\em marked nodal Riemann surface} to designate a surface $\Sigma$
with these additional structures. A {\em stable curve} is a marked nodal
curve whose  corresponding marked nodal Riemann surface has a finite
automorphism group. The main result of this paper extends the universal
unfolding construction from the groupoid of stable Riemann surfaces
to the groupoid of stable marked nodal Riemann surfaces.

\begin{quote}
\bigskip\noindent{\bf Theorem~A.} \em
A marked nodal Riemann surface admits a
universal unfolding if and only if it is stable.
\bigskip
\end{quote}

This theorem is an immediate consequence of
Theorems~\ref{thm:idm} and~\ref{thm:existence} below.
To avoid the intricacies of singularity theory
our precise definitions
(see Sections~\ref{sec:nodal-fam} and~\ref{sec:universal})
involve only what we call {\em nodal families}.
However, it is well known that (near its central fiber) an unfolding
is a submersion if and only if its central fiber is a smooth complex curve
and is a nodal family  if and only if its central fiber is a nodal curve.

Now we describe the proof. First we consider  the case of a Riemann surface
without marked points or nodal points. In this case the sequence~$(EE)$
is a principal bundle if and only if $g\ge 2$, i.e. if and only if any Riemann
surface of genus $g$ is stable. Abbreviate
$$
\cD_0:=\Diff_0(\Sigma), \qquad \cJ:=\cJ(\Sigma), \qquad
\cT:=\cT(\Sigma):=\cJ(\Sigma)/\Diff_0(\Sigma).
$$
Thus $\cT_g:=\cT$ is  Teichm\"uller space
and the principal fiber bundle~$(EE)$ takes the form
$$
\cD_0 \to\cJ\to\cT.
$$
The associated fiber bundle
$$
\pi_{\cT}:\cQ := \cJ\times_{\cD_0}\Sigma\to\cT
$$
has fibers isomorphic to $\Sigma$. It is commonly called the
{\em universal curve of genus $g$ over Teichm\"uller space}.
Choose a Riemann surface $(\Sigma,j_0)$ and a holomorphic slice $B\subset\cJ$
through $j_0$. Let
$$
       \pi_B:Q\to B
$$
be the restriction to $B$ of the pull back of the bundle
$\pi_\cT$ to its total space. As $B$ is a slice, the
projection $\pi_B$ is a trivial bundle (in the smooth sense).
The map $\pi_B$ is a holomorphic submersion.
In Section~\ref{sec:teichmuller} we show that it is a universal
unfolding of $j_0$. Here's why $(\pi_B,j_0)$ is universal.
Let $\pi_A:P\to A$ be a holomorphic submersion whose fiber has genus $g$
and whose central fiber over $a_0\in A$
is isomorphic to $(\Sigma,j_0)$.
As a smooth map $\pi_A$ is trivial so after shrinking $A$ we have a smooth
local trivialization $\tau:A\times\Sigma\to P$.
Write $\tau_a(z):=\tau(a,z)$ for $a\in A$
so $\tau_a$ is a diffeomorphism from $\Sigma$ to $P_a$.
Denote the pull back by $\tau_a$    of the complex
structure on $P_a$ by $j_a$, i.e. $\tau_a:(\Sigma,j_a)\to P_a$
is an isomorphism.  As $B$ is a slice
we can modify  the trivialization $\tau$  so $j_a\in B$.
The equation $\phi(a)=j_a$ defines a map $\phi:A\to B$.
Using the various trivializations we then get
a morphism $(\phi,\Phi)$  from $\pi_A$ to $\pi_B$.
In Section~\ref{sec:teichmuller}  we show that these maps are holomorphic.
In Section~\ref{sec:teichmuller-n} we carry out the
analogous construction for $\cT_{g,n}$.

It is now clear that $(\pi_A,a_0)$ is universal if and only if
$\phi:(A,a_0)\to(B,b_0)$ is the germ of
a diffeomorphism. By the inverse function theorem
this is so if and only if the linear operator
$d\phi(a_0):T_{a_0}A\to T_{b_0}B$
is invertible. This condition can be formulated as the unique
solvability of a partial differential equation on $P_{a_0}$;
we call an unfolding {\em infinitesimally universal} when it
satisfies this unique solvability condition. The crucial point is
that infinitesimal universality is meaningful even for nodal families,
i.e. when there is no analog of the Earle--Eells principal fiber bundle.
But we still have the following

\begin{quote}
\bigskip\noindent{\bf Theorem~B.}
\em A nodal unfolding is universal if and only if
it is infinitesimally universal.
\bigskip
\end{quote}

This is restated as Theorem~\ref{thm:idm} below. Here is the idea of the proof.
Let ${(\pi_A:P\to A,a_0)}$ and $(\pi_B:Q\to B,b_0)$  be nodal unfoldings and
${f_0:P_{a_0}\to Q_{b_0}}$ be an isomorphism of the central fibers.
For simplicity assume there
is at most one critical point in each fiber and no marked points.
Essentially by  the definition of nodal unfolding there is a neighborhood
$N$ of the set of critical points such that for $a\in A$
the  intersection $N_a:=N\cap P_a$ admits an isomorphism
$$
N_a\cong \{(x,y)\in\D^2: xy=z\}
$$
where $\D$ is the closed unit disk in $\C$
and $z=z(a)\in\D$. Thus if $z(a)\ne 0$
the fiber $N_a$ is an annulus whereas if
$z(a)=0$ it is a pair of transverse disks.
In either case the boundary is a disjoint
union  $(\p\D \sqcup\p\D)$ of two copies
of the circle $S^1:=\p\D$.
The map $N\to A$ is therefore not trivializable
as the topology of the fiber changes.
However, the bundle $\p N\to A$ is trivializable;
choose a trivialization $A\times(\p\D \sqcup\p\D) \to \p N$.
Using this trivialization we will define
(see Section~\ref{sec:hardy}) manifolds of maps
$$
\cW:=\bigsqcup_{a\in A}\cW_a, \qquad \cW_a:=\bigsqcup_{b\in B}\cW(a,b),
\qquad \cW(a,b):=\mathrm{Map}(\p N_a,Q_b\setminus C_B)
$$
where $C_B$ is the set of critical points of $\pi_B$ and
$\bigsqcup$ denotes disjoint union. Let $\cU_a\subset\cW$ be the
set of all maps in $\cW_a$ which extend to a holomorphic map
$N_a\to Q$ and $\cV_a\subset\cW$ be the set of all maps in $\cW_a$
which extend to a holomorphic map $P_a\setminus N_a\to Q$.
We will replace $A$ and $\cW$ by smaller neighborhoods of $a_0$
and $f_0|\p N_{a_0}$ as necessary.
We show that $\cU_a$ and $\cV_a$ are submanifolds of $\cW_a$.
It is not too hard to show that the unfolding $(\pi_B,b_0)$ is universal
if and only if the manifolds $\cU_a$ and $\cV_a$ intersect in a unique point:
the morphism $(\phi,\Phi):(\pi_A,a_0)\to(\pi_B,b_0)$ is then defined
so that this intersection point $\gamma$  lies in the fiber $\cW_{\phi(a)}$
and  $\Phi_a$ is the unique holomorphic map extending $\gamma$.
We will see that the unfolding $(\pi_B,b_0)$ is infinitesimally universal
if and only if (for all $(\pi_A,a_0)$ and $f_0$) the corresponding
infinitesimal  condition
$$
T_{\gamma_0} \cW_{a_0} = T_{\gamma_0} \cU_{a_0} \oplus T_{\gamma_0} \cV_{a_0}
$$
holds where $\gamma_0=f_0|\p N_{a_0}$. This Hardy space decomposition
is reminiscent of the construction of the moduli space of holomorphic
vector bundles explained by Pressley \& Segal in~\cite{PS}.

We have already explained why smooth marked Riemann surfaces have
universal unfoldings. It is now easy to construct  a universal unfolding
of a stable marked nodal Riemann surface: it is constructed from a
universal unfolding  for the marked Riemann surface that results by
replacing each nodal point by a marked point. Such an unfolding
is a triple $(\pi, S_*,b_0)$ where $\pi:Q\to B$ is a nodal family,
$S_*$ is a sequence of holomorphic sections of $\pi$ corresponding
to the marked points, and $b_0\in B$. We call a pair $(\pi,S_*)$
a {\em universal family of type $(g,n)$} iff (1)~$(\pi, S_*,b_0)$
is a universal unfolding for each $b_0\in B$ and (2)~every marked
nodal Riemann surface of type $(g,n)$ occurs as the domain
of a desingularization
of some fiber $Q_b$, $b\in B$. Theorem~\ref{thm:stable}
(openness of transversality) says that if $(\pi, S_*,b_0)$
is an infintesimally  universal unfolding so is
$(\pi, S_*,b)$ for $b$ near $b_0$.
Together with Theorems~A and~B this implies

\begin{quote}
\bigskip\noindent{\bf Theorem~C.}
\em If $n>2-2g$ there exists a universal family of type $(g,n)$.
\bigskip
\end{quote}

This is restated as Proposition~\ref{prop:universal} below. It is not asserted that $B$ is connected.
Rather, the universal family should be viewed as a generalization of the
notion of an atlas for a manifold. This generalization is called
an {\em etale groupoid}. The Deligne Mumford orbifold $\bar\cM_{g,n}$ is
then the orbit space of this groupoid and the definitions
are arranged so that the orbifold structure is independent of
the choice of the universal family used to define it.
See Section~\ref{sec:dm-orbifold}.

A consequence of our theorems is that other constructions
of the Deligne--Mumford moduli space (and in particular of
the Riemann moduli space) which have the universal unfolding
property give the same space.
However, in the case of a construction  where
the moduli space is given only a topology
(or a notion of convergence of sequences as in~\cite{HUMMEL})
we show that  the topology determined by our construction
agrees with the topology of the other construction
(see Section~\ref{sec:topology}).  In Section~\ref{sec:compact}
we prove that our  $\bar\cM_{g,n}$ is compact and Hausdorff
by adapting the arguments of the monograph of Hummel~\cite{HUMMEL}.

\medskip\noindent{\bf Acknowledgement.}
JWR  would like to
thank Yongbin Ruan for helpful discussions
and Ernesto Lupercio and Benardo Uribe
for helping him understand the concept of {\em etale groupoid}.

\medskip\noindent{\bf Notation.}
Throughout  the closed unit disk in the complex plane  is denoted by
$$
\D:=\left\{z\in\C\,:\,|z|\le1\right\}
$$
and its interior is denoted by $\INT(\D):=\left\{z\in\C\,:\,|z|<1\right\}$.
Thus $S^1:=\p\D$ is the unit circle. Also
$$
\A(r,R):=\{z\in\C: r\le |z|\le R\}
$$
denotes the closed annulus with inner radius $r$ and outer radius $R$.


\section{Orbifold structures}\label{sec:orbifold}

In this section we review orbifolds.  Our definitions are arranged so as to
suit our ultimate objective of defining an orbifold structure on the 
Deligne--Mumford moduli space. 

\para A \jdef{groupoid} is a category in which every morphism
is an isomorphism. Let $B$ be the set of objects of a groupoid
and $\Gamma$ denote the set of (iso)morphisms. For $a,b\in B$
let $\Gamma_{a,b}\subset \Gamma$ denote the isomorphisms from $a$ to $b$;
the group
$$
\Gamma_a:=\Gamma_{a,a}
$$
is called the \jdef{automorphism group}\footnote{
Also commonly called the {\em isotropy group} or {\em stabilizer group}.
}
of $a$.
The groupoid is called \jdef{stable} iff  every automorphism group is finite.
Define the \jdef{source} and \jdef{target} maps
$s,t:\Gamma\to B$ by
$$
s(g)=a \mbox{ and }  t(g)=b \iff g\in \Gamma_{a,b}.
$$
The map  $e:B\to \Gamma$
which assigns to each object $a$ the identity morphism of $a$ is called
the \jdef{identity section} of the groupoid and
the map  $i:\Gamma\to \Gamma$
which assigns to each morphism $g$ its inverse $i(g)=g^{-1}$ is called
the \jdef{inversion map}.
Define the set $\Gamma\Times{s}{t} \Gamma$ of \jdef{composable pairs} by
$$
\Gamma\Times{s}{t} \Gamma=\{(g,h)\in \Gamma\times \Gamma: s(g)=t(h)\}.
$$
The map  $m:\Gamma\Times{s}{t} \Gamma\to \Gamma$ which
assigns to each composable pair  the composition  $m(g,h)=gh$ is called
the \jdef{multiplication map}.
The five maps $s$, $t$, $e$, $i$, $m$ are called the \jdef{structure maps}
of the groupoid. Note that
$$
\Gamma_{a,b}= (s\times t)^{-1}(a,b).
$$
We denote the \jdef{orbit space}
of the groupoid $(B,\Gamma)$ by $B/\Gamma$:
$$
 B/\Gamma:=\{[b]:b\in B\}, \qquad [b]:=\{ t(g)\in B: g\in\Gamma, \;s(g)=b\}.
$$
\arap

\para\label{lie}
A  \jdef{Lie groupoid} is a groupoid $(B,\Gamma)$ such that
$B$ and $\Gamma$ are smooth manifolds\footnote{
For us a manifold is always second countable and Hausdorff,
unless otherwise specified.
},
the structure maps  are smooth, and
the map  $s:\Gamma\to B$ (and hence also the map $t=s\circ i$)
is a submersion. (The latter condition implies that
$\Gamma\Times{s}{t} \Gamma$ is a submanifold of $\Gamma\times \Gamma$
so that the condition that $m$ be smooth is meaningful.)
A \jdef{homomorphism}
from a Lie groupoid  $(B,\Gamma)$
to a Lie groupoid $(B',\Gamma')$
is a smooth functor,
i.e. a pair of smooth maps $B\to B'$ and $\Gamma\to\Gamma'$,
both denoted by $\iota$, which intertwine the structure maps:
$$
s'\circ\iota=\iota\circ s,\qquad
t'\circ\iota=\iota\circ t,\qquad
e'\circ\iota=\iota\circ e,
$$
$$
i'\circ \iota=\iota\circ i,\qquad m'\circ(\iota\times\iota)=\iota\circ m.
$$
(The first two of these five conditions imply that
$(\iota\times\iota)(\Gamma\;\Times{s}{t}\Gamma)
\subset \Gamma'\;\Times{s'}{t'}\Gamma'$
so that the fifth condition is meaningful.)
Similar definitions are used in the complex category
reading {\em complex} for {\em smooth}  (for manifolds)
or {\em holomorphic} for smooth (for maps).
A Lie groupoid $(B,\Gamma)$ is called \jdef{proper} if
the map $s\times t:\Gamma\to B\times B$ is proper.
\arap

\para\label{etale}
An \jdef{etale  groupoid} is a Lie groupoid $(B,\Gamma)$
such that the map $s:\Gamma\to B$
(and hence also the map $t=s\circ i$) is a local diffeomorphism.
A proper etale groupoid is automatically stable.
A homomorphism $\iota:(B,\Gamma)\to (B',\Gamma')$ of etale
groupoids is called a \jdef{refinement} iff the following holds.
\begin{enumerate}
\item[(i)]
The induced map $\iota_*:B/\Gamma\to B'/\Gamma'$ on orbit spaces
is a bijection.
\item[(ii)]
For all $a,b\in B$, $\iota$ restricts to a bijection
$\Gamma_{a,b}\to\Gamma'_{\iota(a),\iota(b)}$.
\item[(iii)]
The map on objects (and hence also the map on morphisms)
is  a local diffeomorphism.
\end{enumerate}
Two proper etale groupoids are called \jdef{equivalent} iff they have
a common proper refinement.
\arap

\dfn\label{def:orbifoldStructure}
Fix an abstract groupoid $(\cB,\cG)$.  This groupoid
is to be viewed as the ``substrate''  for an additional  structure
to be imposed; initially it does not even have a topology.
Indeed, the definitions are worded so as to allow for the possibility
that $\cB$ is not even a set but a proper class in the sense
of  G\"odel Bernays set theory (see~\cite{KELLEY}).

An \jdef{orbifold  structure} on the groupoid $(\cB,\cG)$
is a functor $\sigma$ from a proper etale groupoid
$(B,\Gamma)$ to $(\cB,\cG)$ such that
\begin{enumerate}
\item[(i)] $\sigma$ induces a bijection
$B/\Gamma\to\cB/\cG$ of orbit spaces, and
\item[(ii)] for all $a,b\in B$, $\sigma$ restricts to a bijection
$\Gamma_{a,b}\to\cG_{\sigma(a),\sigma(b)}$.
\end{enumerate}
A \jdef{refinement}  of orbifold structures  is a refinement
$\iota:(B,\Gamma)\to (B',\Gamma')$ of proper etale groupoids
such that $\sigma=\sigma'\circ\iota$; as before we say that
$\sigma:(B,\Gamma)\to(\cB,\cG)$ is a refinement
of $\sigma':(B',\Gamma')\to(\cB,\cG)$.
Two orbifold structures are called \jdef{equivalent}
iff they have a common refinement.
An \jdef{orbifold} is an abstract groupoid
$(\cB,\cG)$ equipped with an orbifold structure
$\sigma:(B,\Gamma)\to(\cB,\cG)$.
\nfd

\xmpl\label{ex:openCover}
A smooth manifold $M$ is a special case of an orbifold as follows:
View $M=:\cB$ as a trivial groupoid, i.e. the only
morphisms are identity morphisms.
Any countable open cover $\{U_\alpha\}_{\alpha\in I}$
on $M$ determines an etale groupoid $(B,\Gamma)$ with
$$
B:=\bigsqcup_{\alpha\in I} U_\alpha,\qquad
\Gamma:=\bigsqcup_{(\alpha,\beta)\in I\times I} U_\alpha\cap U_\beta,
$$
$$
s(\alpha,p,\beta):=(\alpha,p),\qquad
t(\alpha,p,\beta):=(\beta,p),\qquad
e(\alpha,p):=(\alpha,p,\alpha),
$$
$$
i(\alpha,p,\beta):=(\beta,p,\alpha), \qquad
m((\beta,p,\gamma),(\alpha,p,\beta)):=(\alpha,p,\gamma).
$$
Here $\bigsqcup$ denotes disjoint union.
(The \jdef{disjoint union} $\bigsqcup_{\alpha\in I}X_\alpha$
of an indexed collection
$\{X_\alpha\}_{\alpha\in I}$ of sets is
the set of pairs $(\alpha,x)$ where $\alpha\in I$
and $x\in X_\alpha$.)
A refinement of open covers in the usual sense determines a
refinement of etale groupoids as in \ref{etale}.

If $\{\phi_\alpha,U_\alpha\}_{\alpha\in I}$ is a countable atlas then
an obvious modification of the above construction gives rise to
an orbifold structure on $M$ where $B$ is a disjoint union of open subsets
of Euclidean space, i.e. a manifold structure
is a special case of an orbifold structure.
\lpmx

\xmpl \label{ex:groupAction}
A Lie group action $G\to\Diff(M)$ determines a Lie groupoid $(\cB,\cG)$ where
$\cB=M$, $\cG=\{(g,a,b)\in G\times M\times M:b=g(a)\}$,
and the structure maps are defined by $s(g,a,b):=a$, $t(g,a,b):=b$,
$e(a):=(\id,a,a)$, $i(g,a,b):=(g^{-1},b,a)$, and
$m((h,b,c),(g,a,b)):=(hg,a,c)$.
The orbit space $\cB/\cG$ of this groupoid
is the same as the orbit space $M/G$ of the group action.
The condition that this groupoid be proper is the usual definition
of proper group action, i.e. the map
$G\times M\to M\times M:(g,x)\mapsto(x,g(x))$ is proper.

Assume that the action is almost free
[meaning that the isotropy group $G_p$ of each point of $M$ is finite]
and sliceable
[meaning that there is a slice through every point of $M$;
a slice is a submanifold $S\subset M$ such that there is a neighborhood
$U$ of the identity in $G$ with the property that the map
$U\times S\to M:(g,x)\mapsto g(x)$ is a diffeomorphism onto a
neighborhood of $S$ in $M$].  Now let
$$
B:=\bigsqcup_{\alpha\in I}S_\alpha
$$
be a disjoint union of slices such that every orbit passes through
at least one slice.  Let
$$
\Gamma:=\bigsqcup_{\alpha,\beta\in I}\Gamma_{\alpha\beta},\qquad
\Gamma_{\alpha\beta}:=\{(g,a,b)\in\cG: a\in S_\alpha,b\in S_\beta\}.
$$
Then $\Gamma_{\alpha\beta}$ is a submanifold of $\cG$.
Moreover, if the group action is proper, then the obvious morphism
$\sigma:(B,\Gamma)\to(\cB,\cG)$ is an orbifold structure,
and any two such orbifold structures are equivalent.
Note that, if $G$ is a discrete group acting properly on $M$,
then $S:=\cB=M$ is a slice and $\sigma:=\id$ is an orbifold structure.
\lpmx

\xmpl\label{ex:irrational}
Consider the group action where $G:=\Z$ acts on $M:=S^1$
by $(k,z)\mapsto e^{2\pi ik\omega}z$ and
$\omega\in\R\setminus\Q$ is irrational.
Then the groupoid $(\cB,\cG)$ constructed in
Example~\ref{ex:groupAction} is etale
but not proper. Note that the quotient $\cB/\cG$ is an
uncountable set with the trivial topology (two open sets).
The inclusion of any open set into $S^1$ is a refinement.
\lpmx

\xmpl\label{ex:nonHausdorff}
Consider the group action where the multiplicative group $G:=\R^*$
of nonzero real numbers acts on $M:=\R^2\setminus 0$ by
$
t\cdot (x,y) := (tx,t^{-1}y).
$
The action is free and sliceable but not proper, and
the quotient topology is non Hausdorff
(every neighborhood of $\R^*\cdot(1,0)$ intersects
every neighborhood of $\R^*\cdot(0,1)$).
The groupoid  constructed from the disjoint union
$B:=S_1\bigsqcup S_2$ of the two slices
$S_1:= \{1\}\times\R$, $S_2:=\R\times\{1\}$
is not proper. If we extend the group action, by
adjoining the map $(x,y)\mapsto(y,x)$, the
orbit space is $\R$ which is Hausdorff,
but the new group action is still not proper.
\lpmx

\para
Let $(B,\Gamma)$ be a stable etale groupoid,
$a,b\in B$, and $g\in\Gamma_{a,b}$.
Then there exist neighborhoods $U$ of $a$,
$V$ of $b$ in $B$, and $N$ of $g$ in $\Gamma$
such that $s$ maps $N$ diffeomorphically onto $U$
and $t$ maps $N$ diffeomorphically onto $V$. Define
$s_g:=s|N$, $t_g:=t|N$, and $\phi_g:=t_g\circ s_g^{-1}$.
Thus $\phi_g$ ``extends'' $g\in\Gamma_{a,b}$ to
diffeomorphism $\phi_g:U\to V$.  The following lemma says
that when $a=b$ we may choose $U=V$ independent of $g$ and
obtain an action
$$
\Gamma_a\to\Diff(U): g\mapsto\phi_g
$$
of the finite group $\Gamma_a$ on the open set $U$.
\arap

\begin{lemma} \label{le:invNhbd}
Let $(B,\Gamma)$ be a stable etale groupoid and $a\in B$.
Then there exists a neighborhood $U$ of $a$
and pairwise disjoint neighborhoods $N_g$ (for $g\in\Gamma_a$)
of $g$ in $\Gamma$
such that both $s$ and $t$ map each $N_g$ diffeomorphically onto $U$.
\end{lemma}

\begin{proof}
Choose disjoint open neighborhoods $P_g$ of $g\in\Gamma_a$
such that $s_g:=s|P_g$ and $t_g:=t|P_g$ are diffeomorphisms onto
(possibly different) neighborhoods of $a$.
By stability the group $\Gamma_a$ is finite
so there is a neighborhood $V$ of $a$ in $B$
such that $V\subset s(P_g)\cap t(P_g)$
for $g\in \Gamma_a$. Define $\phi_g:V\to B$ by $\phi_g:=t_g\circ s_g^{-1}$.
Now choose $f,g\in\Gamma_a$ and let  $h:=m(f,g)$.
We show that
\begin{equation}\label{eq:fgh}
\phi_h(x)=\phi_f\circ\phi_g(x)
\end{equation}
for $x$ in a sufficiently small neighborhood  of $a$ in $V$.
For such $x$ define  $y:=\phi_g(x)\in V$,  $z:=\phi_f(y)\in V$,
$g':=s_g^{-1}(x)\in P_g$, and $f':=s_f^{-1}(y)\in P_f$.
As $t(g')=s(f')=y$ we have $(f',g')\in \Gamma\Times{s}{t}\Gamma$, i.e.
$h':=m(f',g')$ is well defined. By continuity,
$h'\in P_h$ and $s(h')=s(g')=x$ and $t(h')=t(f')=z$,
and hence $z=\phi_h(x)$ as claimed.
Using the finiteness of $\Gamma_a$ again
we may choose a neighborhood $W$ of $a$ so
that~(\ref{eq:fgh}) holds for all $f,g\in\Gamma_a$ and all $x\in W$.
Now the intersection
$$
U:=\bigcap_{g\in\Gamma_a} \phi_g(W)\subset V
$$
satisfies $\phi_f(U)=U$ for $f\in\Gamma_a$ so $U$ and
$N_g:=s_g^{-1}(U)$ satisfy the conclusions of the lemma.
\end{proof}

\begin{corollary}\label{cor:UVN}
Let $(B,\Gamma)$ be a stable etale groupoid and $a,b\in B$.
Then there exist    neighborhoods $U$ and $V$ of $a$ and $b$ in $B$
and pairwise disjoint neighborhoods $N_f$ (for $f\in\Gamma_{a,b}$)
of $f$ in $\Gamma$ such that $s$  maps each $N_f$
diffeomorphically onto $U$ and $t$  maps each $N_f$
diffeomorphically onto $V$. The etale groupoid is proper
if and only if these neighborhoods may be chosen so that
in addition
$$
(s\times t)^{-1}(U\times V)=
\bigcup_{f\in\Gamma_{a,b}} N_f.\eqno(*)
$$
\end{corollary}

\begin{proof}
Choose disjoint neighborhoods $P_f$ of $f\in\Gamma_{a,b}$
such that $s_f:=s|P_f$ and $t_f:=t|P_f$ are diffeomorphisms
onto (possibly different) neighborhoods of $a$.
Choose $U$ as in Lemma~\ref{le:invNhbd} so small that
$U\subset s(P_f)$ for all  $f\in\Gamma_{a,b}$ and define
$\phi_f:U\to B$ by
$$
\phi_f:=t_f\circ s_f^{-1}|U.
$$
Define $N_f:=s_f^{-1}(U)$.
As in Lemma~\ref{le:invNhbd} we have $\phi_h=\phi_f\circ\phi_g$
for $g\in\Gamma_a$, $f\in\Gamma_{a,b}$, $h:=m(f,g)$,
so $t_h(N_h)=\phi_h(U)=\phi_f(U)=t_f(N_f)$.
Any two elements $h,f\in\Gamma_{a,b}$ satisfy $h=m(f,g)$
for some $g\in\Gamma_a$ so $V:=t_f(N_f)$ is independent
of the choice of $f\in\Gamma_{a,b}$ used to define it.
The condition that $s\times t$ is proper, is that for any sequence
$\{f_\nu\in\Gamma_{a_\nu,b_\nu}\}_\nu$
such that the  sequences $\{a_\nu\}_\nu$ and $\{b_\nu\}_\nu$
converge to $a$ and $b$ respectively,
the sequence $\{f_\nu\}_\nu$ has a convergent subsequence.
Condition~$(*)$ implies this as $f_\nu$ must lie in some $N_f$
for infinitely many values of $\nu$. The converse follows easily
by an indirect argument.
\end{proof}

\para
Let $(B,\Gamma)$ be an etale groupoid and equip the orbit space $B/\Gamma$
with the quotient topology, i.e. a subset of $B/\Gamma$ is open iff
its preimage under the quotient map $\pi:B\to B/\Gamma$ is open.
If $U\subset B$ is open then so is $\pi^{-1}(\pi(U))=\{t(g):g\in s^{-1}(U)\}$
so $\pi$ is an open map.
If $\iota:(B,\Gamma)\to (B',\Gamma')$ is a refinement of etale groupoids,
then the induced bijection $\iota_*:B/\Gamma\to B'/\Gamma'$ is a homeomorphism.
[The continuity of $\iota_*$ follows from the continuity of $\iota$;
the openness of $\iota_*$ follows from the openness of $\iota$ and
the fact that if $U'\subset B'$ is open then so is ${\pi'}^{-1}(\pi'(U'))$.]
Hence equivalent etale groupoids have homeomorphic orbit spaces.
It follows that the topology induced on $\cB/\cG$ by an  orbifold structure
$\sigma:(B,\Gamma)\to(\cB,\cG)$ depends only the equivalence class.
This topology is called the \jdef{orbifold topology}.
\arap

\begin{corollary}\label{cor:properHausdorff}
For a proper  etale groupoid the quotient topology on $B/\Gamma$ is Hausdorff.
\end{corollary}

\begin{proof}
In other words if $\Gamma_{a_0,b_0}=\emptyset$
then there are neighborhoods $U$ of $a$ and $V$ of $b$
such that $\Gamma(a,b)=\emptyset$ for $a\in U$ and $v\in V$.
This is a special case of Corollary~\ref{cor:UVN}.
\end{proof}


\section{Structures on surfaces}\label{sec:RiemannSurface}

The phrase \jdef{surface} means
{\em oriented smooth (i.e. $C^\infty$) manifold
of (real) dimension two, not necessarily connected}.
Unless otherwise specified all surfaces are
assume to be closed, i.e. compact and without boundary.
The structures we impose on surfaces are complex structures,
nodal structures, and point markings.  Surfaces equipped 
with these structures form the objects of a groupoid.
The objective of this paper is to equip the orbit space 
of this groupoid with an orbifold structure. 

\dfn A \jdef{Riemann surface} is a pair $(\Sigma,j)$
where $\Sigma$ is a surface
and $j:T\Sigma\to T\Sigma$
is a smooth complex structure  on $\Sigma$ which determines
the given orientation of $\Sigma$.
Since a complex structure on a surface is necessarily integrable,
a Riemann surface may be viewed as a smooth complex curve,
i.e. a compact complex manifold of (complex) dimension one.
When there is no danger of confusion
we denote a Riemann surface and its underlying surface
by the same letter.
\nfd

\dfn\label{def:nodalManifold}
A \jdef{nodal surface} is a pair $(\Sigma,\nu)$ consisting of
a surface $\Sigma$ and a set
$$
 \nu=\bigl\{\{y_1,y_2\},\{y_3,y_4\},\ldots,\{y_{2k-1},y_{2k}\}\bigr\}
$$
where $y_1,y_2,\ldots,y_{2k}$ are distinct points of $\Sigma$;
we also say $\nu$ is a \jdef{nodal structure} on  $\Sigma$.
The points $y_1,y_2,\ldots,y_{2k}$
are called the \jdef{nodal points} of the structure
and the points $y_{2j-1}$ and $y_{2j}$ are
called \jdef{equivalent nodal points}.
The nodal structure should be viewed as
an equivalence relation  on $\Sigma$ such that
every equivalence class consists of either one or two points and
only finitely many equivalence classes have two points. Hence we
often abbreviate $\Sigma\setminus\cup\nu$ by
$$
\Sigma\setminus\nu:=\Sigma\setminus\{y_1,y_2,y_3,y_4,\ldots,y_{2k-1},y_{2k}\}.
$$
\nfd

\dfn\label{def:markedNodal}
A \jdef{point marking} of a surface $\Sigma$ is a sequence
$$
       r_*=(r_1,r_2,\ldots,r_n)
$$
of distinct points of $\Sigma$;
the points $r_i$ are called \jdef{marked points}.
A \jdef{marked nodal surface} is a
triple $(\Sigma,r_*,\nu)$ where $(\Sigma,\nu)$ is
a nodal surface and  $r_*$ is a point marking of $\Sigma$
such that no marked point $r_i$ is a nodal point of $(\Sigma,\nu)$;
a \jdef{special point} of the marked nodal surface  is a point
which is either a nodal point or a marked point.
\nfd

\dfn\label{def:signature}
A marked nodal surface  $(\Sigma,r_*,\nu)$
determines a labelled graph
called the \jdef{signature} of $(\Sigma,r_*,\nu)$ as follows.
The set of vertices of the graph label the connected components
of $\Sigma$ and there is one  edge connecting vertices $\alpha$
and $\beta$ for every pair  of equivalent nodal points with one
of the points in $\Sigma_\alpha$ and the other in $\Sigma_\beta$.
More precisely, the number  of edges from $\Sigma_\alpha$ to
$\Sigma_\beta$ is the number of pairs $\{x,y\}$ of equivalent
nodal points with either $x\in \Sigma_\alpha$ and $y\in\Sigma_\beta$
or $y\in \Sigma_\alpha$ and $x\in\Sigma_\beta$. Each vertex $\alpha$
has two labels, the genus of the component $\Sigma_\alpha$
denoted by $g_\alpha$ and the set of indices of marked points
which lie in the component $\Sigma_\alpha$.
\nfd

\rmk\label{rmk:signature}
Two marked nodal surfaces are isomorphic if and only if
they have the same signature.
\kmr

\begin{proof}
In other words,  $(\Sigma,r_*,\nu)$ and  $(\Sigma',r'_*,\nu')$
have the same signature if and only if
there is a diffeomorphism $\phi:\Sigma\to\Sigma'$
such that $\nu'=\phi_*\nu$ where
$$
\phi_*\nu:=\bigl\{\{\phi(y_1),\phi(y_2)\},\{\phi(y_3),\phi(y_4)\},
\ldots,\{\phi(y_{2k-1}),\phi(y_{2k})\}\bigr\}
$$
and $r'_i=\phi(r_i)$ for $i=1,2,\ldots,n=n'$.  This
is because two connected surfaces are diffeomorphic
if and only if they have the same genus and any bijection
between two finite subsets of a connected surface
extends to a diffeomorphism of the ambient manifold.
\end{proof}

\dfn\label{def:arithGenus}
Define the \jdef{Betti numbers} of a graph by the formula
$$
b_i:=\mathrm{rank}\, H_i(K), \qquad i=0,1,
$$
where $H_i(K)$ is the $i$th homology group of the cell complex $K$.
Thus $K$ is connected if and only if $b_0=1$ and
$$
b_0-b_1 = \; \mbox{$\#$ vertices - $\#$ edges}.
$$
Define the \jdef{genus} of the labelled graph by
$$
g:= b_1+ \sum_\alpha g_\alpha.
$$
The \jdef{arithmetic genus} of a nodal surface $(\Sigma,\nu)$
is the genus of the signature of $(\Sigma,\nu)$.
Note that the arithmetic genus can be different
from the \jdef{total genus} $g':=\sum_\alpha g_\alpha$.
\nfd

\dfn\label{def:stableRiemannSurface}
A marked nodal surface $(\Sigma,r_*,\nu)$ said to be of \jdef{type $(g,n)$}
iff the length of the sequence $r_*$ is $n$, the underlying graph $K$ in the
signature is connected, and the arithmetic genus of $(\Sigma,\nu)$ is $g$.
A marked nodal Riemann surface $(\Sigma,r_*,\nu, j)$  is called \jdef{stable}
iff its \jdef{automorphism group}
$$
\Aut(\Sigma,r_*,\nu.j):=
\{\phi\in\Diff(\Sigma): \phi_*j=j,\; \phi_*\nu=\nu,\; \phi(r_*)=r_*\}
$$
is finite. A stable marked nodal Riemann surface is commonly
called a \jdef{stable curve}.
\nfd

\para\label{stableCondition}
A marked nodal Riemann surface of type $(g,n)$  is stable if and only if
the number of special points in each
component of genus zero is at least three and
the number of special points in each component of genus
one is at least one.
This  is an immediate consequence  of the following:
\begin{description}
\item[(i)] An automorphism of a surface of
genus zero is a  M\"obius transformation;
if it  fixes three points it is the identity.
\item[(ii)] A surface of genus one is isomorphic to $\C/\Lambda$
where $\Lambda=\Z\oplus\Z\tau$ and $\tau$ lies in the upper half plane.
\item[(iii)] The automorphisms of the
Abelian group $\Lambda$ of form $z\mapsto az$ where
$a\in\C\setminus0$ form a group of order at most six.
\item[(iv)] The automorphism group of a compact Riemann surface
of genus greater than one is finite.
\end{description}
The proofs of these well known assertions
can be found in any book on Riemann surfaces.
It follows that
for each pair $(g,n)$ of nonnegative integers
there are only finitely many  labelled
graphs which arise as the signature of a
stable marked nodal Riemann surface of  type $(g,n)$.
\arap

\rmk
A  marked nodal surface  has arithmetic genus zero
if and only if each component has genus zero and the graph is a tree.
The automorphism group of a stable  marked nodal Riemann surface
of arithmetic genus zero is trivial.
\kmr


\section{Nodal families}\label{sec:nodal-fam}

In this section we introduce the basic setup which will allow us
to define the charts of the Deligne--Mumford orbifold. 

\para\label{coordinates}
Let  $P$ and $A$ be complex manifolds with  $\dim_\C(P)=\dim_\C(A)+1$
and
$$
\pi:P\to A
$$
be a holomorphic map.  By the holomorphic implicit function theorem
a point $p\in P$ is a regular point of $\pi:P\to A$
if and only if there is a holomorphic coordinate system
$(t_1,\ldots,t_n)$ defined in a neighborhood of
$\pi(p)\in A$, and a function $z$ defined in a
neighborhood of $p$ in $P$ such that
$(z,t_1\circ\pi,\ldots,t_n\circ\pi)$ is
holomorphic coordinate system.
In other words, the point $p$ is a regular point if and only if
the germ of $\pi$ at $p$ is isomorphic to the  germ at $0$ of
the projection
$$
    \C^{n+1}\to\C^n:(z,t_1,\ldots,t_n)\mapsto (t_1,\ldots,t_n).
$$
Similarly,  a point $p\in P$ is a called a \jdef{nodal point} of $\pi$
if and only if the germ of $\pi$ at $p$
is isomorphic to the  germ at $0$ of the map
$$
    \C^{n+1}\to\C^n:(x,y,t_2,\ldots,t_n)\mapsto (xy,t_2,\ldots,t_n),
$$
i.e. if and only if there are holomorphic coordinates
$z,t_2,\ldots,t_n$ on $A$  at $\pi(p)$ and holomorphic functions
$x$ and $y$ defined in a neighborhood of $p$ such that
 $(x,\,y,\,t_2\circ\pi,\,\ldots,\,t_n\circ\pi)$
is a holomorphic coordinate system, $x(p)=y(p)=0$,
and $xy=z\circ\pi$.
At a regular point $p$ we have that
 $\dim_\C\ker (d\pi(p))=1$ and $\dim_\C\coker (d\pi(p))=0$
while at a nodal point we have that
 $\dim_\C\ker (d\pi(p))=2$ and $\dim_\C\coker (d\pi(p))=1$
\arap

\dfn
A \jdef{nodal family} is a surjective proper holomorphic map
$$
\pi:P\to A
$$
between connected complex manifolds such that
$\dim_\C(P)=\dim_\C(A)+1$ and every critical point of $\pi$ is nodal.
We denote the set of critical points of $\pi$ by
$$
    C_\pi:=\{p\in P: d\pi(p) \mbox{ not surjective}\}.
$$
It intersects each fiber
$$
   P_a:=\pi^{-1}(a)
$$
in a finite set. For each regular value $a\in A$ of $\pi$ the fiber $P_a$
is a compact Riemann surface. When $a\in A$ is a critical value of $\pi$
we view the fiber $P_a$ as a nodal Riemann surface as follows.

By the maximum principle the composition $\pi\circ u$
of $\pi$ with a holomorphic map $u:\Sigma\to P$ defined on
a compact Riemann surface $\Sigma$ must be constant,
i.e. $u(\Sigma)\subset P_a$ for some $a$.
A  \jdef{desingularization} of a fiber $P_a$ is a  holomorphic map
$u:\Sigma\to P$
defined  on a compact Riemann surface $\Sigma$ such that
\begin{description}
\item[(1)] $u^{-1}(C_\pi)$ is finite,
\item[(2)] the restriction
of $u$ to $\Sigma\setminus u^{-1}(C_\pi)$ maps this  set
bijectively  to $P_a\setminus C_\pi$.
\end{description}
The restriction of $u$ to $\Sigma\setminus u^{-1}(C_\pi)$
is an  isomorphism between this open Riemann surface and  $P_a\setminus C_\pi$
(because it is holomorphic, bijective, and proper).
\nfd

\begin{lemma} \label{le:desing}
\begin{description}
\item[(i)] Every fiber  of a nodal family admits a desingularization.
\item[(ii)] If $u_1:\Sigma_1\to P$ and $u_2:\Sigma_2\to P$
are two desingularizations of the same fiber, then the map
$$
    u_2^{-1}\circ u_1:
\Sigma_1\setminus u_1^{-1}(C_\pi)\to\Sigma_2\setminus u_2^{-1}(C_\pi)
$$
extends to an isomorphism $\Sigma_1\to\Sigma_2$.
\item[(iii)] A desingularization $u$ of a
fiber of a nodal family is an immersion
and the preimage  $u^{-1}(p)$ of a critical point $p\in C_\pi$ consists of
exactly two points.
\end{description}
\end{lemma}

\begin{proof}
Let $\pi:P\to A$ be a nodal family and $a\in A$.
Each $p\in C_\pi\cap P_a$ has a small neighborhood
intersecting $P_a$ in two transverse embedded holomorphic
disks intersecting at $p$.
Define $\Sigma$ set theoretically as
the disjoint union of $P_a\setminus C_\pi$
with two copies of $P_a\cap C_\pi$ and use these disks as
coordinates; the map $u:\Sigma\to P_a$ is the identity on $P_a\setminus C_\pi$
and sends each pair of nodal points to the point of $C_p$ which gave
rise to it.  Assertion~(ii) follows from
the removable singularity theorem for holomorphic functions
and~(iii) follows from~(ii) and the fact that the maps
$x\mapsto (x,0)$ and $y\mapsto (0,y)$ are immersions.
\end{proof}

\rmk\label{rmk:canonical-desingularization}
We can construct a \jdef{canonical desingularization} of the fiber
by replacing each  point  $p\in P_a\cap C_\pi$ by a point
for each   connected component of   $U\setminus\{p\}$
where $U$ is a suitable neighborhood of $p$ in $P_a$
and extending the smooth and complex structures in the only way possible.
\kmr

\dfn\label{def:fiberIso}
Let $\pi_A:P\to A$ and $\pi_B:Q\to B$ be nodal families.
For
$a\in A$ and $b\in B$ a bijection $f:P_a\to Q_b$
is called a  \jdef{fiber isomorphism} if for some (and hence every)
desingularization $u:\Sigma\to P_a$ the map
$f\circ u:\Sigma\to Q_b$ is a desingularization.
A \jdef{pseudomorphism} from  $\pi_A$ to $\pi_B$  is a commutative diagram
$$
 \Rectangle{P}{\Phi}{Q}{\pi_A}{\pi_B}{A}{\phi}{B}
$$
where $\Phi$ and $\phi$ are smooth and, for each $a\in A$,
the restriction of $\Phi$ to the fiber $P_a$ is  a fiber isomorphism.
A \jdef{morphism} is a pseudo morphism
such that both $\phi$ and $\Phi$ are holomorphic.
For $a\in A$ and $b\in B$ the  notation
$$
   (\Phi,\phi):(\pi_A,a)\to(\pi_B,b)
$$
indicates that the pseudo morphism $(\Phi,\phi)$ satisfies $\phi(a)=b$.
\nfd

\begin{lemma}\label{le:arithGenus}
Let $\pi:P\to A$ be a  nodal family.
Then the arithmetic genus
(see Definition~\ref{def:arithGenus})
of the fiber $P_a$
is a locally constant function of $a\in A$.
\end{lemma}

\begin{proof}
The arithmetic genus is the genus of the surface obtained
by removing a small disk about each nodal point and
identifying corresponding components.
Hence it is equal to the ordinary genus of a regular fiber.
\end{proof}

\dfn
A \jdef{marked  nodal family}  is a pair $(\pi,R_*)$
where $\pi:P\to A$ is a nodal family and
$$
  R_*=(R_1,\ldots,R_n)
$$
is a sequence of complex submanifolds of $P$ which are pairwise disjoint
and such that $\pi|R_i$ maps $R_i$ diffeomorphically onto $A$. It follows
that $R_i$ does not intersect the set $C_\pi$ of critical points.
A desingularization  $u:\Sigma\to P$  of a fiber $P_a$
of a marked nodal  family $(\pi,R_*)$ determines  a point marking $r_*$:
the  point marking $r_*$ is given by the formula
$$
       \{u(r_i)\}= R_i\cap P_a
$$
for $i=1,2,\ldots,n$.
By Lemma~\ref{le:desing} any two desingularizations
of the same fiber give rise to isomorphic marked nodal Riemann surfaces.
Thus the signature (see Definition~\ref{def:signature})
of the fiber $(P_a,P_a\cap R_*)$
is independent of the choice of the desingularization used to define it.
In the context  of marked  nodal families, the term {\em fiber isomorphism}
is understood to entail that the bijection
$f$ preserves the induced point markings;
similarly pseudo morphisms and morphisms of marked  nodal families preserve
the corresponding point markings.
We say that the marked nodal family $(\pi,R_*)$ is of \jdef{type $(g,n)$}
when each fiber is of type $(g,n)$
(see Definition~\ref{def:stableRiemannSurface}).
\nfd

\dfn
A fiber of a marked nodal family $\pi:P\to A$ is called \jdef{stable}
iff its desingularization is stable.
A marked nodal family is called \jdef{stable}
iff each of its fibers is stable.
\nfd

\rmk
It is easy to see that stability is an open condition,
i.e. every stable fiber has a neighborhood
consisting of stable fibers.
However, the open set of stable fibers can have unstable
fibers in its closure. For example, consider
the nodal  family $(\pi,(R_1,R_2,R_3))$
with
$$
   P =\{([x,y,z],a)\in \C P^2\times\C: xy=az^2\},
$$
$A=\C$,
$\pi([x,y,z],a)=a$,
$R_1=\{[1,0,0]\}\times A$,
$R_2=\{[0,1,0]\}\times A$, and
$R_3=\{([1,a,1],a):a\in A\}$.
The
desingularization of the fiber  over $0$
consists of two components of genus zero
and the regular fibers consist of one component
of genus zero. The regular fibers all have
three marked points and are thus stable;
one of the two components of the (desingularized)
singular fiber has fewer than three special points
and is thus unstable.
\kmr


\section{Universal unfoldings}\label{sec:universal}

In this section we formulate the most important definitions
and theorems of this paper.  The key definition is that of a 
universal unfolding.  Once we have established the existence 
of universal unfoldings, the definition of the orbifold structure
on the Deligne--Mumford moduli space (which we carry out
in  the next section) becomes almost tautological. 
The most important theorem asserts that an unfolding 
is universal if and only if it satisfies a suitable infinitesimal
conditon (which is easier to verify). 

\dfn
A \jdef{nodal unfolding} is a triple $(\pi_B, S_*,b)$
consisting of a marked
nodal family $(\pi_B:Q\to B,S_* )$ and a point $b\in B$ of the base $B$.
The fiber $Q_b$ is called the \jdef{central fiber} of the unfolding
and the unfolding is said to be an unfolding of the marked nodal
Riemann surface induced by any desingularization of this central fiber.
The unfolding is called \jdef{universal} iff
for every other nodal unfolding $(\pi_A:P\to A,R_*,a)$ and
any fiber isomorphism $f:P_a\to Q_b$
there is a unique germ of a morphism
$$
(\Phi,\phi): (\pi_A,a)\to (\pi_B,b)
$$
such that $\Phi(R_i)\subset S_i$ for all $i$  and $\Phi|P_a=f$.
The term  {\em germ} means that $\phi$ is defined in a neighborhood
of $a$ in $A$ and $\Phi$ is defined on the preimage of this neighborhood
under $\pi_A$. The term {\em unique} means that if $(\Phi',\phi')$
is another morphism with the same properties then
it agrees with $(\Phi,\phi)$ over a sufficiently small neighborhood of $a$.
\nfd

\dfn \label{def:D}
Let $(\pi:Q\to B, S_*,b)$ be an unfolding of a marked nodal Riemann surface
$(\Sigma,s_*,\nu,j)$ and $u:\Sigma\to Q_b$ be a desingularization.
Let $\cX_{u,b}$ denote the space
$$
\cX_{u,b}:=\left\{(\hat u,\hat b)\in\Omega^0(\Sigma,u^*TQ)\times T_bB\biggl|
\begin{array}{l}   d\pi(u)\hat u=\hat b,\quad
     \hat u(s_i)\in T_{u(s_i)}S_i,\mbox{ and }\\
     u(z_1)=u(z_2)\implies \hat u(z_1)=\hat u(z_2).
\end{array}\right\}
$$
Let $\cY_u$ denote the space
$$
     \cY_u := \left\{\eta\in\Omega^{0,1}(\Sigma,u^*TQ):
     d\pi(u)\eta=0\right\}.
$$
For  $(\hat u,\hat b)\in\cX_{u,b}$ define
$$
D_{u,b}(\hat u,\hat b):= D_u\hat u
$$
where $D_u:\Omega^0(\Sigma,u^*TQ)\to\Omega^{0,1}(\Sigma,u^*TQ)$
is the linearized Cauchy--Riemann operator.
We call the unfolding $(\pi,S_*,b)$ \jdef{infinitesimally universal}
if the operator $D_{u,b}:\cX_{u,b}\to\cY_u$ is bijective
for some (and hence every) desingularization of the central fiber.
Theorems~\ref{thm:stable}, \ref{thm:pseudo},
and~\ref{thm:existence} which follow are proved in
Section~\ref{sec:proof} below.
\nfd

\begin{theorem} [Stability]\label{thm:stable}
Let $(\pi,S_*,b_0)$ be an infinitesimally universal unfolding.
Then $(\pi,S_*,b)$ is infinitesimally universal for $b$
sufficiently near $b_0$.
\end{theorem}

\begin{theorem}[Universal Unfolding]\label{thm:idm}
An unfolding $(\pi,S_*,b)$  is universal
if and only if it is infinitesimally universal.
\end{theorem}

\begin{proof}
We prove `if' in Section~\ref{sec:proof}.
For `only if' we argue as follows.
A composition of morphisms (of nodal unfoldings) is again a morphism.
The only morphism which is the identity on the central
fiber of a universal unfolding  is the identity.
It follows that any two universal unfoldings of the same marked nodal
Riemann surface are isomorphic.
By Theorem~\ref{thm:existence} below there is
an infinitesimally universal unfolding
and by `if'  it  is universal and hence
isomorphic to every other universal unfolding.
Any unfolding isomorphic to an infinitesimally
universal unfolding is itself infinitesimally universal.
\end{proof}

\begin{theorem}[Uniqueness]\label{thm:pseudo}
Let $(\pi_B,S_*,b_0)$ be an infinitesimally universal unfolding.
Then every pseudomorphism from $(\pi_A,R_*,a_0)$ to  $(\pi_B,S_*,b_0)$
is a morphism.
\end{theorem}

\begin{theorem}[Existence]\label{thm:existence}
A marked nodal Riemann surface admits an infinitesimally
universal unfolding if and only if it is stable.
\end{theorem}

\begin{proof}
We prove `if' in Section~\ref{sec:proof}.
For `only if' we argue as follows.
Let $(\Sigma,s_*,\nu,j)$ be  a marked nodal Riemann
surface. Assume it is not stable. Then either $\Sigma$ has genus one
and has no special points or else $\Sigma$ contains a component of
genus zero with at most two special points. In either case there
is an abelian complex Lie group $A$ (namely $A=\Sigma$ in the former case
and $A=\C^*$ in the latter) and an effective holomorphic action
$A\times\Sigma\to\Sigma:(a,z)\mapsto a_\Sigma(z)$.
Let $P:=A\times\Sigma$ and $\pi_A$ be the projection
on the first factor. If $v:\Sigma\to Q$ is any desingularization
of a fiber $Q_b$ of an unfolding $\pi_B:Q\to B$, then
$\Phi_1(a,z):=v(z)$ and $\Phi_2(a,z):=v(a_\Sigma(z))$
are distinct morphisms which extend the fiber isomorphism
$(e,z)\mapsto v(z)$.  Hence $\pi_B$ is not universal.
\end{proof}


\section{Universal families and the Deligne--Mumford moduli space}
\label{sec:dm-orbifold}

In this section we define the orbifold structure on the Deligne--Mumford
moduli space.  The proof of compactness will be relegated to 
Section~\ref{sec:compact}. The results we prove in this section
are easy consequences of Theorems~\ref{thm:stable} 
and~\ref{thm:existence}.

\para \label{cB}
Throughout this section $g$ and $n$ are nonnegative integers
with $n>2-2g$.  Let $\bar\cB_{g,n}$ denote the groupoid whose
objects are stable marked nodal Riemann
surfaces of type $(g,n)$ and whose morphisms
are isomorphisms of marked nodal Riemann surfaces.
The \jdef{Deligne--Mumford moduli space}
is the orbit space $\bar\cM_{g,n}$ of this groupoid:
a point of $\bar\cM_{g,n}$  is an equivalence class\footnote
{
Strictly speaking, the equivalence class is a proper class in the sense
of  set theory as explained in the appendix of~\cite{KELLEY}
for example. One could avoid this problem by
choosing for each stable signature (see Remark~\ref{rmk:signature}
and \ref{stableCondition})
a ``standard marked nodal surface'' with that signature
and restricting the space of objects
of the groupoid  $\bar\cB_{g,n}$ to those having
a standard surface  as substrate.
}
of objects of
$\bar\cB_{g,n}$ where two objects are equivalent if and only if
they are isomorphic.
We will introduce a canonical orbifold structure
(see Definition~\ref{def:orbifoldStructure}) on this groupoid.
The following definition is crucial.
\arap

\dfn \label{def:universal}
 A \jdef{universal marked nodal family} of type $(g,n)$
is a marked nodal family
$(\pi_B:Q\to B,S_*)$  satisfying the following conditions.
\begin{description}
\item[(1)] $(\pi_B,S_*,b)$ is a universal unfolding for every $b\in B$.
\item[(2)] Every stable marked nodal Riemann surface of type $(g,n)$
is the domain of a desingularization of at least one fiber of $\pi_B$.
\item[(3)] $B$ is second countable (but  possibly disconnected).
\end{description}

\nfd

\begin{proposition}\label{prop:universal}
For every pair $(g,n)$ with $n>2-2g$ there is a universal marked nodal family.
\end{proposition}

\begin{proof}
By Theorems~\ref{thm:existence}, \ref{thm:idm}, and~\ref{thm:stable},
each stable marked nodal Riemann surface admits a universal
unfolding satisfying~(1) and~(3). To construct a universal unfolding
that also satisfies~(2) we must cover $\bar\cM_{g,n}$ by countably many
such families. This is possible because $\bar\cM_{g,n}$ is a union of
finitely many strata, one for each stable signature, and each stratum
is a separable topological space.
\end{proof}

\dfn\label{B-Gamma}
Let $(\pi_B:Q\to B,S_*)$ be a universal marked nodal family.
The \jdef{associated groupoid} is the tuple $(B,\Gamma,s,t,e,i,m)$, 
where $\Gamma$ denotes the set of all triples $(a,f,b)$ 
such that $a,b\in B$ and $f:Q_a\to Q_b$ is a fiber isomorphism,
and the structure maps $s,t:\Gamma\to B$, $e:B\to\Gamma$,
$i:\Gamma\to\Gamma$, and 
$m:\Gamma\Times{s}{t} \Gamma\to\Gamma$ 
are defined by
$$
s(a,f,b):=a,\qquad
t(a,f,b):=b,\qquad
e(a):=(a,\id,a),
$$
$$
i(a,f,b):=(b,f^{-1},a),\qquad
m((b,g,c),(a,f,b)):=(a,g\circ f,c).
$$
The associated groupoid is equipped with a functor
$$
B\to\bar\cB_{g,n}:b\mapsto\Sigma_b
$$
to the groupoid $\bar\cB_{g,n}$ of~\ref{cB}. 
In other words, $\iota_b:\Sigma_b\to Q_b$
denotes the canonical desingularization defined in
Remark~\ref{rmk:canonical-desingularization}.
By definition the induced map
$$
B/\Gamma\to\bar\cM_{g,n}:[b]_B\mapsto [\Sigma_b]_{\bar\cB_{g,n}},\qquad
[b]_B:=\{t(f):f\in\Gamma,\;s(f)=b\},
$$
on orbit spaces is bijective.  The next theorem asserts that 
the groupoid $(B,\Gamma)$ equips the moduli space 
$\bar\cM_{g,n}$ with an orbifold structure which is independent of the 
choice of the universal family.  This is the \jdef{orbifold structure}
on the Deligne--Mumford moduli space. 
\nfd

\begin{theorem}\label{thm:Gamma}
{\bf (i)}
Let  $(\pi_B:Q\to B,S_*)$ be universal as in Definition~\ref{def:universal}
and $(B,\Gamma)$ be the associated groupoid of Definition~\ref{B-Gamma}.
Then there is a unique complex manifold structure on $\Gamma$
such that $(B,\Gamma)$ is a complex etale Lie groupoid
with structure maps $s,t,e,i,m$.

\smallskip\noindent{\bf (ii)}
A morphism between universal families $\pi_0:Q_0\to B_0$ and
${\pi_1:Q_1\to B_1}$ induces a refinement
$\iota:(B_0,\Gamma_0)\to(B_1,\Gamma_1)$
of the associated etale groupoids.

\smallskip\noindent{\bf (iii)}
The orbifold structure on $\bar\cM_{g,n}$ introduced in
Definition~\ref{B-Gamma} is independent of the choice 
of the universal marked nodal family $(\pi_B,S_*)$
used to define it. 
\end{theorem}

\begin{proof}
We prove~(i). Uniqueness is immediate since part of the definition of
complex etale Lie groupoid is that $s$ is a local holomorphic
diffeomorphism.  We prove existence.
It follows from the definition of universal unfolding that each
triple $(a_0,f_0,b_0)\in\Gamma$ determines a morphism
$$
\Rectangle{Q|U}{\Phi}{Q|V}{\pi}{\pi}{U}{\phi}{V}
$$
for suitable neighborhoods $U\subset B$ of $a_0$ and $V\subset B$
of $b_0$ such that $\Phi|Q_{a_0}=f_0$. Every such morphism determines
a chart
$
\iota_\Phi:U\to\Gamma
$
given by
$$
\iota_\Phi(a) := (a,\Phi_a,\phi(a)).
$$
(In this context a chart is a bijection between an open set
in a complex manifold and a subset of $\Gamma$.)
By construction each transition map between
two such charts is the identity. This defines the
manifold structure on $\Gamma$.
That the structure maps are holomorphic follows from the identities
$$
s\circ\iota_\Phi=\id,\qquad t\circ\iota_\Phi=\phi,\qquad
e=\iota_\id,
$$
$$
i\circ\iota_\Phi=\iota_{\Phi^{-1}}\circ\phi,\qquad
m\circ\left(\iota_\Psi\circ\phi\times\iota_\Phi\right) = \iota_{\Psi\circ\Phi}.
$$
This proves~(i).  

We prove~(ii).  If $(\phi,\Phi)$ is a morphism from $\pi_0$ to $\pi_1$
then the refinement $\iota:(B_0,\Gamma_0)\to(B_1,\Gamma_1)$ of~(ii)
is given by
$$
(a_0,f_0,b_0) \mapsto
(\phi(a_0),\Phi_{b_0}\circ f_0\circ\Phi_{a_0}^{-1},\phi(b_0)).
$$
This proves~(ii).

We prove~(iii).  
Let $\pi_0:Q_0\to B_0$ and $\pi_1:Q_1\to B_1$ be universal families.
For each $b\in B_0$ choose a neighborhood $U_b\subset B_0$ of $b$
and a morphism $\Phi_b:Q_0|U_b\to Q_1$. Cover $B_0$ by countably many
such neighborhoods $U_{b_i}$.  Then the disjoint union $B$ of the nodal
families $Q_0|U_{b_i}$ defines another universal family $\pi:Q\to B$
equiped with morphisms to both $\pi_0$ and $\pi_1$ (to $\pi_0$ by inclusion
and to $\pi_1$ by construction).
Now each morphism of universal families induces
a refinement of the corresponding orbifold structures.
\end{proof}

\begin{theorem}\label{thm:proper}
Let $(\pi_B:Q\to B,S_*)$ be a universal family.
Then the etale groupoid $(B,\Gamma)$ constructed in Definition~\ref{B-Gamma}
is proper and the quotient topology on $B/\Gamma$ is compact.
\end{theorem}

\begin{proof}
See Section~\ref{sec:compact} below.
\end{proof}

\xmpl\label{ex:GK}
Assume $g=0$.  Then the moduli space $\bar\cM_{0,n}$ of marked
nodal Riemann surfaces of genus zero (called the
\jdef{Grothendieck--Knudsen compactification}) is a compact
connected complex manifold (Knudsen's theorem).
In our formulation this follows from the fact that the automorphism
group of each marked nodal Riemann surface of genus zero
consists only of the identity. In~\cite[Appendix~D]{MS2} the
complex manifold structure on $\bar\cM_{0,n}$ is obtained
from an embedding into a product of 2-spheres via cross ratios.
That the manifold structure in~\cite{MS2} agrees with ours
follows from the fact that the projection
$\pi:\bar\cM_{0,n+1}\to\bar\cM_{0,n}$ (with the complex
manifold structures of~\cite{MS2}) is a universal family
as in Definition~\ref{def:universal}.
\lpmx


\section{Complex structures on the sphere}\label{sec:S2}

In preparation for the construction of universal unfoldings 
(without nodes and marked points) we review the
space of complex structures on a Riemann surface $\Sigma$
in this and the following two sections.  This section treats the case 
of genus zero.  Denote by $\cJ(S^2)$ the space of complex structures
on $S^2$ that induce the standard orientation and by
$\Diff_0(S^2)$ the group of orientation preserving
diffeomorphisms of $S^2$.

\begin{theorem}\label{thm:S2}
There is a fibration
$$
\begin{array}{ccc}
     \mathrm{PSL}_2(\C) & \to & \Diff_0(S^2) \\
     & & \downarrow \\
     & & \cJ(S^2)
\end{array}
$$
where the inclusion $\mathrm{PSL}_2(\C)\to\Diff_0(S^2)$
is the action by M\"obius
transformations and the projection
$\Diff_0(S^2)\to\cJ(S^2)$ sends $\phi$ to $\phi^*i$.
\end{theorem}

The theorem asserts that the map
$\Diff_0(S^2)\to\cJ(S^2)$ has the path lifting property
for smooth paths and that the lifting depends smoothly on the path.
One consequence of this, as observed
in~\cite{EE}, is the celebrated theorem
of Smale~\cite{SMALE} which asserts that
$\Diff_0(S^2)$ retracts onto $\mathrm{SO}(3)$.
Another consequence is that, up to diffeomorphism, there is a unique complex
structure on the $2$-sphere. Yet another consequence is that
a proper holomorphic submersion whose fibers have genus zero
is holomorphically locally trivial. (See Theorem~\ref{thm:teich-mark}.)

\begin{proof}[Proof of Theorem~\ref{thm:S2}.]
Choose a smooth path
$
      [0,1]\to\cJ(S^2):t\mapsto j_t.
$
We will find an isotopy $t\mapsto\psi_t$
of $S^2$ such that
\begin{equation}\label{eq:psij}
      {\psi_t}^*j_t = j_0.
\end{equation}
Suppose that the unknown isotopy $\psi_t$ is generated
by a smooth family of vector fields $\xi_t\in\Vect(S^2)$
via
$$
     \frac{d}{dt}\psi_t=\xi_t\circ\psi_t,\qquad
     \psi_0=\id.
$$
Then~(\ref{eq:psij}) is equivalent to 
$\psi_t^*( \cL_{\xi_t}j_t + \hat j_t)=0$ 
and hence to
\begin{equation}\label{eq:xij}
     \cL_{\xi_t}j_t + \hat j_t = 0,
\end{equation}
where $\hat j_t:=\frac{d}{dt} j_t\in\Cinf(\End(TS^2))$.
As usual we can think of $\hat j_t$ as a $(0,1)$-form on $S^2$
with values in the complex line bundle
$$
      E_t := (TS^2,j_t).
$$
The vector field $\xi_t$ is a section of this line bundle.
This line bundle is holomorphic
and its  Cauchy-Riemann operator
$$
      \bar\p_{j_t}:\Cinf(E_t)\to\Omega^{0,1}(E_t)
$$
has the form
$$
      \bar\p_{j_t}\eta
      = \frac{1}{2}\left(\nabla \eta + j_t\circ\nabla \eta\circ j_t\right)
$$
where $\nabla$ is the Levi-Civita connection of
the Riemannian metric  $\omega(\cdot,j_t\cdot)$ on $S^2$
and $\omega\in\Omega^2(S^2)$ denotes the standard volume form.
Now, for every vector field $\eta\in\Vect(S^2)$, we have
\begin{eqnarray*}
      (\cL_{\xi_t}j_t)\eta
&= &
      \cL_{\xi_t}(j_t\eta) - j_t\cL_{\xi_t}\eta  \\
&= &
      [j_t\eta,\xi_t] - j_t[\eta,\xi_t]  \\
&= &
      \Nabla{\xi_t}(j_t\eta) - \Nabla{j_t\eta}\xi_t
      - j_t\Nabla{\xi_t}\eta + j_t\Nabla{\eta}\xi_t \\
&= &
      j_t\Nabla{\eta}\xi_t - \Nabla{j_t\eta}\xi_t \\
&= &
      2j_t(\bar\p_{j_t}\xi_t)(\eta).
\end{eqnarray*}
The penultimate equality uses the fact that
$j_t$ is integrable and so $\nabla j_t=0$.
Hence equation~(\ref{eq:xij}) can be expressed
in the form
\begin{equation}\label{eq:xij'}
     \bar\p_{j_t}\xi_t
     =  \frac{1}{2}j_t\hat j_t. 
\end{equation}
Now the line bundle $E_t$ has Chern number
$c_1(E_t)=2$ and hence, by the Riemann-Roch theorem,
the Cauchy-Riemann operator $\bar\p_{j_t}$ has real
Fredholm index six and is surjective for every $t$.
Denote by
$$
      {\bar\p_{j_t}}^*
      :\Omega^{0,1}(E_t)\to\Cinf(E_t)
$$
the formal $L^2$-adjoint operator of $\bar\p_{j_t}$.
By elliptic regularity, the formula
$$
      \xi_t
      :=  \frac{1}{2}{\bar\p_{\xi_t}}^*
        \left(\bar\p_{j_t}{\bar\p_{j_t}}^*\right)^{-1}
        \left(j_t\hat j_t\right) 
$$
defines a smooth family of vector fields on $S^2$
and this family obviously satisfies~(\ref{eq:xij'}).
Hence the isotopy $\psi_t$ generated by $\xi_t$
satisfies~(\ref{eq:psij}).
\end{proof}

\begin{lemma}\label{le:markedPT}
Let $\C\to\cJ(S^2):s+it\mapsto j_{s,t}$ be holomorphic
and $\C\to\Diff(S^2):s+it\mapsto\phi_{s,t}$ be the unique family
of diffeomorphisms satisfying
$$
\phi_{s,t}^*j_{s,t}=i, \qquad
\phi_{s,t}(0)=0,\qquad \phi_{s,t}(1)=1,\qquad \phi_{s,t}(\infty)=\infty.
$$
Then the map
$$
\C\times S^2\to\C\times S^2 :(s+it,z)\mapsto (s+it,\phi_{s,t}(z))
$$
is holomorphic with respect to the standard complex structure
at the source and the complex structure
$$
J(s,t,z):=\left(\begin{array}{cc} i & 0\\ 0 & j_{s,t}(z)\end{array}\right)
$$
at the target.
\end{lemma}

\begin{proof} Define $\xi_{s,t},\eta_{s,t}\in\Vect(S^2)$ by
$$
 \p_s\phi_{s,t}=\xi_{s,t}\circ\phi_{s,t}, \qquad
 \p_y\phi_{s,t}=\eta_{s,t}\circ\phi_{s,t}.
$$
Differentiating the identity $\phi_{s,t}^*j_{s,t}=i$ gives
$\p_sj+\cL_\xi j=\p_tj+\cL_\eta j=0$.   Since $s+it\mapsto j_{s,t}$ is holomorphic
we have
$$
 0=\p_sj+j\p_tj= -\cL_\xi j-j\cL_\eta j=  -\cL_{\xi+j\eta}j
$$
where the last equality uses the integrability of $j$. Thus $\xi_{s,t}+j_{s,t}\eta_{s,t}$
is a holomorphic vector field vanishing at three points so $\xi_{s,t}+j_{s,t}\eta_{s,t}=0$
for all $s,t$. Hence by definition of $\xi$ and $\eta$ we have
$$
\p_s\phi+j\p_t\phi=0
$$
as required.
\end{proof}


\section{Complex structures on the torus}\label{sec:T2}

Continuing the preparatory discussion of the previous 
section we treat the case of genus one. 
Denote by $\cJ(\T^2)$ the space of complex structures
on the $2$-torus $\T^2:=\R^2/\Z^2$ that induce the standard
orientation and by $\Diff_0(\T^2)$ the group of diffeomorphisms of
$\T^2$ that induce the identity on homology. Denote the elements
of the upper half plane $\H$ by $\lambda=\lambda_1+i\lambda_2$
and consider the map $j:\H\to\cJ(\T^2)$, given by
\begin{equation}\label{eq:jlambda}
j(\lambda) := \frac{1}{\lambda_2}\left(\begin{array}{cc}
-\lambda_1 & -\lambda_1^2-\lambda_2^2 \\
1 & \lambda_1
\end{array}\right).
\end{equation}
Thus $j(\lambda)$ is the pullback of the standard
complex structure under the diffeomorphism
$$
f_\lambda:\T^2\to\frac{\C}{\Z+\lambda\Z},\qquad
f_\lambda(x,y):=x+\lambda y.
$$
A straight forward calculation shows that the map
$j:\H\to\cJ(\T^2)$ is holomorphic as is the map
$$
(\lambda,z+\Z+\lambda\Z)\mapsto(j(\lambda),f_\lambda^{-1}(z)+\Z^2)
$$
from $\left\{(\lambda,z+\Z+\lambda\Z)\,:\,\lambda\in\H,\,z\in\C\right\}$
to $\cJ(\T^2)\times\T^2$. The next theorem shows that the map
$j:\H\to\cJ(\T^2)$ is a global slice for the action of $\Diff_0(\T^2)$.

\begin{theorem}\label{thm:T2}
There is a proper fibration
$$
\begin{array}{ccc}
     \T^2 & \to & \Diff_0(\T^2)\times\H \\
     & & \downarrow \\
     & & \cJ(\T^2)
\end{array}
$$
where the inclusion $\T^2\to\Diff_0(\T^2)$ is the action by translations
and the projection $\Diff_0(\T^2)\times\H\to\cJ(\T^2)$ sends
$(\phi,\lambda)$ to $\phi^*j_\lambda$.
\end{theorem}

The theorem asserts that the map $\Diff_0(\T^2)\times\H\to\cJ(\T^2)$
has the path lifting property for smooth paths and that the lifting
depends smoothly on the path. One consequence of this is that
$\Diff_0(\T^2)$ retracts onto $\T^2$. Another consequence is that
every complex structure on $\T^2$ is diffeomorphic to
$j_\lambda$ for some $\lambda\in\H$.

\begin{proof}[Proof of Theorem~\ref{thm:T2}.]
The uniformization theorem asserts that for every $j\in\cJ(\T^2)$ there is
a unique volume form $\omega_j\in\Omega^2(\T^2)$ with $\int_{\T^2}\omega_j=1$
such that the metric $g_j=\omega_j(\cdot,j\cdot)$ has constant curvature zero.
(A proof can be based on the Kazdan--Warner equation.) Hence it follows from
the Cartan--Ambrose--Hicks theorem that, for every positive real number $\mu$,
there is an orientation preserving diffeomorphism $\psi_j:\C\to\R^2$, unique
up to composition with a rotation, such that
$$
\psi_j^*g_j=\mu g_0,\qquad\psi(0)=0.
$$
Here $g_0$ denotes the standard metric on $\C$.  We can choose $\mu$ and the
rotation such that
$
\psi_j(1)=(1,0).
$
This determines $\psi_j$ (and $\mu$) uniquely. The orientation
preserving condition shows that
$
\lambda_j:=\psi_j(i) \in\H.
$
Moreover, it follows from the invariance of
$g_j$ under the action of $\Z^2$ that
$$
\psi_j(\Z+\lambda\Z) = \Z^2.
$$
Hence $\psi_j$ induces an isometry
of flat tori $(\C/\Z+\lambda_j\Z,g_0)\to(\T^2,g_j)$
which will still be denoted by $\psi_j$.  Let $\phi_j$ be the
precomposition of this isometry with the map
$\T^2\to\C^2/\Z+\lambda_j\Z:(x,y)\mapsto x+\lambda_jy$.
Then $\phi_j\in\Diff_0(\T^2)$ and
$
\phi_j^*j=j(\lambda_j).
$
Thus we have proved that the map
$$
\Diff_{00}(\T^2)\times\H\to\cJ(\T^2):(\phi,\lambda)\mapsto\phi^*j(\lambda)
$$
is a bijection, where $\Diff_{00}(\T^2)$ denotes the subgroup
of all diffeomorphisms $\phi\in\Diff_0(\T^2)$ that satisfy $\phi(0)=0$.
That the map $\Diff_{00}(\T^2)\times\H\to\cJ(\T^2)$ is actually
a diffeomorphism follows by examining the linearized operator at points
$(\phi,\lambda)$ with $\phi=\id$ and noting that it is a bijection
(between suitable Sobolev completions). This proves the theorem.
\end{proof}


\section{Complex structures on surfaces of higher genus}\label{sec:cx}

Continuing the preparatory discussion of the previous 
two sections we treat the case of genus bigger than one. 
Let $\Sigma$ be a compact connected oriented $2$-manifold
of genus $g>1$ and $\cJ(\Sigma)$ be the Frech\'et
manifold of complex structures $j$ on $\Sigma$,
i.e. $j$ is an automorphism of $T\Sigma$ such that $j^2=-\one$.
The identity component $\Diff_0(\Sigma)$ of the
group of orientation preserving diffeomorphisms
acts on $\cJ(\Sigma)$ by $j\mapsto\phi^*j$.
The orbit space
$$
     \cT(\Sigma):=\cJ(\Sigma)/\Diff_0(\Sigma)
$$
is called the \jdef{Teichm\"uller space} of $\Sigma$.
For $j\in\cJ(\Sigma)$ the tangent space $T_j\cJ(\Sigma)$
is the space the space of endomorphisms
$\hat j\in\Omega^0(\Sigma,\End(T\Sigma))$ that
anti-commute with $j$, i.e. $j\hat j+\hat jj=0$.
Thus
$$
     T_j\cJ(\Sigma) = \Omega^{0,1}_j(\Sigma,T\Sigma).
$$
Define an almost complex structure on $\cJ(\Sigma)$
by the formula $\hat j\mapsto j\hat j$.  The next theorem
shows that $\cT(\Sigma)$ is a complex manifold of dimension
$3g-3$.

\begin{theorem}\label{thm:slice}
For every $j_0\in\cJ(\Sigma)$ there exists a holomorphic
local slice through $j_0$. More precisely, there is an open
neighborhood $B$ of zero in $\C^{3g-3}$ and a holomorphic
map $\iota:B\to\cJ(\Sigma)$ such that the map
$$
B\times\Diff_0(\Sigma)\to\cJ(\Sigma):(b,\phi)\mapsto\phi^*\iota(b)
$$
is a diffeomorphism onto a neighborhood of the orbit of $j_0$.
\end{theorem}

\begin{proof}
We first show that each orbit of the action of
$\Diff_0(\Sigma)$ is an almost complex
submanifold of $\cJ(\Sigma)$.
(The complex structure on $\cJ(\Sigma)$
is integrable because $\cJ(\Sigma)$ is the space
of sections of a bundle over $\Sigma$ whose fibers
are complex manifolds.  However, we shall not use this fact.)
The Lie algebra of $\Diff_0(\Sigma)$ is the space of vector
fields
$$
     \Vect(\Sigma)=\Omega^0(\Sigma,T\Sigma).
$$
Its infinitesimal action on $\cJ(\Sigma)$ is
given by
$$
     \Vect(\Sigma)\to T_j\cJ(\Sigma):
     \xi\mapsto \cL_\xi j= 2j\bar\p_j\xi.
$$
Thus the tangent space of the orbit of $j$ is the
image of the Cauchy--Riemann operator
$\bar\p_j:\Omega^0(\Sigma,T\Sigma)\to\Omega^{0,1}_j(\Sigma,T\Sigma)$.
Since $j$ is integrable the operator $\bar\p_j$
is complex linear and so its image is invariant
under multiplication by $j$.

By the Riemann--Roch theorem the operator $\bar\p_j$
has complex Fredholm index $3-3g$.  It is injective because
its kernel is the space of holomorphic sections of a holomorphic
line bundle of negative degree.  Hence its cokernel has dimension
$3g-3$.  Let $B\subset\Omega^{0,1}_{j_0}(\Sigma,T\Sigma)$ be an
open neighborhood of zero in a complex subspace of dimension
$3g-3$ which is a complement of the image of $\bar\p_{j_0}$
and assume that $\one+\eta$ is invertible for every $\eta\in B$.
Define $\iota:B\to\cJ(\Sigma)$ by
$$
\iota(\eta) := (\one+\eta)^{-1}j_0(\one+\eta).
$$
Then
$$
d\iota(\eta)\hat\eta = \left[\iota(\eta),(\one+\eta)^{-1}\hat\eta\right]
$$
and an easy calculation shows that $\iota$ is holomorphic, i.e.
$d\iota(\eta)j_0\hat\eta=\iota(\eta)d\iota(\eta)\hat\eta$
for all $\eta$ and $\hat\eta$.

Let $p>2$ and denote by $\Diff_0^{2,p}(\Sigma)$ and
$\cJ^{1,p}(\Sigma)$ the appropriate Sobolev completions.
Consider the map
$$
     \Diff_0^{2,p}(\Sigma)\times B\to\cJ^{1,p}(\Sigma):
     (\phi,\eta)\mapsto\phi^*\iota(\eta).
$$
This is a smooth map between Banach manifolds and,
by construction, its differential at $(\id,0)$ is bijective.
Hence, by the inverse function theorem,
it restricts to a diffeomorphism from an open neighborhood
of $(\id,0)$ in $\Diff_0^{2,p}(\Sigma)\times B$ to an open
neighborhood of $j_0$ in $\cJ^{1,p}(\Sigma)$.
The restriction of this diffeomorphism to the space of
smooth pairs in $\Diff_0(\Sigma)\times B$ is a
diffeomorphism onto an open neighborhood of $j_0$ in $\cJ(\Sigma)$.
To see this, note that every element of $\Diff_0^{2,p}(\Sigma)$
is a $C^1$-diffeomorphism and that every $C^1$-diffeomorphism
of $\Sigma$ that intertwines two smooth complex structures
is necessarily smooth. Shrink $B$ so that $\{\id\}\times B$
is a subset of the neighborhood just constructed.
The action of $\Diff_0(\Sigma)$ on $\cJ(\Sigma)$
is free and Lemma~\ref{le:proper} below asserts that it is proper.
Hence, by a standard argument, we may shrink $B$
further so that  the local diffeomorphism
$$
     \Diff_0\times B\to\cJ(\Sigma):
     (\phi,\eta)\mapsto\phi^*\iota(\eta)
$$
is injective; it is the required diffeomorphism
onto an open neighborhood of the orbit of $j_0$.
\end{proof}

\begin{lemma}\label{le:proper} Let $\Sigma$ be a surface and
$j_k,j_k'\in\cJ(\Sigma)$ and $\phi_k\in\Diff(\Sigma)$
be sequences such that $j_k'$ converges
to $j'\in\cJ(\Sigma)$ and $j_k=\phi_k^*j_k'$
converges to $j\in\cJ(\Sigma)$.
Then $\phi_k$ has a subsequence which converges
in $\Diff(\Sigma)$.
\end{lemma}

\begin{proof}
Fix an embedded closed disk
$D\subset\Sigma$ and two points $z_0\in\INT(D)$, $z_1\in\p D$.
Let $\D\subset\C$ denote the closed unit disk.
By the Riemann mapping theorem, there is a
unique diffeomorphism $u_k:\D\to D$ such that
$$
u_k^*j_k=i,\qquad u_k(0)=z_0,\qquad u_k(1)=z_1.
$$
The standard bubbling and elliptic  bootstrapping arguments for $J$-holomorphic
curves (see~\cite[Appendix~B]{MS2}) show that $u_k$ converges in the
$\Cinf$-topology.
The same arguments show that the sequence $u_k':=\phi_k\circ u_k$
of $j_k'$-holomorphic disks has a subsequence which converges
on every compact subset of the interior of $\D$.
Thus we have proved that the restriction of $\phi_k$ to any
embedded disk in $\Sigma$ has a convergent subsequence.
Hence $\phi_k$ has a convergent subsequence.
The limit $\phi$ satisfies $\phi^*j'=j$ and has degree one.
Hence $\phi$ is a diffeomorphism.
\end{proof}


\section{The Teichm\"uller space $\cT_g$}\label{sec:teichmuller}

In this section we prove Theorems~\ref{thm:stable}-\ref{thm:existence} 
for $g>1$ in the case of surfaces without nodes or marked points.

\para
Let $A$ be an open set in $\C^m$ and $\Sigma$ be a surface.
An almost complex structure
on $A\times\Sigma$ with respect to which the projection
$A\times\Sigma\to A$ is holomorphic has the form
$$
J = \left(\begin{array}{cc}
i & 0 \\ \alpha & j
\end{array}\right),
$$
where $j:A\to\cJ(\Sigma)$ is a smooth function with values in the space of
(almost) complex structures on $\Sigma$ and
$\alpha\in\Omega^1(A,\Vect(\Sigma))$ is a smooth $1$-form
on $A$ with values in the space of vector fields on $\Sigma$
such that
$$
\alpha(a,i\hat a) + j(a)\alpha(a,\hat a) = 0
$$
for all $a\in A$ and $\hat a\in T_aA$. This means that the $1$-form $\alpha$
is complex anti-linear with respect to the complex structure on the vector
bundle $A\times\Vect(\Sigma)\to A$ determined by $j$.
From an abstract point of view it is useful to think of
$\alpha$ as a connection on the (trivial) principal bundle
$A\times\Diff(\Sigma)$ and of $j:A\to\cJ(\Sigma)$ as a section
of the associated fiber bundle $A\times\cJ(\Sigma)$.
This section is holomorphic with respect to the Cauchy--Riemann
operator associated to the connection $\alpha$ if and only if
\begin{equation}\label{eq:integrable}
dj(a)\hat a + j(a)dj(a)i\hat a
+ j(a)\cL_{\alpha(a,\hat a)}j(a) =0
\end{equation}
for all $a\in A$ and $\hat a\in T_aA$.
(For a finite dimensional analogue see for
example~\cite{CGS}.)
\arap

\begin{lemma}\label{le:integrable}
$J$ is integrable if and only if $j$ and $\alpha$
satisfy~(\ref{eq:integrable}).
\end{lemma}

\begin{proof}
It suffices to consider the case $m=1$,
so $A\subset\C$ with coordinate $s+it$.
Then the complex structure $J$
on $A\times\Sigma$ has the form
\begin{equation}\label{eq:JASigma}
J = \left(\begin{array}{ccc}
0 & -1 & 0 \\
1 & 0 & 0 \\
-j\xi & -\xi & j
\end{array}\right),
\end{equation}
where $A\to\cJ(\Sigma):s+it\mapsto j_{s,t}$ and
$A\to\Vect(\Sigma):s+it\mapsto\xi_{s,t}$ are smooth maps.
The equation~(\ref{eq:integrable}) has the form
\begin{equation}\label{eq:integrable1}
\p_sj + j\p_tj+\cL_\xi j = 0.
\end{equation}
To see that this is equivalent to integrability of $J$ evaluate the
Nijenhuis tensor $N_J(X,Y):=[JX,JY]-J[X,JY]-J[JX,Y]-[X,Y]$
on a pair of vectors of the form $X=(1,0,0)$, $Y=(0,0,\hat z)$.
The condition $N_J(X,Y)=0$ for all such vectors is equivalent
to~(\ref{eq:integrable1}) and it is easy to see that $N_J=0$
if and only if $N_J((1,0,0),(0,0,\hat z))=0$ for all $\hat z\in T\Sigma$.
The latter assertion uses the facts that $N_J$ is bilinear,
$N_J(X,Y)=-N_J(Y,X)=JN_J(JX,Y)$, and every complex structure
on a $2$-manifold is integrable. This proves the lemma.
\end{proof}

Let $A$ be a complex manifold and $\iota:A\to\cJ(\Sigma)$
be a holomorphic map.  Consider the fibration
$$
\pi_\iota:P_\iota:=A\times\Sigma\to A
$$
with almost complex structure
\begin{equation}\label{eq:Jiota}
      J_\iota(a,z) := \left(\begin{array}{cc}
      i & 0 \\ 0 & \iota(a)(z)\end{array}\right).
\end{equation}
Here we denote by $i$ the complex structure on $A$.
By Lemma~\ref{le:integrable} the almost complex structure
$J_\iota$ on $P_\iota$ is integrable.

\begin{lemma}\label{le:infuniv}
Let $a\in A$.  Then the pair $(\pi_\iota,a)$ is an infinitesimally
universal unfolding if and only if the restriction of $\iota$
to a sufficiently small neighborhood of $a$ is a local slice
as in Theorem~\ref{thm:slice}.
\end{lemma}

\begin{proof}
Let $u:\Sigma\to P_a$ be the diffeomorphism $u(z):=(a,z)$
and denote $j:=\iota(a)$. Then the linearized operator $D_{u,a}$
(at the pair $(u,a)$ for the equation $\bar\p_ju=0$ with $j=\iota(a)$)
has domain $\cX_{u,a}=\Omega^0(\Sigma,T\Sigma)\times T_aA$,
target space $\cY_u=\Omega^{0,1}_j(\Sigma,T\Sigma)$
and is given by
$$
D_{u,a}(\hat u,\hat a)
= \dbar_j\hat u -\frac12 j d\iota(a)\hat a.
$$
(See the formula in~\cite[page~176]{MS2} with $v=\id$.)
This operator is bijective if and only if $d\iota(a)$ is injective
and its image in $T_j\cJ=\Omega^{0,1}_j(\Sigma,T\Sigma)$
is a complement of $\mathrm{im}\,\dbar_j=T_j(\Diff_0(\Sigma)^*j)$
(see the proof of Theorem~\ref{thm:slice}).
This proves the lemma.
\end{proof}

\begin{theorem}\label{thm:teich}
Theorems~\ref{thm:stable}-\ref{thm:existence} hold for Riemann surfaces
of genus $g>1$ without nodes and marked points.
\end{theorem}

\begin{proof} Let $\Sigma$ be a surface of genus $g$. Abbreviate
$$
\cD_0:=\Diff_0(\Sigma), \qquad \cJ:=\cJ(\Sigma), \qquad
\cT:=\cT(\Sigma):=\cJ(\Sigma)/\Diff_0(\Sigma).
$$
Thus $\cT_g:=\cT$ is  Teichm\"uller space.
Consider the principal fiber bundle
$$
\cD_0 \to\cJ\to\cT.
$$
The associated fiber bundle
$$
\pi_{\cT}:\cQ := \cJ\times_{\cD_0}\Sigma\to\cT
$$
has fibers isomorphic to $\Sigma$.

\medskip\noindent{\bf Step~1.}
{\em $\cQ$ and $\cT$ are complex manifolds and $\pi_T$ is a proper
holomorphic submersion.
}

\medskip\noindent 
By Lemma~\ref{le:integrable} with $A=\cJ$ and the map $A\to\cJ$ 
equal to the identity, the space $\cJ\times\Sigma$ is a complex manifold.  
Since $\cD_0$ acts by holomorphic diffeomorphisms, 
so is the (finite dimensional) quotient $\cQ$.

\medskip\noindent{\bf Step~2.}
{\em  The projection $\pi_\cT$ is an
infinitesimally universal unfolding of each of its fibers.
}

\medskip\noindent  Choose $[j_0]\in\cT$.
Let $B$ be an open neighborhood of $0$ in $\C^{3g-3}$ and
$\iota:B\to\cJ$ be a local holomorphic slice such that
$\iota(0)=j_0$ (see~\cite{TROMBA} or Section~\ref{sec:cx}).
Then the projection
$Q_\iota\to B$ is a local coordinate chart on $\cQ\to\cT$.
Hence Step~2 follows from Lemma~\ref{le:infuniv}.

\medskip\noindent{\bf Step~3.}
{\em Every pseudomorphism from $(\pi_A,a_0)$ to $(\pi_\cT,[j_0])$
is a morphism.
}

\medskip\noindent
Let $(\phi,\Phi)$ be a pseudomorphism
from $(\pi_A:P\to A,a_0)$ to $(\pi_\cT,[j_0])$ and
$\iota:B\to\cJ$ be as the proof of in Step~2.
Define $(\psi,\Psi)$ to be the composition of $(\phi,\Phi)$
with the obvious morphism from $(\pi_\cT,[j_0])$
to $(Q_\iota,0)$. Using the maps $\Psi_a:P_a\to\Sigma$
given by $\Psi(p)=:(\psi(a),\Psi_a(p))$ for $p\in P_a$
we construct a trivialization
$$
\tau:A\times\Sigma\to P,\qquad
\tau(a,z):=\tau_a(z):=\Psi_a^{-1}(z).
$$
Then the pullback of the complex structure on $P$
under $\tau$ has the form
$$
      J(a,z) := \left(\begin{array}{cc}
      i & 0 \\ \alpha & j(a)(z)\end{array}\right)
$$
where $j:=\iota\circ\psi:A\to\cJ$ and
$\alpha\in\Omega^{0,1}_j(A,\Vect(\Sigma))$
Since $J$ is integrable it follows from Lemma~\ref{le:integrable}
that $j$ and $\alpha$ satisfy~(\ref{eq:integrable}).
Since the local slice is holomorphic the term
$dj(a)\hat a+j(a)dj(a)i\hat a$ is tangent to the slice while
the last summand $j(a)\cL_{\alpha(a,\hat a)}j(a)=-\cL_{\alpha(a,i\hat a)}j(a)$
is tangent to the orbit of $j(a)$ under $\cD_0$. It follows that
both terms vanish for all $a\in A$ and $\hat a\in T_aA$.
Hence $\alpha=0$ and the map $j:A\to\cJ$ is holomorphic.
Hence $\psi:A\to B$ is holomorphic and hence so is $\Psi$.

\medskip\noindent{\bf Step~4.}
{\em  $\pi_\cT$ is a universal unfolding of each of its fibers.
}

\medskip\noindent Choose an unfolding $(\pi_A:P\to A,a_0)$
and a holomorphic diffeomorphism $u_0:(\Sigma,j_0)\to P_{a_0}$.
Then $u_0^{-1}$ is a fiber isomorphism  from $P_{a_0}$ to $\cQ_{[j_0]}$.
Trivialize $P$ by a map $\tau:A\times\Sigma\to P$ such that $\tau_{a_0}=u_0$.
Define $j:A\to\cJ$ so that $j(a)$ is the pullback of the complex structure
on $P_a$ under $\tau_a$. Then $j(a_0)=j_0$.
Define $\phi:A\to\cT$ and $\Phi:P\to\cQ$ by
$$
\phi(a) := [j(a)],\qquad \Phi(p) := [j(a),z],\qquad p=:\tau(a,z)
$$
for $a\in A$ and $p\in P_a$.  This is a pseudomorphism and hence,
by Step~3, it is a morphism.

To prove uniqueness, choose a local holomorphic slice $\iota:B\to\cJ$
such that $\iota(0)=j_0$. Choose two
morphisms  $(\psi,\Psi),(\phi,\Phi):(\pi_A,a_0)\to(\pi_\iota,0)$
such that $\Phi_{a_0}=\Psi_{a_0}=u_0^{-1}:P_{a_0}\to \Sigma$.
If $a$ is near $a_0$ then
$$
    \Psi_a\circ\Phi_a^{-1}:(\Sigma,\iota(\phi(a)))\to (\Sigma,\iota(\psi(a)))
$$
is a diffeomorphism close to the identity
and hence isotopic to the identity. Hence by the local slice
property $\phi(a)=\psi(a)$ and $\Psi_a\circ\Phi_a^{-1}=\id$.

\medskip\noindent{\bf Step~5.}
{\em Let $j_0$ be a complex structure  on $\Sigma$.
Every infinitesimally universal unfolding
$(\pi_B: Q\to B,b_0)$ of  $(\Sigma_0,j_0)$
is isomorphic to $(\pi_\cT, [j_0])$.
}

\medskip\noindent  As in Step~3 we may assume that
$Q=B\times\Sigma$ with complex structure
$$
J(b,z) := \left(\begin{array}{cc}
i & 0 \\ 0 & \iota(b)(z)\end{array}\right)
$$
where $\iota:B\to\cJ$ is holomorphic and $\iota(b_0)=j_0$.
By Lemma~\ref{le:integrable} this almost complex structure is integrable.
Since $(\pi_B:Q\to B,b_0)$  is   infinitesimally universal,
it follows from Lemma~\ref{le:infuniv} that the restriction of
$\iota$ to a neighborhood of $b_0$ is a local slice.
Hence $(\pi_B,b_0)$ is isomorphic to $(\pi_\cT,[j_0])$
by the local  slice property.
\end{proof}

\rmk\label{rmk:proper}
The universal unfolding $\pi_{\cT}:\cQ\to\cT$
of Theorem~\ref{thm:teich} determines
an etale groupoid $(B,\Gamma)$ with
$
B:=\cT=\cJ(\Sigma)/\Diff_0(\Sigma)
$
and
$$
\Gamma:=\left\{[j,\phi,j']:j,j'\in\cJ(\Sigma),\,
\phi\in\Diff(\Sigma),\,j=\phi^*j'\right\}.
$$
Here $[j,\phi,j']$ denotes the equivalence class
under the diagonal action of $\Diff_0(\Sigma)$
by $\psi^*(j,\phi,j'):=(\psi^*j,\psi^{-1}\circ\phi\circ\psi,\psi^*j')$.
By Lemma~\ref{le:proper} this etale groupoid is proper.
\kmr


\section{The Teichm\"uller space $\cT_{g,n}$}\label{sec:teichmuller-n}

In this section we prove Theorems~\ref{thm:stable}-\ref{thm:existence} 
for all stable marked Riemann surfaces without nodes.  
Let $(\Sigma,s_*,j_0)$ be a stable marked Riemann surface 
of type $(g,n)$ without nodes. 
We will construct an infinitesimally universal unfolding
$(\pi_B,S_*,b_0)$ of $(\Sigma,j_0,s_*)$, prove that it is universal,
and prove that every infinitesimally universal unfolding 
of  $(\Sigma,s_*,j_0)$ is isomorphic to the one we've constructed.

\para\label{universalcurve}
Let $n$ and $g$ be nonnegative integers such that $n>2-2g$ and
let $\Sigma$ be a surface of genus $g$. Abbreviate
$$
\cG:=\Diff_0(\Sigma), \qquad
\cP:=\cJ(\Sigma)\times(\Sigma^n\setminus\Delta),\qquad
\cB:=\cP/\cG,
$$
where $\Delta\subset\Sigma^n$ denotes the fat diagonal, i.e. set of all
$n$-tuples of points in $\Sigma^n$ where at least two components are equal.
Thus $\cB=\cT_{g,n}$ is  the Teichm\"uller space of Riemann surfaces of
genus $g$ with $n$ distinct marked points.
Consider the principal fiber bundle
$$
\cG\to\cP\to\cB.
$$
The associated fiber bundle
$$
\pi_{\cB}:\cQ := \cP\times_{\cG}\Sigma\to\cB
$$
has fibers isomorphic to $\Sigma$ and is equipped with $n$ disjoint sections
$$
\cS_i := \left\{[j,s_1,\dots,s_n,z]\in\cQ\,:\,z=s_i\right\},\qquad
i=1,\dots,n.
$$
It is commonly called the  \jdef{universal curve} of genus $g$
with $n$ marked points.
\arap

\para
Let $(j_0,r_*)\in\cP$, $A$ be a complex manifold, $a_0\in A$, and
$$
\iota=(\iota_0,\iota_1,\dots,\iota_n):A\to\cP
$$
be a holomorphic map such that
\begin{equation}\label{eq:iota0}
\iota_0(a_0)=j_0,\qquad \iota_i(a_0)=r_i,\qquad i=1,\dots,n.
\end{equation}
Define the unfolding $(\pi_\iota:P_\iota\to A,R_{\iota,*},a_0)$ by
\begin{equation}\label{eq:Qiota}
P_\iota := A\times\Sigma,\qquad
J_\iota(a,z) := \left(\begin{array}{cc}
      \sqrt{-1} & 0 \\ 0 & \iota_0(a)(z)\end{array}\right)
\end{equation}
where $\sqrt{-1}$ denotes the complex structure on $A$ and
\begin{equation}\label{eq:Siota}
R_{\iota,i}:=\left\{(a,\iota_i(a))\,:a\in A\right\},\qquad
i=1,\dots,n.
\end{equation}
\arap

\begin{lemma}\label{le:infuniv-n}
The unfolding $(\pi_\iota,R_{\iota,*},a_0)$ is infinitesimally
universal if and only if the
restriction of $\iota$ to a sufficiently small neighborhood of $a_0$
is a (holomorphic) local slice for the action of $\cG$ on $\cP$.
\end{lemma}

\begin{proof}
Let $u_0:(\Sigma,j_0)\to A$ be the holomorphic embedding
$u_0(z):=(a_0,z)$. Then the operator $D_{u_0,a_0}$ has domain
$$
\cX_0:=\left\{(\xi,\hat a)\in\Omega^0(\Sigma,T\Sigma)\times T_{a_0}A\,:\,
\xi(r_i)=d\iota_i(a_0)\hat b\right\}
$$
target space $\cY_0:=\Omega^{0,1}_{j_0}(\Sigma,T\Sigma)$ and is given by
$$
D_{u_0,a_0}(\hat u,\hat a) 
= \dbar_{j_0}\hat u -\frac12 j_0 d\iota_0(a_0)\hat a.
$$
Now the tangent space of the group orbit
$\cG^*(j_0,r_*)$ at $(j_0,r_*)$ is given by
$$
T_{(j_0,r_*)}\cG^*(j_0,r_*)
= \left\{(2j_0\bar\p_{j_0}\hat u,-\hat u(r_1),\dots,-\hat u(r_n))\,:\,
\xi\in\Omega^0(\Sigma,T\Sigma)\right\}.
$$
(See the proof of Theorem~\ref{thm:S2} for the formula
$\cL_{\hat u} j_0=2j_0\bar\p_{j_0}\hat u$.)
Hence the operator $D_{u_0,a_0}$ is injective if and only if
$\mathrm{im}\,d\iota(a_0)\cap T_{(j_0,r_*)}\cG^*(j_0,r_*)=0$
and $d\iota(a_0)$ is injective. It is surjective if and only if 
$\mathrm{im}\,d\iota(a_0)+T_{(j_0,r_*)}\cG^*(j_0,r_*)=T_{(j_0,r_*)}\cP$.
This proves the lemma.
\end{proof}

\begin{theorem}\label{thm:teich-mark}
Theorems~\ref{thm:stable}-\ref{thm:existence} hold for marked Riemann
surfaces without nodes.
\end{theorem}

\begin{proof}   {\bf Step~1.}
{\em $\cQ$ and $\cB$ are complex manifolds,
the projection $\pi_\cB$ is a proper holomorphic submersion,
and $\cS_1,\dots,\cS_n$ are complex submanifolds of $\cQ$.
}

\medskip\noindent
Apply Lemma~\ref{le:integrable} to the complex manifold $A=\cJ=\cJ(\Sigma)$,
replace the fiber $\Sigma$ by $\Sigma^n\setminus\Delta$, and replace $\iota$
by the map $\cJ\to\cJ(\Sigma^n\setminus\Delta)$ which assigns to each complex
structure $j\in\cJ(\Sigma)$ the corresponding product structure on
$\Sigma^n\setminus\Delta$.   Then (the proof of) Lemma~\ref{le:integrable}
shows that $\cP=\cJ\times(\Sigma^n\setminus\Delta)$ is a complex manifold.
The group $\cG=\Diff_0(\Sigma)$ acts on this space by the
holomorphic diffeomorphisms
$$
(j,s_1,\dots,s_n)\mapsto (f^*j,f^{-1}(s_1),\dots,f^{-1}(s_n))
$$
for $f\in\cG$. The action is free and admits holomorphic local slices
for all $g$ and $n$. It follows that the quotient $\cB=\cP/\cG$ is a
complex manifold.  The same argument shows that the total space $\cQ$
is a complex manifold and that the projection $\pi_\cB:\cQ\to\cB$ is
holomorphic. That it is a proper submersion is immediate
from the definitions.

Here are more details on the holomorphic local slices for the action
of $\cG$ on $\cP$. In the case $g>1$ we will find a holomorphic
local slice $\iota:B\to\cP$, defined on $B:=B_0\times\INT(\D)^n$,
which has the form
$$
\iota(b_0,b_1,\dots,b_n)=(\iota_0(b_0),\iota_1(b_0,b_1),\dots,\iota_n(b_0,b_n)).
$$
Here $\iota_0:B_0\to\cJ$ is a holomorphic local slice as in
Theorem~\ref{thm:slice}.  For $i=1,\dots,n$, the map
$(b_0,b_i)\mapsto(b_0,\iota_i(b_0,b_i))$ is holomorphic
with respect to the complex structure $J_{\iota_0}$ on
$Q_0:=B_0\times\Sigma$ defined by~(\ref{eq:Qiota})
and restricts  to a holomorphic
embedding from $b_0\times\INT(\D)$ to $(\Sigma,j)$
with $j=\iota_0(b_0)$.  That such maps $\iota_i$ exist
and can be chosen with disjoint images follows from
Lemma~\ref{le:markedPT}.

In the case $g=1$ and $n\ge 1$ with $\Sigma=\T^2:=\R^2/\Z^2$
an example of a holomorphic local slice is the map
$\iota:B=B_0\times B_1\times\cdots\times B_{n-1}\to\cP$
given by
$$
\iota(\lambda_0,b_1,\dots,b_{n-1})
:= (j(\lambda_0),f_{\lambda_0}^{-1}(b_1),
\dots,f_{\lambda_0}^{-1}(b_{n-1}),f_{\lambda_0}^{-1}(s_n))
$$
where $B_0\subset\H$ and $B_i\subset\C$ are open sets such that closures of
the $n-1$ sets $B_i+\Z+\lambda_0\Z\subset T_{\lambda_0}:=\C/\Z+\lambda_0\Z$
are pairwise disjoint, none of these sets
contains the point $s_n+\Z+\lambda_0\Z$,
the complex structure $j(\lambda_0)\in\cJ(\T^2)$
is defined by~(\ref{eq:jlambda}), and the isomorphism
$f_{\lambda_0}:(\T^2,j(\lambda_0))\to T_{\lambda_0}$ is defined by
$f_{\lambda_0}(x,y):=x+\lambda_0y$.
That any such map is a holomorphic local slice for the action of
$\cG=\Diff_0(\T^2)$ follows from Theorem~\ref{thm:T2}.

In the case $g=0$ and $n\ge 3$ with $\Sigma=S^2$ an example of a
holomorphic local slice is the map $\iota:B=\INT(\D)^{n-3}\to\cP$
given by
$$
\iota(b_1,\dots,b_{n-s})
:= (j_0,\iota_1(b_1),\dots,\iota_{n-3}(b_{n-3}),s_{n-2},s_{n-1},s_n)
$$
$s_{n-2},s_{n-1},s_n$ are distinct points in $S^2$,
$j_0\in\cJ(S^2)$ denotes the standard complex structure,
and the $\iota_i:\INT(\D)\to S^2$ are holomorphic
embeddings for  $1\le i\le n-3$ such that the closures
of their images are pairwise disjoint and do not contain
the points $s_{n-2},s_{n-1},s_n$.  That any such map is a
holomorphic local slice for the action of $\cG=\Diff_0(S^2)$
follows from Theorem~\ref{thm:S2}.

Thus we have constructed holomorphic local slices for the action of
$\cG=\Diff_0(\Sigma)$ on $\cP=\cJ(\Sigma)\times(\Sigma^n\setminus\Delta)$
in all cases. Holomorphic slices for the action of $\cG$ on $\cP\times\Sigma$
can be constructed in a similar fashion. It then follows from the symmetry
of the construction under permutations of the components in $\Sigma$ that
the sections $\cS_i$ are complex submanifolds of $\cQ$.  This proves Step~1.

\medskip\noindent{\bf Step~2.}
{\em  The pair $(\pi_\cB,\cS_*)$ is an
infinitesimally universal unfolding of each of its fibers.
}

\medskip\noindent  Choose $[j_0,s_*]\in\cB$.
Let $B$ be an open neighborhood of $b_0=0$ in $\C^{3g-3+n}$ and
$\iota=(\iota_0,\iota_1,\dots,\iota_n):B\to\cP$ be a local holomorphic
slice satisfying~(\ref{eq:iota0}). Then the unfolding
$(\pi_\iota:Q_\iota\to B,S_{\iota,*},b_0)$ defined as
in~(\ref{eq:Qiota}) and~(\ref{eq:Siota})  is isomorphic to
$(\pi_\cB,\cS_*,[j_0,s_*])$. Hence Step~2 follows from
Lemma~\ref{le:infuniv-n}.

\medskip\noindent{\bf Step~3.}
{\em Every pseudomorphism from $(\pi_A,R_*,a_0)$
to $(\pi_\cB,\cS_*,[j_0,s_*])$ is a morphism.
}

\medskip\noindent
Let $(\phi,\Phi)$ be a pseudomorphism
from $(\pi_A:P\to A,R_*,a_0)$ to $(\pi_\cB,\cS_*,[j_0,s_*])$
and $\iota=(\iota_0,\iota_1,\dots,\iota_n):B\to\cP$ be as the proof of
in Step~2. Define $(\psi,\Psi)$ be the composition of $(\phi,\Phi)$
with the obvious morphism from $(\pi_\cB,\cS_*,[j_0,s_*])$
to the unfolding $(Q_\iota,S_{\iota,*},b_0)$, defined
as in~(\ref{eq:Qiota}) and~(\ref{eq:Siota}). Using the maps
$\Psi_a:P_a\to\Sigma$ given by $\Psi(p)=:(\psi(a),\Psi_a(p))$ for $p\in P_a$
we construct a trivialization
$$
\tau:A\times\Sigma\to P,\qquad
\tau(a,z):=\tau_a(z):=\Psi_a^{-1}(z).
$$
Then the pullback of the section $R_i$ is given by
$$
\tau^{-1}(R_i) = \left\{(a,\sigma_i(a))\,:\,a\in A\right\},\qquad
\sigma_i:=\iota_i\circ\psi:A\to\Sigma
$$
and the pullback of the complex structure on $P$
under $\tau$ has the form
$$
      J(a,z) := \left(\begin{array}{cc}
      \sqrt{-1}& 0 \\ \alpha & j(a)(z)\end{array}\right)
$$
where $j:=\iota_0\circ\psi:A\to\cJ$ and
$\alpha\in\Omega^{0,1}_j(A,\Vect(\Sigma))$.
Since $J$ is integrable it follows from Lemma~\ref{le:integrable}
that $j$ and $\alpha$ satisfy
$$
dj(a)\hat a+j(a)dj(a)\sqrt{-1}\hat a-\cL_{\alpha(a,\sqrt{-1}\hat a)}j(a)=0.
$$
Since $\tau^{-1}(R_i)$ is a complex submanifold of $A\times\Sigma$,
we have
$$
d\sigma_i(a)\hat a + j(a)d\sigma_i(a)\sqrt{-1}\hat a
+ \alpha(a,\sqrt{-1}\hat a)(\sigma_i(a)) = 0
$$
for $i=1,\dots,n$. Since $\iota$ is a local holomorphic slice
these two equations together imply that
$$
dj(a)\hat a+j(a)dj(a)\sqrt{-1}\hat a=0,\qquad
d\sigma_i(a)\hat a + j(a)d\sigma_i(a)\sqrt{-1}\hat a=0,
$$
and $\cL_{\alpha(a,\sqrt{-1}\hat a)}j(a)=0$ and
$\alpha(a,\sqrt{-1}\hat a)(\sigma_i(a))=0$
for all $\hat a\in T_aA$.  Since $n>2-2g$ it follows that $\alpha\equiv0$.
Moreover, the map $(j,\sigma_1,\dots,\sigma_n)=\iota\circ\psi:A\to\cP$
is holomorphic.  Since $\iota$ is a holomorphic local slice, this
implies that $\psi$, and hence also $\Psi$, is holomorphic.

\medskip\noindent{\bf Step~4.}
{\em The pair $(\pi_\cB,\cS_*)$ is a universal unfolding
of each of its fibers.
}

\medskip\noindent Choose $[j_0,s_*]\in\cB$ and let $(\pi_A:P\to A,a_0)$
which admits an isomorphism $u_0:(\Sigma,j_0)\to P_{a_0}$
such that $u_0(s_i):=P_{a_0}\cap R_i$.
Then $u_0^{-1}$ is a fiber isomorphism  from $P_{a_0}$ to $\cQ_{[j_0,s_*]}$.
Trivialize $P$ by a map $\tau:A\times\Sigma\to P$ such that $\tau_{a_0}=u_0$.
Define $j:A\to\cJ$ and $\sigma_i:A\to\Sigma$
so that $j(a)$ is the pullback of the complex structure
on $P_a$ under $\tau_a$ and $\tau^{-1}(R_i)=\{(a,\sigma_i(a))\,:\,a\in A\}$.
Then $j(a_0)=j_0$ and $\sigma_i(a_0)=s_i$.
Define $\phi:A\to\cB$ and $\Phi:P\to\cQ$ by
$$
\phi(a) := [j(a),\sigma_*(a)],\qquad
\Phi(p) := [j(a),\sigma_*(a),z],\qquad p=:\tau(a,z)
$$
for $a\in A$ and $p\in P_a$.  This is a pseudomorphism and hence,
by Step~3, it is a morphism.

To prove uniqueness, choose a local holomorphic
slice $\iota=(\iota_0,\iota_1,\dots,\iota_n):B\to\cP$
such that $\iota(b_0)=(j_0,s_1,\dots,s_n)$.
Choose two morphisms
$$
(\psi,\Psi),(\phi,\Phi):(\pi_A,R_*,a_0)\to(\pi_\iota,S_{\iota,*},b_0)
$$
such that $\Phi_{a_0}=\Psi_{a_0}=u_0^{-1}:P_{a_0}\to \Sigma$.
If $a$ is near $a_0$ then
$$
\Psi_a\circ\Phi_a^{-1}:(\Sigma,\iota_0(\phi(a)))
\to (\Sigma,\iota_0(\psi(a)))
$$
is a diffeomorphism isotopic to the identity
that sends $\iota_i(\phi(a))$ to $\iota_i(\psi(a))$ for $i=1,\dots,n$.
Hence by the local slice property $\phi(a)=\psi(a)$ and $\Phi_a=\Psi_a$.

\medskip\noindent{\bf Step~5.}
{\em Let $j_0$ be a complex structure on $\Sigma$ and
$s_1,\dots,s_n$ be distinct marked points on $\Sigma$.
Every infinitesimally universal unfolding
$(\pi_B: Q\to B,S_*,b_0)$ of  $(\Sigma_0,s_*,j_0)$
is isomorphic to $(\pi_\cB,\cS_*,[j_0,s_*]])$.
}

\medskip\noindent
As in Step~3 we may assume that $Q=B\times\Sigma$ with complex structure
$$
      J(b,z) := \left(\begin{array}{cc}
      i & 0 \\ 0 & \iota_0(b)(z)\end{array}\right)
$$
and $S_i=\left\{(b,\iota_i(b)\,:\,b\in B\right\}$
where $\iota=(\iota_0,\iota_1,\dots,\iota_n):B\to\cP$ is holomorphic
and $\iota(b_0)=(j_0,s_1,\dots,s_n)$.  By Lemma~\ref{le:integrable}
the almost complex structure $J$ is integrable.
Since $(\pi_B:Q\to B,S_*,b_0)$  is infinitesimally universal,
the restriction of $\iota$ to a neighborhood of
$b_0$ is a local slice by Lemma~\ref{le:infuniv-n}.
Hence $(\pi_B,S_*,b_0)$ is isomorphic to $(\pi_\cB,\cS_*,[j_0,s_0])$
by the local  slice property.
\end{proof}

\rmk
The etale groupoid associated to the universal marked
curve of \ref{universalcurve} is proper as in Remark~\ref{rmk:proper}.
\kmr


\section{Nonlinear Hardy spaces}\label{sec:hardy1}

In this section we characterize infinitesimally universal unfoldings in terms
of certain ``nonlinear Hardy spaces'' associated to a desingularization.
The idea is to decompose a Riemann surface $\Sigma$ as
a union of submanifolds $\Omega$ and $\Delta$, intersecting
in their common boundary, and to identify holomorphic maps
on $\Sigma$ with pairs of holomorphic maps defined on $\Omega$
and $\Delta$ that agree on that common boundary.

\para  \label{UV}
Throughout this section we assume that
$$
(\pi_B:Q\to B,S_*,b_0),
$$
is a nodal unfolding of a marked nodal Riemann surface
$(\Sigma,s_*,\nu,j)$ and that
$$
w_0:\Sigma \to Q_{b_0}
$$
is a desingularization.
Let $C_B\subset Q$ denote the set of critical points of $\pi_B$.
Let $U$ be a neighborhood of $C_B$
equipped with nodal coordinates.
This means
$$
    U=U_1\cup\cdots \cup U_k
$$
where the sets $U_i$ have pairwise disjoint closures,
each $U_i$ is a connected neighborhood of one of the components
of $C_B$, and for $i=1,\ldots,k$ there is a holomorphic coordinate
system
$$
     (\zeta_i,\tau_i):B\to \C\times\C^{d-1}, \qquad d:=\dim_\C B
$$
and holomorphic functions $\xi_i,\eta_i:U_i\to\C$ such that
$$
   (\xi_i,\eta_i,\tau_i\circ\pi_B):U_i\to\C\times\C\times\C^{d-1}
$$
is a holomorphic coordinate system and $\xi_i\eta_i=\zeta_i\circ\pi_B$.
Assume that $\bar U\cap S_*=\emptyset$.
Let $V\subset Q$ be an open set such that
$$
     Q=U\cup V, \qquad \bar V\cap C_B=\emptyset,
$$
and $U_i\cap V$ intersects each fiber $Q_b$ in two open annuli
with $|\xi_i|>|\eta_i|$ on one component and $|\xi_i|<|\eta_i|$ on the other.
Introduce the abbreviations
$$
    W:=U\cap V, \quad W_i:=U_i\cap V, \quad W_{i,1}:=\{|\xi_i|>|\eta_i|\},
    \quad W_{i,2}:=\{|\xi_i|<|\eta_i|\},
$$
$$
      U_b:= U\cap Q_b, \qquad  V_b:=V\cap Q_b, \qquad W_b:=W\cap Q_b.
$$
\arap

\para We consider a decomposition
$$
\Sigma=\Omega\cup \Delta,\qquad
\p\Omega=\p \Delta=\Omega\cap \Delta,
$$
into submanifolds with boundary such that  $\Delta$ is a disjoint union
$$
\Delta=\Delta_1\cup\dots\cup \Delta_k
$$
where, for each $i$, the set $\Delta_i$ is either an embedded
closed annulus  or it is the union
of two disjoint embedded closed disks centered at two equivalent nodal points
and
$$
w_0(\Omega)\subset V,\qquad w_0(\Delta_i)\subset U_i
$$
for $i=1,\ldots,k$.
It follows that every pair of equivalent nodal points appears in some $\Delta_i$.
In case $\Delta_i$ is a disjoint union of two disks, say
$\Delta_i=\Delta_{i,1}\cup\Delta_{i,2}$, choose holomorphic
diffeomorphisms $x_i:\Delta_{i,1}\to\D$ and $y_i:\Delta_{i,2}\to\D$
which send the nodal point to $0$.
In case $\Delta_i$ is an annulus
choose a holomorphic diffeomorphism $x_i:\Delta_i\to\A(\delta_i,1)$
and define $y_i:\Delta_i\to\A(\delta_i,1)$ by $y_i=\delta_i/x_i$.
In both cases choose the names so that
$$
w_0\bigl(x_i^{-1}(S^1)\bigr)\subset W_{i,1},\qquad
w_0\bigl(y_i^{-1}(S^1)\bigr)\subset W_{i,2}.
$$
The curves $\xi_i\circ w_0\circ x_i^{-1}$and $ \eta_i\circ w_0\circ y_i^{-1}$
from $S^1$ to $\C\setminus0$
both have winding number one about the origin.
\arap

\para\label{WUV0}
Fix a constant $s>1/2$.
Since $\pi_B|W$ is a submersion the space
$H^s(\p \Delta,W_b)$ is, for each $b\in B$, 
a submanifold of the Hilbert manifold of all $H^s$ maps
from $\p\Delta$ to $W$. Define an  open subset
$$
\cW(b)\subset H^s(\p \Delta,W_b)
$$
by the condition that for $\gamma\in H^s(\p \Delta,W_b)$
we have $\gamma\in\cW(b)$ iff
$$
\gamma\bigl(x_i^{-1}(S^1)\bigr)\subset W_{i,1},\qquad
\gamma\bigl(y_i^{-1}(S^1)\bigr)\subset W_{i,2},
$$
and the curves $\xi_i\circ \gamma\circ x_i^{-1}$ and
$\eta_i\circ \gamma\circ y_i^{-1}$ from $S^1$ to $\C\setminus0$
both have winding number one about the origin.
Introduce sets
\begin{equation*}
\begin{split}
\cZ(b)&:=\{v\in\Hol^s(\Omega,V_b) :
\, v|\p\Delta\in\cW(b) \mbox{ and }
\, v(s_*\cap\Omega)=S_*\cap Q_b
\},
\\
\cN(b)&:=\{u\in\Hol^s(\Delta,U_b):
\, u|\p\Delta\in\cW(b)
\mbox{ and $u$ preserves $\nu$}\}.
\end{split}
\end{equation*}
Here $\Hol^s(X,Y)$ denotes the set of continuous maps from $X$ to $Y$
which are holomorphic on the interior of $X$ and restrict to $H^s$-maps
on the boundary, and the phrase ``$u$ preserves $\nu$'' means
that $\{x,y\}\in\nu\implies u(x)=u(y)\in C_B$.
Define the \jdef{nonlinear Hardy spaces} by
$$
\cU(b):=\{u|\p\Delta\,:\, u\in\cN(b))\},
\qquad
\cV(b):=\{v|\p\Delta\,:\, v\in\cZ(b)\}.
$$
Define
$$
\cW_0:=\bigsqcup_{b\in B}\cW(b),\qquad
\cV_0:=\bigsqcup_{b\in B}\cV(b),\qquad
\cU_0:=\bigsqcup_{b\in B}\cU(b),
$$
so that $(\gamma,b)\in\cW_0\iff\gamma\in\cW(b)$, etc.
The desingularization $w_0:\Sigma\to Q_{b_0}$ determines a point
$$
\gamma_0:= w_0|\p \Delta
\in \cU_0\cap\cV_0\subset\cW_0.
$$
\arap

\begin{lemma}\label{le:ugamma}
For every $(\gamma,b)\in\cU_0\cap\cV_0$
there is a unique desingularization $w:\Sigma\to Q_b$
with $w|\p \Delta=\gamma$.
\end{lemma}

\begin{proof}
Uniqueness is an immediate consequence of unique continuation.
To prove existence, let $(\gamma,b)\in\cU_0\cap\cV_0$ be given.
Then, by definition of $\cU_0$ and $\cV_0$,  there is a continuous
map $w:\Sigma\to Q_b$ which is holomorphic in $\INT(\Omega)$
and in $\INT(\Delta)$ with $w(s_*)=S_*\cap Q_b$ and $w(z_0)=w(z_1)$
for every nodal pair $\{z_0,z_1\}\in\nu$.  The map $w:\Sigma\to Q$
is of class $H^{s+1/2}$ and is therefore holomorphic on all
of $\Sigma$. We must prove that  if $z_0\ne z_1$ we have
$w(z_0)=w(z_1)$ if and only if $\{z_0,z_1\}\in\nu$.
Assume first that there are no nodes, i.e. $Q_b\cap C_B=\emptyset$.
Then $\Delta$ is a union of disjoint annuli and,
by the winding number assumption,
the restriction of $w$ to $\Delta$ is an embedding
into $Q_b$ and $w(\Delta)\cup V_b=Q_b$. Hence there is a point
$q\in w(\Delta)\setminus  V$. Hence the degree
of $w$ at $q$ is one and hence $w:\Sigma\to Q_b$ is a holomorphic
diffeomorphism.  Now assume $Q_b\cap C_B\ne\emptyset$.
Then $w^{-1}(C_B)=\cup\nu$, by the winding number assumption,
and so the restriction $w:\Sigma\setminus\cup\nu\to Q_b\setminus C_B$
is proper. Hence the degree of this restriction is constant on each
component of $\Sigma\setminus\cup\nu$.  Now each component
of $\Sigma$ contains a component $\Delta'$ of $\Delta$ that is
diffeomorphic a disc. By the winding number assumption,
the restriction of $w$ to $\Delta'$ is an embedding.
Moreover, the images under $w$ of the components of
$\Delta\setminus C_B$ are disjoint and there is a point
$q\in w(\Delta')\setminus(C_B\cup V)$. Since
$w(\Omega)\subset V$, the degree of the restriction
$w:\Sigma\setminus\cup\nu\to Q_b\setminus C_B$
at any such point $q$ is one.  Hence the degree of
the restriction is one at every point and hence
$w:\Sigma\setminus\cup\nu\to Q_b\setminus C_B$
is a holomorphic diffeomorphism.  This proves the lemma.
\end{proof}

\begin{theorem}\label{thm:infuniv}
Assume $s>7/2$. Then the spaces $\cU_0$ and $\cV_0$ 
are complex Hilbert submanifolds of $\cW_0$.
Moreover, the  unfolding $(\pi_B,S_*,b_0)$ is 
infinitesimally universal if and only if
$$
T_{\gamma_0}\cW_0
=T_{\gamma_0}\cU_0\oplus T_{\gamma_0}\cV_0.
$$
\end{theorem}

\begin{proof}
{\em We prove that $\cU_0$ is a complex Hilbert
submanifold of $\cW_0$.}

\medskipnoindent
Denote by $H^s$ the Hilbert space of all power series
$$
\zeta(z)=\sum_{n\in\Z} \zeta_nz^n
$$
whose norm
$$
\left\|\zeta\right\|_s:=\sqrt{\sum_{n\in\Z} (1+|n|)^{2s}|\zeta_n|^2}
$$
is finite.
Choose the indexing so that
$z_i(b_0)=0$ for $i\le\ell$ and $z_i(b_0)\ne 0$ for $i>\ell$.
Consider the map
$$
\cW_0\to (H^s)^{2k}\times B:
\gamma\mapsto   (\alpha_1,\beta_1,\dots,\alpha_k,\beta_k,b)
$$
where $\gamma\in\cW(b)$ and
$\alpha_i=\xi_i\circ\gamma\circ x_i^{-1}$
and $\beta_i=\eta_i\circ\gamma\circ y_i^{-1}$.
This maps $\cW_0$ diffeomorphically onto an open set in a Hilbert space.
The map sends $\cU_0\subset\cW_0$ to the subset of all
tuples $(\alpha_1,\beta_1,\dots,\alpha_k,\beta_k,b)$ such that all
nonpositive coefficients of $\alpha_i$ and $\beta_i$ vanish for $i\le\ell$
and such that
$
\beta_i(y) = \zeta_i(b)/\alpha_i(\zeta_i(b)/y)
$
for $i>\ell$. Thus the tuple $(\alpha_1,\beta_1,\dots,\xi_\ell,\eta_\ell)$
is restricted to a closed subspace of $(H^s)^{2\ell}$ and, for $i>\ell$,
the component $\beta_i$ can be expressed as a holomorphic function
of $\alpha_i$ and $b$ (provided that $\zeta_i(b)=0$ for $i\le\ell$ and
$0<|\zeta_i(b)|<2\delta_i$ for $i>\ell$). This shows that $\cU_0$
is a complex Hilbert submanifold of $\cW_0$.

\bigskip

We will show that the restriction map
$$
\cZ_0:=\bigsqcup_{b\in B}\cZ(b)\to \cW_0:v\mapsto v|\p\Omega
$$
is a holomorphic embedding.
Since the image is  precisely $\cV_0$ by definition,
this will show that $\cV_0$ is a complex Hilbert submanifold of $\cW_0$.
Denote by $\cB$ the space of all pairs $(v,b)$,
where $b\in B$ and $v:\Omega\to V_b$ is an
$H^{s+1/2}$ map satisfying
$$
v(s_*) = S_*\cap Q_b.
$$
The space $\cB$ is a complex Hilbert manifold whose tangent space
at $(v,b)$ is the Sobolev space
$$
T_{v,b}\cB = \left\{(\hat v,\hat b)
\in H^{s+1/2}(\Omega,v^*TQ)\times T_bB\,:\,
d\pi_B(v)\hat v\equiv\hat b,\,\hat v(s_i)\in T_{v(s_i)}S_i\right\},
$$
i.e.  $\hat v$ is a section of class $H^{s+1/2}$ of the pullback
tangent bundle $v^*TQ$ that projects to a constant tangent vector
of $B$ and at the marked points is tangent to $S_*$.
Consider the complex Hilbert space bundle
$\cE\to\cB$ whose fiber
$$
\cE_{v,b} := H^{s-1/2}(\Omega,\Lambda^{0,1}T^*\Omega\otimes v^*TQ_b)
$$
over $(v,b)\in\cB$ is the Sobolev space of $(0,1)$-forms
on $\Omega$ of class $H^{s-1/2}$ with values in the
vertical pullback tangent bundle $v^*TQ_b$.
The Cauchy--Riemann operator $\dbar$ is a section
of this bundle and its zero set is the space $\cZ_0\subset\cB$
defined above. The vertical derivative of $\dbar$ at a zero $(v,b)$
is the restriction
$$
D_{v,b}: T_{v,b}\cB\to\cE_{v,b}
$$
of the Cauchy--Riemann operator of the holomorphic vector
bundle $v^*TQ\to\Omega$ to the subspace
$T_{v,b}\cB\subset H^{s+1/2}(\Omega,v^*TQ)$.
This operator is split surjective; a right inverse
can be constructed from an appropriate Lagrangian
boundary condition  (see~\cite[Appendix~C]{MS2}).
Hence $\cZ_0$ is a complex submanifold of~$\cB$.

\bigskip

{\em We show that the restriction map is an
injective  holomorphic immersion.}

\medskipnoindent
By unique continuation at boundary
points, the restriction map is injective, i.e. two elements
of $\Hol^s(\Omega,V_b)$ that agree
on the boundary agree everywhere. By the same
reasoning the derivative of the restriction
map (also a restriction map) is injective.
Its image is the Hardy space of all $H^s$ sections
of the vector bundle $v^*TQ|\p\Omega\to\p\Omega$ that
project to a constant tangent vector of $B$
and extend to holomorphic sections of $v^*TQ_b$
that are tangent to $S_i$ at $s_i$.
By standard Sobolev theory this Hardy space is closed
and the extended section is of class $H^{s+1/2}$.
Hence by the open mapping theorem, the linearized
restriction has a left inverse. Hence the restriction map
$\cZ_0\to \cW_0$ is a holomorphic immersion.

\bigskip

{\em We show that the restriction map $\cZ_0\to \cW_0$ is proper.}

\medskipnoindent
Suppose
that $v_k\in\cZ(b_k)$,
that $\gamma_k:=v_k|\p\Omega$,
that $\gamma_k$ converges to $\gamma\in\cW(b)$,
and
that $\gamma=v|\p\Omega$ where $v\in\cZ(b)$.
We prove in four steps that $v_k$   converges to $v$
in $H^{s+1/2}(\Omega,Q)$.

\medskip\noindent{\bf Step~1.}
{\em
We may assume without loss of generality that each $v_k$
is an embedding for every $k$.
}

\medskip\noindent
After shrinking $\cV_0$ we obtain that $\gamma:\p\Omega\to Q_b$
is an embedding for every $(\gamma,b)\in\cV_0$.
(This makes sense because $s>3/2$, so $\gamma$ is continuously
differentiable.) If $\gamma=v|\p\Omega$ is an embedding
and $v\in\Hol^s(\Omega,V_b)$
then $v$ is an embedding. This is because $\# v^{-1}(q)$
(the number of preimages counted with multiplicity)
for $q\in Q_b\setminus v(\p\Omega)$ can only change
as $q$ passes through the image of $\gamma$.
As $\gamma$ is an embedding $\# v^{-1}(q)$ is either
zero or one.  Hence $v$ is an embedding.

\medskip\noindent{\bf Step~2.}
{\em
A subsequence of $v_k$ converges in the $\Cinf$ topology
on every compact subset of $\INT(\Omega)$.
}

\medskip\noindent
If the first derivatives of $v_k$ are uniformly bounded
then $v_k|\INT(\Omega)$ has a $\Cinf$ convergent
subsequence (see~\cite[Appendix~B]{MS2}).
Moreover, a nonconstant holomorphic
sphere in $Q$ bubbles off whenever the first
derivatives of $v_k$ are not bounded.  But bubbling cannot
occur in $V$.  To see this argue as follows.
Suppose $z_k$ converges to $z_0\in\INT(\Omega)$ and
the derivatives of $v_k$ at $z_k$ blow up.  Then the standard
bubbling argument (see~\cite[Chapter~4]{MS2}) applies.
It shows that, after passing to a subsequence and modifying
$z_k$ (without changing the limit), there are $(i,j_k)$-holomorphic
embeddings $\eps_k$ from the disk $\D_k\subset\C$, centered
at zero with radius $k$, to $Q$ such that $\eps_k(0)=z_k$,
the family of disks $\eps_k(\D_k)$ converges to $z_0$,
and $v_k\circ\eps_k$ converges to a nonconstant $J$-holomorphic
sphere $v_0:S^2=\C\cup\infty\to Q_b$. (The convergence is
uniform with all derivatives on every compact subset of $\C$.)
The image of $v_0$ must intersect the nodal set $Q_b\cap C_B$.
Hence there is a point $a\in\C=S^2\setminus\{\infty\}$ such that
$v_0(a)\in Q\setminus\bar V$. This implies $v_k(\eps_k(a))\notin V$ 
for $k$ sufficiently large, contradicting the fact that
$v_k(\Omega)\subset V$.

\medskip\noindent{\bf Step~3.}
{\em
A subsequence of $v_k$ converges to $v$
in the $C^0$ topology.
}

\medskip\noindent
By Arz\'ela--Ascoli it suffices to show that the sequence
$v_k$ is bounded in $C^1$.  We treat this as a Lagrangian
boundary value problem.  Choose $M\subset Q_b$ to be
a submanifold with boundary that contains the image
of $v$ in its interior. Choose a smooth family
of embeddings
$$
\iota_a:M\to Q_a\setminus C_B,\qquad a\in B,
$$
such that $\iota_b:M\to Q_b\setminus C_B$ is the inclusion.
Then the image of $\iota_{b_k}$ contains the image of $v_k$
for $k$ sufficiently large. Think of $M$ as a symplectic manifold
and define the Lagrangian submanifolds $L\subset M$
and $L_k\subset M$ by
$$
L:=\gamma(\p\Omega),\qquad
L_k := \iota_{b_k}^{-1}\circ\gamma_k(\p\Omega).
$$
Since $s>7/2$ the sequence $\iota_{b_k}^{-1}\circ\gamma_k:\p\Omega\to M$
converges to $\gamma$ in the $C ^3$ topology. Hence there is a sequence
of diffeomorphisms $\phi_k:M\to M$ such that $\phi_k$ converges
to the identity in the $C^3$ topology and
$$
\phi_k\circ\iota_{b_k}^{-1}\circ\gamma_k = \gamma,\qquad
\phi_k(L_k)=L.
$$
Define
$$
\tilde v_k := \phi_k\circ\iota_{b_k}^{-1}\circ v_k,\qquad
\tilde J_k := (\phi_k\circ\iota_{b_k}^{-1})_*J_{b_k},
$$
where $J_a$ denotes the complex structure on $Q_a$.
Then $\tilde  J_k$ converges to $\tilde J:= J_b$ in the $C^2$
topology,
$\tilde v_k:\Omega\to M$ is a $\tilde J_k$-holomorphic
curve such that $\tilde v_k(\p\Omega)\subset L$ and, moreover,
\begin{equation}\label{eq:allEqual}
\tilde v_k|\p\Omega = \gamma:\p\Omega\to L
\end{equation}
for all $k$. We must prove that the first derivatives of
$\tilde v_k$ are uniformly bounded.  Suppose by contradiction
that that there is a sequence $z_k\in\Omega$ such that
$$
c_k := \left|d\tilde v_k(z_k)\right|
= \left\|d\tilde v_k\right\|_{L^\infty}\to\infty.
$$
Now apply the standard rescaling argument:
Assume w.l.o.g. that $z_k$ converges to $z_0\in\Omega$,
choose a coordinate chart from a neighborhood of $z_0$
to upper half plane (sending $z_k$ to $\zeta_k$),
and compose the resulting $\tilde J_k$-holomorphic curve
with the rescaling map $\eps_k(\zeta):=\zeta_k+\zeta/c_k$.
Let $d_k$ be the Euclidean distance of $\zeta_k$ from the
boundary of the upper half plane. There are two cases.
If $c_k\cdot d_k\to\infty$ then a nonconstant holomorphic
sphere bubbles off and the same argument as in Step~2
leads to a contradiction. If the sequence $c_k\cdot d_k$ is bounded
then, by~\cite[Theorem~B.4.2]{MS2}, the rescaled sequence has
a subsequence that converges in the $C ^1$ topology to a
holomorphic curve $\tilde w:\{\zeta\in\C:\mathrm{im}\,\zeta\ge 0\}\to M$
with $\tilde w(\R)\subset L$. The choice of the rescaling factor shows
that the derivative of $\tilde w$ has norm one at some point and so
$\tilde w$ is nonconstant. On the other hand, since $\eps_k$
converges to a constant, condition~(\ref{eq:allEqual}) implies that the
restriction of this holomorphic curve to the boundary is constant;
contradiction.

\medskip\noindent{\bf Step~4.}
{\em
A subsequence of $v_k$ converges to $v$
in the $H^{s+1/2}$ topology.
}

\medskip\noindent
Choose local coordinates on source and target near the image
of the boundary and use the fact that the Hardy space
for the annulus $\A(r,1)$ is closed in $H^s(\p \A(r,1))$.
(In the notation of~\ref{Hrs} below this Hardy space
is the diagonal in $H_r^s\times H^s$.)

Thus we have proved that every subsequence of $v_k$
has a further subsequence converging to $v$ in $H^{s+1/2}$.
Hence the sequence $v_k$ itself converges to $v$ in  the
$H^{s+1/2}$ topology.

\bigskip

{\em We prove that the unfolding $(\pi_B,S_*,b_0)$
is infinitesimally universal if and only if
$\cU_0$ and $\cV_0$ intersect transversally
at $\gamma_0:=w_0|\p \Delta$.}

\medskipnoindent
Recall the linear operator
$$
     D_{w_0,b_0}:\cX_{w_0,b_0}\to\cY_{w_0}
$$
of Definition~\ref{def:D}. We shall prove the following
assertions.
\begin{description}
\item[(I)]
The operator $D_{w_0,b_0}$ is injective if and only if
$$
T_{\gamma_0}\cU_0\cap T_{\gamma_0}\cV_0=0.
$$
\item[(II)]
The operator $D_{w_0,b_0}$ is surjective if and only if
$$
T_{\gamma_0}\cU_0 + T_{\gamma_0}\cV_0
=T_{\gamma_0}\cW_0.
$$
\end{description}
The first assertion is obvious, because the intersection of
the tangent spaces  $T_{\gamma_0}\cU_0$ and
$T_{\gamma_0}\cV_0$ is, by definition, precisely the set of
all maps $\xi|\p \Delta:\Gamma\to\gamma_0^*TQ$, where
$\Gamma:=\p\Delta$ and $\xi:\Sigma\to w_0^*TQ$ runs over all elements
in the kernel of the operator $D_{w_0,b_0}$. Since $D_{w_0}\xi=0$
it follows from unique continuation that $\xi$ vanishes
identically if and only if its restriction to the disjoint
union $\Gamma$ of circles vanishes. (The fibers are
connected and so $\Gamma$ intersects each component of
$\Sigma$ in at least one circle.) This proves~(I).

To prove~(II), assume first that $D_{w_0,b_0}$ is onto
and let $(\hat\gamma,\hat b)\in T_{(\gamma_0,b_0)}\cW_0$.
Choose any two vector fields $\xi_u$ along
$u_0:=w_0|\Delta$ and $\xi_v$ along $v_0:=w_0|\Omega$
(in the appropriate Sobolev spaces) that project each to a
constant vector $\hat b_u:=d\pi_B(u_0)\xi_u\in T_bB$,
respectively $\hat b_v:=d\pi_B(v_0)\xi_v\in T_bB$,
and satisfy
$$
     \xi_u|\p \Delta-\xi_v|\p \Delta
     = \hat\gamma:\p \Delta\to\gamma_0^*TQ,\qquad
     \hat b_u-\hat b_v=\hat b.
$$
Next define $\eta\in\cY_{w_0}$ by
$$
     \eta|\Delta:=D_{w_0}\xi_u,\qquad
     \eta|\Omega:=D_{w_0}\xi_v.
$$
Using surjectivity of $D_{w_0,b_0}$ we find a vector field
$\hat w_0$ along $w_0:\Sigma\to Q$ such that $\hat b_0=d\pi_B(w_0)\hat w_0$
is a constant vector in $T_bB$ and
$$
     \eta = D_{w_0,b_0}(\hat w_0,\hat b_0).
$$
Abbreviate
$$
     \hat u:=\xi_u-\hat w_0|\Delta,\qquad
     \hat v:=\xi_v-\hat w_0|\Omega.
$$
Then
$$
     (\hat u,\hat b_u-\hat b_0)\in T_{(\gamma_0,b_0)}\cU_0,\qquad
     (\hat v,\hat b_v-\hat b_0)\in T_{(\gamma_0,b_0)}\cV_0,
$$
and
$$
\hat u|\p \Delta-\hat v|\p \Delta=\hat\gamma,\qquad
(\hat b_u-\hat b_0)-(\hat b_v-\hat b_0)=\hat b.
$$
This shows that  $T_{(\gamma_0,b_0)}\cW_0
=T_{(\gamma_0,b_0)}\cU_0+ T_{(\gamma_0,b_0)}\cV_0$.

Conversely, if we assume transversality then
surjectivity of the operator $D_{w_0,b_0}$
follows by reversing the above argment.
Namely, choose vertical vector fields $\xi_u$ along $u_0$ and $\xi_v$
along $v_0$ that vanish at the marked and nodal points such that
$$
     D_{w_0}\xi_u=\eta|\Delta,\qquad
     D_{w_0}\xi_v=\eta|\Omega.
$$
Now use the transversality hypothesis to modify
$\xi_u$ and $\xi_v$ on each half so that they agree
over the common boundary.  It then follows that the
operator $D_{w_0,b_0}$ is onto.  This completes
the proof of Theorem~\ref{thm:infuniv}.
\end{proof}

\rmk The strategy for the proof of the universal unfolding theorem is
to assign to each unfolding $(\pi_A:P\to A,R_*,a_0)$
of the marked nodal Riemann surface $(\Sigma,s_*,\nu,j)$
a family
of Hilbert submanifolds $\cU_a, \cV_a\subset\cW_a$ as in~\ref{WUV0}
parametrized by $a\in A$.  Transversality will then imply that for
each $a$ near $a_0$ there is a unique intersection point
$(\gamma_a,b_a)\in\cU_a\cap\cV_a$ near $(\gamma_0,b_0)$.
Then the fiber isomorphisms $f_a:P_a\to Q_{b_a}$ determined by
the $\gamma_a$ as in Lemma~\ref{le:ugamma} will fit together to
determine the required morphisms $P\to Q$ of nodal families.
The key point is to show that the submanifolds
$\cU_a$ fit together to form a complex submanifold
$$
\cU:=\bigsqcup_{a\in A}\cU_a
\quad\subset\quad
\cW:=\bigsqcup_{a\in A}\cW_a
$$
(see Theorem~\ref{thm:UV} below).  We begin
by studying a local model near a given nodal point
in the next section.
\kmr


\section{The local model}\label{sec:local}

\para\label{standardNode}
Consider the \jdef{standard node} defined as the map
$$
N\to\INT(\D): (x,y)\mapsto xy, \quad N:=\{(x,y)\in\D^2: |xy|<1\}.
$$
For $a\in\INT(\D)$ and $b\in\C$ denote
$$
N_a:=\{(x,y)\in\D^2:xy=a\}, \qquad Q_b:=\{(x,y)\in\C^2:xy=b\}.
$$
We study the set of all quadruples $(a,\xi,\eta,b)$
where $a,b\in\C$ are close to $0$ and
$$
(\xi,\eta):N_a\to Q_b
$$
is a holomorphic map.
If $a\ne0$ this means
$\xi(z),\eta(z)$ are holomorphic functions
on the annulus $|a|\le |z|\le 1$ that are close to the identity,
and satisfy the condition
\begin{equation}\label{eq:ab}
       xy=a\implies \xi(x)\eta(y)=b
\end{equation}
for  $|a|\le |x|\le 1$ and $|a|\le|y|\le 1$.
If $a\ne 0$, this condition implies that $b\ne 0$.
When $a=0$, the functions  $\xi$ and $\eta$ are defined on the
closed unit disk and vanish at the origin; hence $b=0$.
\arap

\para
Fix $s>1/2$.
Let $H^s$ be the Sobolev space of of all power series
\begin{equation}\label{eq:zeta}
\zeta(z)=\sum_{n\in\Z} \zeta_nz^n.
\end{equation}
which have $s$  derivatives in $L^2$ on the unit circle,
i.e. such that the norm
$$
\left\|\zeta\right\|_s:=\sqrt{\sum_{n\in\Z} (1+|n|)^{2s}|\zeta_n|^2}
$$
is finite. For $r>0$  the \jdef{rescaling map} $z\mapsto rz$
 maps the unit circle to the circle of radius $r$.
 Denote by $\zeta_r$ the result of conjugating $\zeta$ by this map,
 i.e.
 $$
     \zeta_r(z):=r^{-1}\zeta(rz).
 $$
The norm
$\|\zeta_r\|_s$ is finite if and only if the series $\zeta$
converges to an $H^s$ function on the circle of radius $r$.
\arap

\para
For $\delta>0$ define the open set
$\cW_\delta\subset \C\times H^s\times H^s\times\C$  by
$$
\cW_\delta:=\{(a,\xi,\eta,b) :
    \|\xi-\id\|_s<\delta,\; \|\eta-\id\|_s<\delta,\; |a|<\delta\}.
$$
Define $\cU_\delta\subset\cW_\delta$ to be the set of
those quadruples $(a,\xi,\eta,b)\in\cW_\delta$
which satisfy~(\ref{eq:ab}). More precisely
if $(a,\xi,\eta,b)\in\cW_\delta$, then
for $a\ne 0$ we have
$$
(a,\xi,\eta,b)\in\cU_\delta
\qquad\iff\qquad
\begin{cases}
\|\xi_{|a|}\|_s<\infty,\;
\|\eta_{|a|}\|_s<\infty,\;
\mbox{ and } &\\ \\
\xi(x)\eta(ax^{-1})=b
\mbox{ for } |a|\le|x|\le 1 &
\end{cases}
$$
while for $a=0$ we have
$$
(0,\xi,\eta,b)\in\cU_\delta
\qquad\iff\qquad
b=0 \mbox{ and } \xi(0)=\eta(0)=0.
$$
Thus $\cU_\delta$ is the space of (boundary values of) 
local holomorphic fiber isomorphisms in the standard model. 
The main result of this section is that $\cU_\delta$ is a 
manifold:
\arap

\begin{theorem}\label{thm:localmodel}
Let $s$ be a positive integer. Then,
for $\delta>0$ sufficiently small, the set $\cU_\delta$
is a complex submanifold of the open set
$\cW_\delta\subset \C\times H^s\times H^s\times\C$.
\end{theorem}
 
 The proof occupies the rest of this section.
Using the Hardy space decomposition defined
in~\ref{Hrs} we formulate three propositions
which define a map $\cT$ whose graph lies in $\cU_\delta$.
We then prove six lemmas, then we prove the three propositions,
and finally we prove that the graph of $\cT$ is exactly equal to $\cU_\delta$.

\para\label{Hrs}
A holomorphic function $\zeta(z)$ defined on an annulus centered
at the origin has a Laurent expansion of the form~(\ref{eq:zeta}).
We write $\zeta=\zeta_++\zeta_-$ where
$$
\zeta_+(z):=\sum_{n>0} \zeta_nz^n,\qquad
\zeta_-(z):=\sum_{n\le0} \zeta_nz^n.
$$
For $r>0$ and $s>1/2$ introduce the norm
$$
\left\|\zeta\right\|_{r,s}
:=\sqrt{\sum_{n\in\Z} (1+|n|)^{2s}r^{2n-2}|\zeta_n|^2}
$$
so that $\zeta_+$ converges inside the circle
of radius $r$ if $\|\zeta_+\|_{r,s}<\infty$
and $\zeta_-$ converges outside the circle of
radius $r$ if $\|\zeta_-\|_{r,s}<\infty$.
Let
$$
H^s_r:=\{\zeta:\|\zeta\|_{r,s}<\infty\}
$$
and $H^s_{r,\pm}$ be the subspace of those
$\zeta$ for which $\zeta=\zeta_\pm$ so we have
the Hardy space decomposition
$$
H^s_r=H^s_{r,+}\oplus H^s_{r,-}.
$$
Then $H^s=H^s_1$ and $\|\cdot\|_s=\|\cdot\|_{1,s}$. We abbreviate
$$
H^s_\pm:=H^s_{1,\pm}.
$$
We view the ball of radius $\delta$ about $\id$ in the Hilbert space  $H^s_r$
as a space  of $H^s$-maps from  the circle of radius $r$
to a neighborhood of this circle;
the norm on $H^s_r$ is defined so that conjugation by the rescaling map
$z\mapsto rz$  induces an isometry $H^s_r\to H^s:\zeta\to\zeta_r$, i.e.
\begin{equation}\label{eq:rescale}
\|\zeta\|_{r,s}= \|\zeta_r\|_s, \qquad \zeta_r(z):=r^{-1}\zeta(rz).
\end{equation}
\arap

\begin{proposition}[Existence]\label{prop:localmodel1}
For every $s>1/2$ there are positive constants $\delta$
and $c$ such that the following holds.
If $a\in\C$ with $0<r:=\sqrt{|a|}\le1$
and $\xi_+,\eta_+\in H^s_+$ satisfy
\begin{equation}\label{eq:deltaPlus}
\left\|\xi_+-\id\right\|_s<\delta,\qquad
\left\|\eta_+-\id\right\|_s < \delta
\end{equation}
then there exists a triple
$(b,\xi_-,\eta_-)\in\C\times H^s_{r^2,-}\times H^s_{r^2,-}$
such that $\xi:=\xi_+ +\xi_-$ and $\eta:=\eta_+ +\eta_-$
satisfy the equation
\begin{equation}\label{eq:xi-eta-a-b}
\xi(x)\eta(ax^{-1})=b
\end{equation}
for $r^2\le|x|\le1$ and
\begin{equation}\label{eq:continuous}
\left|ba^{-1}-\xi_1\eta_1\right|
+\left\|\xi_-\right\|_{r,s}+\left\|\eta_-\right\|_{r,s}
\le 2cr\bigl(\left\|\xi_+-\id\right\|_s+\left\|\eta_+-\id\right\|_s\bigr).
\end{equation}
\end{proposition}

\begin{proposition}[Uniqueness]\label{prop:localmodel2}
For every positive integer $s$  there exist positive constants $\delta$
and $\eps$ such that the following holds.
If $a,b,b'\in\C$, $\xi_+,\eta_+\in H^s_+$ and
$\xi_-,\eta_-,\xi'_-,\eta'_-\in H^s_{r,-}$
with $0<r:=\sqrt{|a|}<1$ satisfy~(\ref{eq:deltaPlus}) and
$$
\left\|\xi_-\right\|_{r,s}<\eps,\qquad
\left\|\eta_-\right\|_{r,s}<\eps,\qquad
\sup_{|x|=r}\left|\xi'_-(x)\right|<r\eps,\qquad
\sup_{|y|=r}\left|\eta'_-(y)\right| < r\eps,
$$
and if $(a,\xi:=\xi_++\xi_-,\eta:=\eta_++\eta_-,b)$
and $(a,\xi':=\xi_++\xi'_-,\eta':=\eta_++\eta'_-,b')$
satisfy~(\ref{eq:xi-eta-a-b}) for $|x|=r$
then
$
(\xi_-,\eta_-,b)=(\xi'_-,\eta'_-,b').
$
\end{proposition}

\para\label{cT}
Fix a positive integer $s$.  Choose positive constants $\delta$
and $\eps$ such that Proposition~\ref{prop:localmodel2} holds.
Shrinking $\delta$ if necessary we may assume that
Proposition~\ref{prop:localmodel1} holds with the
same constant $\delta$ and a suitable constant $c>0$.
Let
$$
H^s_+(\id,\delta):=\{\zeta\in H^s_+:\|\zeta-\id\|_s<\delta\}
$$
and define
$$
\cT:\D\times H^s_+(\id,\delta)\times H^s_+(\id,\delta)
\to \C\times H^s_-\times H^s_-
$$
by the conditions that $\cT(a,\xi_+,\eta_+)=(b,\xi_-,\eta_-)$
is the triple constructed in Proposition~\ref{prop:localmodel1}
for $a\ne0$ and
$$
      \cT(0,\xi_+,\eta_+):=(0,0,0).
$$
(In defining $\cT$ we used the fact that
$H^s_{r,-}\subset H^s_-$ for $r\le 1$.)
\arap

\begin{proposition}\label{prop:localmodel3}
The map $\cT$ is continuous. It is holomorphic for $|a|<1$.
\end{proposition}

\begin{lemma}[A priori estimates]\label{le:apriori}
There is a constant $c>0$ such that, for $\delta>0$ sufficiently small,
the following holds. If $(a,\xi,\eta,b)\in\cU_\delta$ and $a\ne 0$
then
$$
\left|ba^{-1}-1\right|<c\delta,\quad
\sup_{|a|\le|x|\le1}\left|\xi(x)x^{-1}-1\right|\le c\delta,\quad
\sup_{|a|\le|y|\le1}\left|\eta(y)y^{-1}-1\right|\le c\delta.
$$
\end{lemma}

\begin{proof}
Rewrite  $\xi(x)\eta(ax^{-1})=b $ as
$$
      a\frac{\xi(x)}{x^2}=ab\frac{x^{-2}}{\eta(ax^{-1})}.
$$
Using the substitution $y=ax^{-1}$,  $dy=-ax^{-2}\,dx$ we get
$$
a\oint_{|x|=1}\frac{\xi(x)\,dx}{x^2}=
b\oint_{|x|=1}\frac{ax^{-2}\,dx}{\eta(ax^{-1})} =
b\oint_{|y|=|a|}\frac{dy}{\eta(y)} =
b\oint_{|y|=1}\frac{dy}{\eta(y)}
$$
where all the contour integrals are counter clockwise.
By the Sobolev embedding theorem, there is a constant $c$ 
such that $|\zeta(z)|\le c\|\zeta\|_s$ for $\zeta\in H^s$ and $|z|=1$.
This gives the estimate
$$
\left|\frac{\xi(x)}{x^2}-\frac{1}{x}\right|
=|\xi(x)-x|\le c\|\xi-\id\|_s\le c\delta
$$
for $|x|=1$. If $\left\|\eta-\id\right\|_s<\delta\le 1/2c$ then, 
by the Sobolev embedding theorem again, 
$|\eta(y)-y|<1/2$ and so $|\eta(y)|>1/2$ for $|y|=1$.
Hence
$$
\left|\frac{1}{\eta(y)}-\frac{1}{y}\right|
= \frac{\left|\eta(y)-y\right|}{\left|\eta(y)\right|}
\le 2c\|\eta-\id\|_s\le 2c\delta
$$
for $|y|=1$. Hence the contour integrals are within
$4\pi c\delta$ of $2\pi i$ and so, enlarging $c$,
$b/a$ is within $c\delta$ of $1$ as required.

By symmetry the third inequality follows from the second;
we prove the second. Using the Sobolev inequality we have
$$
\sup_{|x|=1}\left|\frac{\xi(x)}{x}-1\right|\le c\delta,\qquad
\sup_{|y|=1}\left|\frac{\eta(y)}{y}-1\right|\le c\delta.
$$
Now let $y:=ax^{-1}$ and $|x|=|a|$.  Then $|y|=1$ and
$$
\frac{\xi(x)}{x}-1 = \frac{b}{a}\frac{y-\eta(y)}{\eta(y)} + \frac{b}{a}-1.
$$
Hence
$$
\sup_{|x|=|a|}\left|\frac{\xi(x)}{x}-1\right|
\le \left|\frac{b}{a}\right|\sup_{|y|=1}
    \left|\frac{y-\eta(y)}{\eta(y)}\right|
    + \left|\frac{b}{a}-1\right|
\le c\delta.
$$
By the maximum principle this implies
$$
\sup_{|a|\le|x|\le1}\left|\frac{\xi(x)}{x}-1\right|\le c\delta.
$$
This proves the lemma.
\end{proof}

\para\label{cFr}
The proofs of Propositions~\ref{prop:localmodel1}, \ref{prop:localmodel2},
and~\ref{prop:localmodel3} are based on a version
of the implicit function theorem for the map
$$
\cF_r:\C\times H^s_r\times H^s_r\to H^s
$$
defined by
\begin{equation}\label{eq:Fr}
    \cF_r(\lambda,\xi,\eta)(z) := r^{-2}\xi(rz)\eta(rz^{-1})-\lambda
\end{equation}
for $|z|=1$ and $r>0$.
The zeros of $\cF_r$ are solutions of~(\ref{eq:xi-eta-a-b}) with
$a=r^2$  and $b=\lambda a$.
Note that $\cF_r(1,\id,\id)=0$ for every $r>0$.
The differential of $\cF_r$ at the point $(1,\id,\id)$
will be denoted by
$$
\cD_r:=d\cF_r(1,\id,\id):\C\times H^s_r\times H^s_r\to H^s.
$$
Thus
$$
\cD_r(\hat\lambda,\hat\xi,\hat\eta)(z)
= r^{-1}z^{-1}\hat\xi(rz)+r^{-1}z\hat\eta(rz^{-1})-\hat\lambda.
$$
We shall need six lemmata.
They are routine consequences of well known facts and rescaling.
To ease the exposition we relegate the proofs
of the first five to the end of the section and omit
the proof of the sixth entirely.
(The proof of the sixth is just the proof
of the implicit function theorem keeping track of the estimates.)
\arap

\begin{lemma}[Sobolev Estimate]\label{le:s}
Denote by $\A(r,R)\subset\C$ the closed annulus $r\le|z|\le R$.
For every $s>1/2$ there is a constant $c>0$ such that
$$
     \left\|\zeta\right\|_{L^\infty(\A(r,R))}
     \le c\bigl(r\|\zeta_-\|_{r,s}+R\|\zeta_+\|_{R,s}\bigr)
$$
for all $r$ and $R$ and every holomorphic function
$\zeta(z)$ on the annulus $r\le|z|\le R$.
\end{lemma}

\begin{proof}
The constant is
$$
     c=\sqrt{1+\sum_{n=1}^\infty n^{-2s}}.
$$
If $r\le|x|\le R$, then
\begin{eqnarray*}
|\zeta(z)|&\le &
\ds\sum_{n=-\infty}^\infty\left| \zeta_n\right| \left| z\right|^n
\\ &\le &
\sum_{n> 0}|\zeta_n|\,R^n+\sum_{n\le0} |\zeta_n|\,r^n
\\ &= &
\sum_{n>0}(1+|n|)^s|\zeta_n|\,R^n\,(1+|n|)^{-s}
+\sum_{n\le0}(1+|n|)^s |\zeta_n|\,r^n(1+|n|)^{-s}
\\ &\le &
c\bigl(R\|\zeta_+\|_{R,s}+r\|\zeta_-\|_{r,s}\bigr)
\end{eqnarray*}
where the last step is by the Cauchy--Schwarz
inequality.
\end{proof}

\begin{lemma}[Product Estimate]\label{le:prod}
For every positive integer $s$ there is a positive constant $C$ such
that, for any two functions $\xi,\eta\in H^s$, we have
$$
\|\xi\eta\|_s\le C\|\xi\|_s\|\eta\|_s,\qquad
\|\xi\eta\|_s\le C\left(
\|\xi\|_s\|\eta\|_{L^\infty(S^1)}
+\|\xi\|_{L^\infty(S^1)}\|\eta\|_s
\right).
$$
\end{lemma}

\begin{proof} 
The second inequality implies the first by the Sobolev 
embedding theorem, but the first is easy to prove more 
generally for $s>1/2$ so we provide a separate proof.
The constant is
$$
C=\sup_{k\in\Z}\sqrt{\sum_{n\in\Z}
\frac{(1+|k|)^{2s}}{(1+|k-n|)^{2s}(1+|n|)^{2s}}}.
$$
To see that this constant is finite assume $k>0$
and consider the sum over the four regions
$n\le 0$, $0\le n\le k/2$, $k/2\le n\le k$, and $n\ge k$.
By equation~(\ref{eq:rescale}) it is enough to prove the
inequality when $r=1$. Now
\begin{eqnarray*}
\|\xi\|_s^2 &=&
\sum_{k\in\Z} (1+|k|)^{2s}
\left|\sum_{n\in\Z}\xi_{k-n}\eta_n\right|^2 \\
&\le &
\sum_{k\in\Z}(1+|k|)^{2s}\left(
\sum_{n\in\Z}|\xi_{k-n}|\,|\eta_n|
\right)^2 \\
&\le &
\sum_{k\in\Z}\left(
\sum_{n\in\Z}\frac{(1+|k|)^{2s}}{(1+|k-n|)^{2s}(1+|n|)^{2s}}
\right)\cdot \\
&&
\cdot\left(\sum_{n\in\Z}(1+|k-n|)^{2s}|\xi_{k-n}|^2\,
(1+|n|)^{2s}\,|\eta_n|^2\right)
\\&\le&
C^2\sum_{k\in\Z}\sum_{n\in\Z}
(1+|k-n|)^{2s}|\xi_{k-n}|^2
(1+|n|)^{2s}|\eta_n|^2 \\
&=&
C^2\left\|\xi\right\|_s^2
\left\|\eta\right\|_s^2.
\end{eqnarray*}
The third step follows from the Cauchy--Schwarz inequality.
This proves the first inequality of Lemma~\ref{le:prod}.

We prove the second inequality in the case where $s$ is an integer.
The proof is based on the interpolation estimate
of Gagliardo--Nirenberg (see Friedman~\cite{GN}).
In our case it has the form
\begin{equation}\label{eq:GN}
\left\|\zeta\right\|_{W^{k,p}}
\le c_{k,s}\left\|\zeta\right\|_s^\alpha
\left\|\zeta\right\|_{L^\infty(S^1)}^{1-\alpha},\qquad
\alpha:=\frac{k}{s},\qquad
p:=\frac{2s}{k}
\end{equation}
for $\zeta\in H^s$ and $0\le k\le s$.
Let $p:=2s/k$ and $q:=2s/\ell$ where $k+\ell=s$.
Then $1/p+1/q=1/2$ and hence, by H\"older's inequality
and~(\ref{eq:GN}), we have
\begin{eqnarray*}
\left\|\p^k\xi\cdot\p^\ell\eta\right\|_{L^2}
&\le&
\left\|\xi\right\|_{W^{k,p}}\left\|\eta\right\|_{W^{\ell,q}} \\
&\le&
c\left\|\xi\right\|_s^{k/s}
\left\|\xi\right\|_{L^\infty}^{1-k/s}
\left\|\eta\right\|_s^{\ell/s}
\left\|\eta\right\|_{L^\infty}^{1-\ell/s}\\
&=&
c
\left\|\xi\right\|_s^{k/s}
\left\|\eta\right\|_{L^\infty}^{k/s}
\left\|\eta\right\|_s^{\ell/s}
\left\|\xi\right\|_{L^\infty}^{\ell/s} \\
&\le &
\frac{ck}{s}\left\|\xi\right\|_s\left\|\eta\right\|_{L^\infty}
+\frac{c\ell}{s}\left\|\eta\right\|_s\left\|\xi\right\|_{L^\infty}
\end{eqnarray*}
where $c:=c_{k,s}c_{\ell,s}$.
The last inequality follows from
$ab\le a^p/p+b^q/q$ for $a,b\ge 0$
and $p,q\ge1$ with $1/p+1/q=1$.
The desired estimate follows by summing over all
pairs $(k,\ell)$ with $k+\ell=s$ and using the product
rule for differentiation.
\end{proof}

\begin{lemma}[Linear Estimate]\label{le:inverse}
For $\hat\xi_-,\hat\eta_-\in H^s_{r,-}$
and $\hat\lambda\in\C$ we have
$$
        \|\hat\xi_-\|_{r,s}^2
        + \|\hat\eta_-\|_{r,s}^2
        + |\hat\lambda|^2
        \le \|\cD_r(\hat\lambda,\hat\xi_-,\hat\eta_-)\|_s^2.
$$
\end{lemma}

\begin{proof}
The formula
$$
\cD_r(\hat\lambda,\hat\xi_-,\hat\eta_-)(z)
= r^{-1}z^{-1}\hat\xi_-(r z)
+ r^{-1}z\hat\eta_-(r z^{-1})
- \hat\lambda
$$
shows that
$$
\cD_r(\hat\lambda,\hat\xi_-,\hat\eta_-)=
\cD_1(\hat\lambda,(\hat\xi_-)_r,(\hat\eta_-)_r).
$$
Hence, by~(\ref{eq:rescale}) it suffices to prove the lemma for $r=1$.
Then
\begin{eqnarray*}
\cD_1(\hat\lambda,\hat\xi_-,\hat\eta_-)(z)
&= &
\sum_{n<0}\hat\xi_{n+1} z^n
- \hat\lambda + \sum_{n>0}\hat\eta_{1-n} z^n
\end{eqnarray*}
so
\begin{eqnarray*}
\bigl\|\cD_1(\hat\lambda,\hat\xi_-,\hat\eta_-)\bigr\|_s^2
&=&
\sum_{n<0}(1+|n|)^{2s}|\hat\xi_{n+1}|^2
+ |\hat\lambda|^2
+ \sum_{n>0}(1+|n|)^{2s}|\hat\eta_{1-n}|^2 \\
&=&
\sum_{n\le 0}(2+|n|)^{2s}|\hat\xi_n|^2
+ |\hat\lambda|^2
+ \sum_{n\le 0}(2+|n|)^{2s}|\hat\eta_n|^2 \\
&\ge&
\|\hat\xi_-\|_s^2  + |\hat\lambda|^2 +  \|\hat\eta_-\|_s^2.
\end{eqnarray*}
This proves the lemma.
\end{proof}

\begin{lemma}[Approximate Solution]\label{le:zero}
For every $s>1/2$ there is a constant $c>0$ such that
$$
      \left\|\cF_r(\xi_1\eta_1,\xi_+,\eta_+)\right\|_s
      \le cr\left(\left\|\xi_+-\id\right\|_s
      + \left\|\eta_+-\id\right\|_s\right)
$$
for every pair
$\xi_+,\eta_+\in H^s_+$ with $\left\|\xi_+\right\|_s\le 1$,
$\left\|\eta_+\right\|_s\le 1$, and every $r\in(0,1]$.
\end{lemma}

\begin{proof}
The constant is  $c=4\sqrt3C$
where $C$ is the constant of Lemma~\ref{le:prod}.
We first prove the inequality
\begin{equation}\label{eq:ap-approx}
\left\|\cF_1(\xi_1\eta_1,\xi_+,\eta_+)\right\|_s
\le 2\sqrt{3}C\left(\left\|\xi_+-\xi_1\id\right\|_s
      + \left\|\eta_+-\eta_1\id\right\|_s\right).
\end{equation}
Since
\begin{eqnarray*}
\cF_1(\xi_1\eta_1,\xi_+,\eta_+)(z)
&= &
\xi_+(z)\eta_+(z^{-1})-\xi_1\eta_1 \\
&= &
\sum_{k\ne 0}\left(\sum_{n>0}\xi_{n+k}\eta_n\right)z^k
+ \sum_{n>1}\xi_n\eta_n
\end{eqnarray*}
we have
\begin{eqnarray*}
      \left\|\cF_1(\xi_1\eta_1,\xi_+,\eta_+)\right\|_s^2
&=&
      \sum_{k\ne 0}(1+|k|)^{2s}
      \left|\sum_{n>0}\xi_{n+k}\eta_n\right|^2
      + \left|\sum_{n>1}\xi_n\eta_n\right|^2 \\
&\le&
      \sum_{k>0}(1+k)^{2s}
      \left(\sum_{n>0}\left|\xi_{n+k}\eta_n\right|\right)^2
      + \left(\sum_{n>1}\left|\xi_n\eta_n\right|\right)^2 \\
&&
      +\, \sum_{k>0}(1+k)^{2s}
        \left(\sum_{n>k}\left|\xi_{n-k}\eta_n\right|\right)^2 \\
&\le&
     3 C^2
      \left(\left\|\xi_+-\xi_1\id\right\|_s^2\left\|\eta_+\right\|_s^2
      + \left\|\xi_+\right\|_s^2\left\|\eta_+-\eta_1\id\right\|_s^2\right).
\end{eqnarray*}
The last inequality follows from Lemma~\ref{le:prod};
note that each sum omits either  $\xi_1$ or $\eta_1$ or both.
The inequality~(\ref{eq:ap-approx}) follows by taking the square root
of the last estimate and using the fact that
$\left\|\xi_+\right\|_s\le2$ and $\left\|\eta_+\right\|_s\le2$.

The formula
$$
      \cF_r(\xi_1\eta_1,\xi_+,\eta_+)(z)=
      r^{-2}\xi_+(r z)\eta_+(r z^{-1})-\xi_1\eta_1
$$
shows that
$$
    \cF_r(\xi_1\eta_1,\xi_+,\eta_+)= \cF_1(\xi_1\eta_1,(\xi_+)_r,(\eta_+)_r).
$$
Note that the operation $\xi\mapsto\xi_r$ leaves the coefficient
$\xi_1$ unchanged. Hence, by~(\ref{eq:ap-approx}), we have
\begin{eqnarray*}
      \left\|\cF_r(\xi_1\eta_1,\xi_+,\eta_+)\right\|_s
&= &
      \left\|\cF_1(\xi_1\eta_1,(\xi_+)_r,(\eta_+)_r)\right\|_s  \\
&\le &
      2\sqrt{3}C\bigl(\left\|(\xi_+)_r-\xi_1\id\right\|_s
      + \left\|(\eta_+)_r-\eta_1\id\right\|_s\bigr)  \\
&= &
      2\sqrt{3}C\bigl(\left\|\xi_+-\xi_1\id\right\|_{r,s}
      + \left\|\eta_+-\eta_1\id\right\|_{r,s}\bigr)  \\
&\le &
      2\sqrt{3}Cr\bigl(\left\|\xi_+-\xi_1\id\right\|_s
      + \left\|\eta_+-\eta_1\id\right\|_s\bigr)  \\
&\le &
      4\sqrt{3}Cr\bigl(\left\|\xi_+-\id\right\|_s\
      + \left\|\eta_+-\id\right\|_s\bigr).
\end{eqnarray*}
This proves the lemma.
\end{proof}

\begin{lemma}[Quadratic Estimate]\label{le:quad}
For every $s>1/2$ there is a constant $c>0$ such that
$$
     \bigl\|\left(d\cF_r(\lambda,\xi,\eta)-\cD_r\right)
       (\hat\lambda,\hat\xi,\hat\eta)\bigr\|_s
     \le c\bigl(\|\eta-\id\|_{r,s}\|\hat\xi\|_{r,s}
          + \|\xi-\id\|_{r,s}\|\hat\eta\|_{r,s}
          \bigr)
$$
for all $\xi,\eta,\hat\xi,\hat\eta\in H^s_r$
and $\lambda,\hat\lambda\in\C$.
\end{lemma}

\begin{proof}
We have
$$
         d\cF_r(\lambda,\xi,\eta)(\hat\lambda,\hat\xi,\hat\eta)(z)
         = r^{-2}\hat\xi(r z)\eta(r z^{-1})
            + r^{-2}\xi(r z)\hat\eta(r z^{-1}) - \hat\lambda
$$
and hence
$$
       \left(d\cF_r (\lambda,\xi,\eta)-\cD_r \right)
       (\hat\lambda,\hat\xi,\hat\eta)(z)
       = r^{-2}\hat\xi(r z)(\eta-\id)(r z^{-1})
            + r^{-2}(\xi-\id)(r z)\hat\eta(r z^{-1}).
$$
So the result follows from Lemma~\ref{le:prod} with $c=C$.
\end{proof}

\begin{lemma}[Inverse Function Theorem]\label{le:ift}
Let $f:U\to V$ be a smooth map between Banach spaces
and $D:U\to V$ be a Banach space isomorphism.
Let $u_0\in U$ and suppose that there is a constant
$\rho>0$ such that
\begin{equation}\label{eq:ift1}
\left\|D^{-1}\right\|\le 1,\qquad \left\|f(u_0)\right\|_V\le\frac{\rho}{2}
\end{equation}
and, for every $u\in U$,
\begin{equation}\label{eq:ift2}
\left\|u-u_0\right\|_U\le\rho\qquad\implies\qquad
\left\|df(u)-D\right\|\le\frac{1}{2}.
\end{equation}
Then there is a unique element $u\in U$ such that
$$
\left\|u-u_0\right\|_U\le\rho,\qquad f(u)=0.
$$
Moreover, $\left\|u-u_0\right\|_U\le 2\left\|f(u_0)\right\|_V$.
\end{lemma}

\begin{proof}  
Standard.
\end{proof}

\begin{proof}[Proof of Proposition~\ref{prop:localmodel1}.]
Throughout we fix a constant $s>1/2$ and a constant $c\ge1$ such that the
assertions of Lemmata~\ref{le:s}, \ref{le:zero} and~\ref{le:quad} hold with
these constants $s$ and $c$. Choose positive constants $\eps$, $\rho$,
and $\delta$ such that
\begin{equation}\label{eq:delta}
3c\eps < \frac12,\qquad
c\sqrt{2\delta^2+\rho^2}\le \frac{1}{2},\qquad
2c\delta\le\frac{\eps}{2},\qquad \rho:=\sqrt{3}\eps.
\end{equation}
We prove the assertion with these constants $c$ and $\delta$.

Assume first that $a$ is a positive real
number and denote $r:=\sqrt{a}$. Fix a pair
$(\xi_+,\eta_+)\in H^s_+\times H^s_+$
satisfying~(\ref{eq:deltaPlus}).
Let
$$
U:=\C\times H^s_{r,-}\times H^s_{r,-},\qquad V:=H^s
$$
and consider the map $f:U\to V$ defined by
\begin{equation}\label{eq:f}
f(u) := \cF_r(\lambda,\xi_++\xi_-,\eta_++\eta_-),
\qquad u:=(\lambda,\xi_-,\eta_-).
\end{equation}
Let $D:=\cD_r:U\to V$ and $u_0:=(\xi_1\eta_1,0,0)$.
Then, by Lemma~\ref{le:inverse},
\begin{equation}\label{eq:D-1}
\left\|D^{-1}\right\|_{\cL(V,U)} \le 1
\end{equation}
and, by Lemma~\ref{le:zero} and~(\ref{eq:delta}),
\begin{equation}\label{eq:u0}
\left\|f(u_0)\right\|_V
\le cr\left(\left\|\xi_+-\id\right\|_s+\left\|\eta_+-\id\right\|_s\right)
\le 2c\delta
\le\frac{\eps}{2}
\le\frac{\rho}{2}.
\end{equation}
In this notation  the operator $df(u)-D:U\to V$
is the restriction of the operator
$d\cF_r(\lambda,\xi,\eta)-\cD_r:\C\times H^s_r\times H^s_r\to H^s$
to the subspace $U$, so by Lemma~\ref{le:quad} we have
$$
\left\|df(u)-D\right\|_{\cL(U,V)}
\le c\sqrt{\|\eta-\id\|_{r,s}^2 + \|\xi-\id\|_{r,s}^2}
$$
for $u=(\lambda,\xi_-,\eta_-)\in U$
and $\xi:=\xi_++\xi_-$ and $\eta:=\eta_++\eta_-$.
Note that $\left\|\zeta\right\|_{r,s}\le\left\|\zeta\right\|_s$
for $\zeta\in H^s_+$ and $0<r\le1$.
Hence
$$
      \left\|df(u)-D\right\|_{\cL(U,V)}
      \le c\sqrt{\|\xi_+-\id\|_s^2 + \|\eta_+-\id\|_s^2
      + \|\xi_-\|_{r,s}^2 + \|\eta_-\|_{r,s}^2}
$$
for $u=(\lambda,\xi_-,\eta_-)\in U$.
Since $\|\xi_+-\id\|_s\le\delta$ and $\|\eta_+-\id\|_s\le\delta$
we have
\begin{equation}\label{eq:dfuD}
\left\|u-u_0\right\|_U\le\rho
\qquad\implies\qquad
\left\|df(u)-D\right\|_{\cL(U,V)}
\le c\sqrt{2\delta^2+\rho^2}
\le \frac12
\end{equation}
for every $u:=(\lambda,\xi_-,\eta_-)\in U$.
Here we have used~(\ref{eq:delta}).

It follows from~(\ref{eq:D-1}), (\ref{eq:u0}), and~(\ref{eq:dfuD})
that the assumptions of Lemma~\ref{le:ift} are satisfied.
Hence there is a unique point $u\in U$ such that
$$
\left\|u-u_0\right\|_U\le \rho,\qquad f(u)=0,
$$
and this unique point satisfies
\begin{equation}\label{eq:uu0}
\left\|u-u_0\right\|_U\le2\left\|f(u_0)\right\|_V\le\eps.
\end{equation}
Thus, for every $(\xi_+,\eta_+)\in H^s_+\times H^s_+ $
satisfying~(\ref{eq:deltaPlus}), we have found a
unique triple $(\lambda,\xi_-,\eta_-)\in U$
such that $\xi:=\xi_++\xi_-$ and $\eta:=\eta_++\eta_-$ satisfy
$$
\cF_r(\lambda,\xi,\eta)=0,\qquad
\|\xi_-\|_{r,s}\le\eps,\qquad
\|\eta_-\|_{r,s}\le\eps,\qquad
|\lambda-\xi_1\eta_1|\le\eps.
$$
That the quadruple $(a,\xi,\eta,b)$ also satisfies the
estimate~(\ref{eq:continuous}) follows from~(\ref{eq:u0})
and~(\ref{eq:uu0}).

Next we prove that this quadruple $(a,\xi,\eta,b)$ satisfies
$\xi(z)\ne 0$ and $\eta(z)\ne 0$ for $r\le|z|\le1$.
To see this note that
$$
\left\|(\zeta/\id)_+\right\|_s
=\sqrt{\sum_{n\ge 2}n^{2s}\left|\zeta_n\right|^2}
\le\left\|\zeta_+\right\|_s
$$
and
$$
r\left\|(\zeta/\id)_-\right\|_{r,s}
= \sqrt{\left|\zeta_1\right|^2
+ \sum_{n\le0}(2-n)^{2s}r^{2n-2}\left|\zeta_n\right|^2}
\le \left\|\zeta_+\right\|_s + 2\left\|\zeta_-\right\|_{r,s}.
$$
Hence
\begin{eqnarray*}
\sup_{r\le|x|\le 1}\left|\xi(x)x^{-1}-1\right|
&\le &
c\left(\left\|(\xi/\id-1)_+\right\|_s
+ r\left\|(\xi/\id-1)_-\right\|_{r,s}\right) \\
&\le &
2c\left(\left\|\xi_+-\id\right\|_s
+ \left\|\xi_-\right\|_{r,s}\right) \\
&\le&
2c(\delta+\eps) \\
&\le &
1/2.
\end{eqnarray*}
Here the first inequality follows from Lemma~\ref{le:s}
and the last uses the fact that $2c\eps\le1/3$
and $2c\delta\le\eps/2\le1/6$. Thus we have proved that
$\xi$ and $\eta$ do not vanish on the closed annulus
$r\le|z|\le 1$. Now extend $\xi$ and $\eta$ to the annulus
$r^2\le|z|\le1$ by the formulas
$$
\xi(x) := \frac{b}{\eta(ax^{-1})},\qquad
\eta(y) := \frac{b}{\xi(ay^{-1})},\qquad
r^2\le|x|,|y|\le r.
$$
The resulting functions $\xi$ and $\eta$ are continuous across
the circle of of radius $r$ by~(\ref{eq:xi-eta-a-b}).  Hence
they are holomorphic on the large annulus $r^2<|z|< 1$.
Since~(\ref{eq:xi-eta-a-b}) holds on the middle circle
$|x|=r$ it holds on the annulus $r^2\le|x|\le1$.
This proves the proposition for positive real numbers $a$.

To prove the proposition for general $a$ we use the
following ``rotation trick''. Fix a constant $\theta\in\R$.
Given $\xi,\eta\in H^s_+\oplus H^s_{r,-}$ and $a,b\in\C$
define $\tilde\xi,\tilde\eta\in H^s_+\oplus H^s_{r,-}$ and
$\tilde a,\tilde b\in\C$ by
$$
\tilde\xi(z) := e^{-i\theta}\xi(e^{i\theta}z),\quad
\tilde\eta(z) := e^{-i\theta}\eta(e^{i\theta}z),\quad
\tilde a := e^{-2i\theta}a,\quad
\tilde b := e^{-2i\theta}b.
$$
Then $a,b,\xi,\eta$ satisfy~(\ref{eq:xi-eta-a-b})
if and only if $\tilde a,\tilde b,\tilde\xi,\tilde\eta$
satisfy~(\ref{eq:xi-eta-a-b}).  Hence the result for
general $a$ can be reduced to the special case
case by choosing $\theta$ such that
$\tilde a := e^{-2i\theta}a$ is a positive real number.
This proves the proposition.
\end{proof}

\begin{proof}[Proof of Proposition~\ref{prop:localmodel2}.]
The general case can be reduced to the case $a>0$ by the
rotation trick in the proof of Proposition~\ref{prop:localmodel1}.
Hence we assume $a=r^2$ and $r>0$.
Choose positive constants $c$ and $C$ such that the
assertions of Lemmata~\ref{le:s}, \ref{le:prod}, and~\ref{le:quad}
hold with these constants. Choose $\delta$ and $\eps$
such that
$$
2c\delta\le 1,\qquad 8C(1+c)\eps < 1.
$$
Let $(a,\xi,\eta,b)$ and $(a',\xi',\eta',b')$ satisfy the
assumptions of Proposition~\ref{prop:localmodel2} with these
constants $\delta$ and $\eps$
and denote
$$
\lambda:=b/a=b/r^2,\qquad \lambda':=b'/a'=b'/r^2.
$$
Then
\begin{equation}\label{eq:unique}
r^{-2}\xi(rz)\eta(rz^{-1})=\lambda,\qquad
r^{-2}\xi'(rz)\eta'(rz^{-1})=\lambda',\qquad
|z|=1.
\end{equation}
Denote by $L_r:\C\times H^s_{r,-}\times H^s_{r,-}\to H^s$
the linear operator given by
$$
L_r(\hat\lambda,\hat\xi_-,\hat\eta_-)(z)
:=
r^{-2}\hat\xi_-(rz)\eta_+(rz^{-1})
+ r^{-2}\xi_+(rz)\hat\eta_-(rz^{-1})
- \hat\lambda.
$$
In the notation of~\ref{cFr} the operator $L_r$ is the restriction
of the differential of $\cF_r$ at $(\lambda,\xi_+,\eta_+)$
(for any $\lambda$) to the subspace $\C\times H^s_{r,-}\times H^s_{r,-}$.
Since $c\delta\le1/2$ it follows from
Lemmata~\ref{le:inverse} and~\ref{le:quad}
that the operator $L_r$ is invertible and the norm
of the inverse is bounded by $2$:
\begin{equation}\label{eq:Linverse}
|\hat\lambda|^2 + \|\hat\xi_-\|_{r,s}^2 + \|\hat\eta_-\|_{r,s}^2
\le 4 \|L_r(\hat\lambda,\hat\xi_-,\hat\eta_-)\|_s^2.
\end{equation}
Let us denote by $Q_r:H^s_{r,-}\times H^s_{r,-}\to H^s$ the
quadratic form
$$
Q_r(\xi_-,\eta_-)(z) := r^{-2}\xi_-(rz)\eta_-(rz^{-1}).
$$
Then, by Lemma~\ref{le:prod}, we have
\begin{equation}\label{eq:Qquad}
\left\|Q_r(\xi_-,\eta_-)\right\|_s
\le Cr^{-1}\left(
\left\|\xi_-\right\|_{r,s}\sup_{|y|=r}\left|\eta_-(y)\right|
+ \left\|\eta_-\right\|_{r,s}\sup_{|x|=r}\left|\xi_-(x)\right|
\right).
\end{equation}
Now let
$$
\hat\lambda:=\lambda'-\lambda,\qquad
\hat\xi_-:=\xi'-\xi,\qquad\hat\eta_-:=\eta'-\eta.
$$
Then the difference of the two equations in~(\ref{eq:unique})
can be expressed in the form
\begin{eqnarray*}
L_r(\hat\lambda,\hat\xi_-,\hat\eta_-)
&= &
Q_r(\xi_-,\eta_-)-Q_r(\xi'_-,\eta'_-) \\
&= &
-Q_r(\hat\xi_-,\eta_-) - Q_r(\xi_-,\hat\eta_-) - Q_r(\hat\xi_-,\hat\eta_-).
\end{eqnarray*}
Abbreviate
$$
\hat\zeta:=(\hat\lambda,\hat\xi_-,\hat\eta_-),\qquad
\|\hat\zeta\|_{r,s}
:= \sqrt{|\hat\lambda|^2 + \|\hat\xi_-\|_{r,s}^2 + \|\hat\eta_-\|_{r,s}^2}.
$$
Then
\begin{eqnarray*}
\|\hat\zeta\|_{r,s}
&\le &
2\|L_r\hat\zeta\|_s \\
&\le &
2\left(\|Q_r(\hat\xi_-,\eta_-)\|_s
+\|Q_r(\xi_-,\hat\eta_-)\|_s
+\|Q_r(\hat\xi_-,\hat\eta_-)\|_s\right) \\
&\le &
2C\bigl(\|\xi_-\|_{r,s}+\|\eta_-\|_{r,s}\bigr)\|\hat\zeta\|_{r,s} \\
&&
+\, 2Cr^{-1}\left(
\sup_{|x|=r}|\hat\xi_-(x)|
+ \sup_{|y|=r}|\hat\eta_-(y)|
\right)
\|\hat\zeta\|_{r,s} \\
&\le &
2C(1+c)\bigl(\|\xi_-\|_{r,s}+\|\eta_-\|_{r,s}\bigr)\|\hat\zeta\|_{r,s} \\
&&
+\, 2Cr^{-1}\left(
\sup_{|x|=r}|\xi'_-(x)| + \sup_{|y|=r}|\eta'_-(y)|
\right)
\|\hat\zeta\|_{r,s} \\
&\le &
8C(1+c)\eps \|\hat\zeta\|_{r,s}.
\end{eqnarray*}
Here the first inequality follows from~(\ref{eq:Linverse}),
the second from the triangle inequality, the third from
Lemma~\ref{le:prod} and~(\ref{eq:Qquad}),
the fourth from Lemma~\ref{le:s}, and the last from the
assumptions of the proposition. Since
$8C(1+c)\eps<1$ it follows that $\hat\zeta=0$.
This proves the proposition.
\end{proof}

\begin{proof}[Proof of Proposition~\ref{prop:localmodel3}.]
{\bf Step~1.}
{\it  The map $\cT$ is continuous.
}

\medskip\noindent
Continuity for $a>0$ is an easy consequence of the
proof of Proposition~\ref{prop:localmodel1}.  The map
$f$ defined in equation~(\ref{eq:f}) depends continuously 
on the parameters $\xi_+,\eta_+,r$ and $r=\sqrt{a}$ 
depends continuously on $a$.  For complex nonzero 
$a$ we can choose $\theta$ in the rotation trick
to depend continuously on $a$.
To prove continuity for $a=0$ we deduce
from~(\ref{eq:continuous}) that
there is a constant $c>0$ such that
$$
\|\cT(a,\xi_+,\eta_+)\| := \sqrt{|b|^2+\|\xi_-\|_s^2 + \|\eta_-\|_s^2}
\le c|a|.
$$
Here we used the fact that $\|\zeta\|_s\le r\|\zeta\|_{r,s}$
for $\zeta\in H^s_{r,-}$.

\medskip\noindent{\bf Step~2.}
{\it Let $\xi_+,\eta_+\in H^s_+(\id,\delta)$ and $a\in\D\setminus0$ and
denote $\xi:=\xi_++\xi_-$ and $\eta:=\eta_++\eta_-$,
where $ (b,\xi_-,\eta_-):=\cT(a,\xi_+,\eta_+)$.
Then the linear operator
$$
    L:\C\times H^s_{r,-}\times H^s_{r,-}\to H^s_r
$$
defined by
$$
L(\hat{b},\hat{\xi}_-,\hat{\eta}_-)(x)
:= \hat{\xi}_-(x)\eta(ax^{-1})+\xi(x)\hat{\eta}(ax^{-1})-\hat{b}
$$
is invertible.}

\medskip\noindent
In the notation of the proof of Proposition~\ref{prop:localmodel1} we have
that $L$ is conjugate to the operator $df(u)$. Specifically,
$$
L(\hat{\lambda}r^2,\hat{\xi}_-,\hat{\eta}_-)_r
=rdf(u)(\hat{\lambda},\hat{\xi}_-,\hat{\eta}_-)
$$
when $a=r^2$. The operator $df(u)$ is
invertible by~(\ref{eq:D-1}) and~(\ref{eq:dfuD}).
For general $a$ use the rotation trick
from the end of the proof of Proposition~\ref{prop:localmodel1}.

\medskip\noindent{\bf Step~3.}
{\it  The map $\cT$ is continuously differentiable for $0<|a|<1$.}

\medskip\noindent
We formulate a related problem. Define a partial rescaling operator
$$
R_r:H^s\to H^s_+\oplus H^s_{r,-}, \qquad
(R_r\xi)(z):= \xi_+(z)+r\xi_-(r^{-1}z).
$$
The operator $R_r$ is a Hilbert space isometry for every $r\in(0,1]$.
Let
$$
\cX:=\C\times\C\times H^s\times H^s,\qquad
\cY:= H^s.
$$
There is a splitting  $\cX = \cX_+\oplus \cX_-$ where
$$
\cX_+:= \C\times H^s_+\times H^s_+, \qquad
\cX_-:= \C\times H^s_-\times H^s_-.
$$
Let $\cU\subset\cX$ denote the open set $\{0<|a|<1\}$ and
define $\tilde\cF:\cU\to\cY$ by
$$
\tilde\cF(a,b,\xi,\eta)(z):=(R_r\xi)(rz)\cdot
(R_r\eta)(ar^{-1}z^{-1})-b,\qquad r:=\sqrt{|a|}.
$$
Define $\tilde\cT:\cU\cap\cX_+\to\cX_-$ by
$$
\tilde\cT(a,\xi_+,\eta_+):=\bigl(b,(\xi_-)_r,(\eta_-)_r\bigr), \qquad
(b,\xi_-,\eta_-):=\cT(a,\xi_+.\eta_+).
$$
(Recall that $\zeta_r(z):=r^{-1}\zeta(rz)$
for $\zeta\in H^s_r$ and $|z|=1$.)
By construction the graph of $\tilde\cT$ is
contained in the zero set of $\tilde\cF$.
By step~2 the derivative of $\tilde\cF$ in
the direction $\tilde\cX_-$ is an invertible operator
at every point in the graph of $\tilde\cT$.
Hence, by the implicit function theorem,
$\tilde\cT$ is continuously differentiable.
Define the map $\cR:\cU\to\cU$ by
$$
      \cR(a,b,\xi,\eta) := (a,b,R_r\xi,R_r\eta),\qquad r:=\sqrt{|a|}.
$$
This map is continuously differentiable and
$$
      \mathrm{graph}(\cT) = \cR\circ\mathrm{graph}(\tilde\cT).
$$
Here $\mathrm{graph}(\cT)$ denotes the map
$(a,\xi_+,\eta_+)\mapsto(a,b,\xi,\eta)$
given by $(b,\xi_-,\eta_-):=\cT(a,\xi_+,\eta_+)$.
Similarly for $\mathrm{graph}(\tilde\cT)$.
Hence $\mathrm{graph}(\cT)$ is continuously
differentiable for $0<|a|<1$ and so is $\cT$.

\medskip\noindent{\bf Step~4.}
{\it  The map $\cT$ is holomorphic for $0<|a|<1$.}

\medskip\noindent
As $\cT$ is differentiable we have
$$
d\cT(a,\xi_+,\eta_+)(\hat{a},\hat\xi_+,\hat\eta_+)
= (\hat{b},\hat{\xi}_-,\hat{\eta}_-)
$$
for $\hat a\in\C$ and $\hat{\xi}_+,\hat{\eta}_+\in H^s_+$
where $\hat{\xi}_-,\hat{\eta}_-\in H^s_{r,-}$
and $\hat{b}\in\C$ are determined by the equation
\begin{equation}\label{eq:daT}
L(\hat{b},\hat\xi_-,\hat\eta_-)(x)
= -\hat\xi_+(x)\eta(ax^{-1}) - \xi(x)\hat\eta_+(ax^{-1})
-\xi(x)\eta'(ax^{-1})\hat{a}x^{-1}
\end{equation}
for $|x|=r:=\sqrt{|a|}$.
Here $L:\C\times H^s_{r,-}\times H^s_{r,-}\to H^s$ is the
operator of Step~2. Since $L$ is complex linear so is
$d\cT(a,\xi_+,\eta_+)$.

\medskip\noindent{\bf Step~5.}
{\it  The map $\cT$ is holomorphic for $|a|<1$.}

\medskip\noindent
That $\cT$ is holomorphic near $a=0$ follows from Step~4,
continuity, and the Cauchy integral formula.
More precisely, suppose $X$ and $Y$ are complex Hilbert
spaces and $\cT:\C\times X\to Y$ is a continuous map
which is holomorphic on $(\C\setminus 0)\times X$.
Then
$$
\cT(a,x) = \frac{1}{2\pi}\int_0^{2\pi}
\cT(a+e^{i\theta}\hat a,x+e^{i\theta}\hat x)\,d\theta
$$
and
$$
d\cT(a,x)(\hat a,\hat x)
= \frac{1}{2\pi}\int_0^{2\pi}
e^{-i\theta}\cT(a+e^{i\theta}\hat a,x+e^{i\theta}\hat x)\,d\theta
$$
for $a,\hat a\in\C$ and $x,\hat x\in X$ with $a\ne 0$.
In the case at hand $\cT(a,x)$ converges uniformly to
zero as $|a|$ tends to zero (see the proof of Step~1).
By the Cauchy integral formula, this implies that $d\cT(a,x)$
converges uniformly to zero in the operator norm as
$|a|$ tends to zero. This proves the proposition.
\end{proof}

\begin{proof}[Proof of Theorem~\ref{thm:localmodel}.]
Fix a positive integer $s$ and
choose $\delta$, $c$, and $\eps$
such that Propositions~\ref{prop:localmodel1}
and~\ref{prop:localmodel2} hold. Shrinking
$\delta$ we may assume $4c\delta<\eps$.
We prove that the graph of $\cT$ intersects $\cW_\delta$
in $\cU_\delta$. By definition
$$
\mathrm{graph}(\cT)\cap\cW_\delta\subset\cU_\delta.
$$
To prove the converse choose $(a,\xi,\eta,b)\in\cU_\delta$.
If $a=0$ then $\xi_-=\eta_-=0$ and $b=0$
so $(a,\xi,\eta,b)$ belongs to the graph
of $\cT$. Hence assume $a\ne 0$ and let $r:=\sqrt{|a|}$.
Then $\left\|\xi_+-\id\right\|_s<\delta$
and $\left\|\eta_+-\id\right\|_s<\delta$.
So, by Proposition~\ref{prop:localmodel1},
there is an element $(a,\xi',\eta',b')
\in\cW_\delta\cap\mathrm{graph}(\cT)$
satisfying $\xi'_+=\xi_+$, $\eta'_+=\eta_+$, and
$$
\left\|\xi'_-\right\|_{r,s}\le 4cr\delta<\eps,\qquad
\left\|\eta'_-\right\|_{r,s}\le 4cr\delta<\eps.
$$
We claim that $\xi=\xi'$, $\eta=\eta'$, and $b=b'$.
By Proposition~\ref{prop:localmodel2}
it suffices to show that
$$
\sup_{|x|=r}\left|\xi_-(x)\right|<r\eps,\qquad
\sup_{|y|=r}\left|\eta_-(y)\right|<r\eps.
$$
By symmetry we need only prove the first inequality.
By the triangle inequality
\begin{equation}\label{eq:tri}
\sup_{|x|=r}\left|\xi_-(x)\right|\le
\sup_{|x|=r}\left|\xi(x)-x\right| +
\sup_{|x|=r}\left|\xi_+(x)-x\right|.
\end{equation}
By Lemma~\ref{le:apriori} we estimate the first term on the right by
\begin{equation}\label{eq:firstTerm}
\sup_{|x|=r}\left|\xi(x)-x\right|\le cr\delta.
\end{equation}
For the second term we have by Lemma~\ref{le:s}
$$
\sup_{|x|=r}\left|\xi_+(x)-x\right|=
r\sup_{|z|=1}\left|(\xi_+-\id)_r(z)\right| \le
rc\left\|(\xi_+-\id)_r\right\|_s=cr\left\|\xi_+-\id\right\|_{r,s}
$$
But the series for $\xi_+-\id$ has only positive powers
and $r\le 1$ so
$$
\left\|\xi_+-\id\right\|_{r,s}\le \left\|\xi_+-\id\right\|_s
\le \left\|\xi-\id\right\|_s\le\delta.
$$
Combining the last two lines gives
\begin{equation}\label{eq:secondTerm}
\sup_{|x|=r}\left|\xi_+(x)-x\right|\le cr\delta.
\end{equation}
Now use~(\ref{eq:tri}), (\ref{eq:firstTerm}), and~(\ref{eq:secondTerm})
and shrink $\delta$ so $2c\delta<\eps$.
\end{proof}

We close this section with two lemmas that will be useful in the sequel.

\begin{lemma}\label{le:smooth}
Fix $s>1/2$ and choose $\delta>0$ as in
Theorem~\ref{thm:localmodel}.
Let $A\subset\INT(\D)\times\C^m$
be an open set and
$$
A\to\cU_\delta:(a,t)\mapsto(a,\xi_{a,t},\eta_{a,t},b_{a,t})
$$
be a holomorphic map. Then the map
$$
\Phi:\{(x,y,t)\in\C^{2+m}\,:\,x,y\in\INT(\D),\,(xy,t)\in A\}
\to\C\times\C
$$
given by
$$
\Phi(x,y,t):=\Phi_t(x,y):=(\xi_{xy,t}(x),\eta_{xy,t}(y))
$$
is holomorphic.
\end{lemma}

\begin{proof}
The evaluation map
$$
H^s\cap H^s_{r^2}\times
\left\{z\in\C\,:\,r^2<|z|<1\right\}\to\C:
(\zeta,z)\mapsto\zeta(z)
$$
is holomorphic.
It follows that the map $(x,y,t)\mapsto\Phi_t(x,y)$ is
holomorphic in the domain $xy\ne 0$.
We prove that $\Phi$ is continuous.
Suppose $x_i\to x\ne 0$, $y_i\to 0$, and $t_i\to t$.
Then $\xi_{x_iy_i,t_i}$ converges to $\xi_{0,t}$,
uniformly in a neighbourhood of $x$, and hence
$\xi_{x_iy_i,t_i}(x_i)$ converges to $\xi_{0,t}(x)$.
Moreover, if $c$ and $\delta$ are the constants of
Lemma~\ref{le:apriori}, then
$$
|\eta_{x_iy_i,t_i}(y_i)|\le(c\delta+1)|y_i|
$$
and so $\eta_{x_iy_i,t_i}(y_i)$ converges to $\eta_{0,t}(0)=0$.
Hence $\Phi_{t_i}(x_i,y_i)$ converges to
$\Phi_t(x,0)=(\xi_{0,t}(x),0)$.  Hence $\Phi$ is
continuous at every point $(x,0,t)$ with $x\ne 0$.
By symmetry, $\Phi$ is continuous at every point
$(0,y,t)$ with $y\ne 0$.  That $\Phi$ is continuous
at very point $(0,0,t)$ follows again from
Lemma~\ref{le:apriori}. Since $\Phi$ is continuous
and is holomorphic in $xy\ne 0$, it follows from the
Cauchy integral formula that $\Phi$ is holmorphic.
\end{proof}

\begin{lemma} \label{le:nodalExtend}
Let  $\xi_0,\eta_0:\INT(\D)\to\C$ be two holomorphic functions
satisfying $\xi_0(0)=\eta_0(0)=0$ and $\xi_0'(0)\ne0$, $\eta_0'(0)\ne 0$.
Then there   are neighborhoods $U_1$ and $U_2$ of $(0,0)$ in $\C^2$
and $B_1$ and $B_2$ of $0$ in $\C$ and
holomorphic diffeomorphisms $\Phi:=(\xi,\eta):U_1\to U_2$
and $\zeta:B_1\to B_2$
such that
$$
     \xi(x,0)=\xi_0(x), \qquad \eta(0,y)=\eta_0(y), \qquad \xi(x,y)\eta(x,y)=\zeta(xy)
$$
for $x,y$ near $0$.
\end{lemma}

\begin{proof} Replacing $\xi_0$ and $\eta_0$
by $\xi_0'(0)^{-1}\xi_0$ and $\eta_0'(0)^{-1}\eta_0$
we may assume w.l.o.g. that $\xi_0'(0)=\eta_0'(0)=1$.
Replacing $\xi_0(x)$ and $\eta_0(y)$ by  $\eps^{-1}\xi_0(\eps x)$ and
$\eps^{-1}\eta_0(\eps y)$ we may assume w.l.o.g. that the power series
for $\xi_0$ and $\eta_0$ lie in $H^s_+(\id,\delta)$ with $\delta>0$ as in~\ref{cT}.
For $z\in\D$ define $\alpha_z,\beta_z\in H^s_-$ by
$
(\zeta(z),\alpha_z,\beta_z):=\cT(z,\xi_0,\eta_0)
$
and then define
$\xi(x,y):=\xi_0(x)+\alpha_{xy}(x)$  and $\eta(x,y):=\eta_0(y)+\beta_{xy}(y)$.
Then $\Phi$ is holomorphic by Lemma~\ref{le:smooth}.
A direct calculation shows that $d\Phi(0,0)$ is the identity
so $\Phi$ is a local diffeomorphism.
The desired identities follow from the definition of $\cT$.
To prove that $\zeta'(0)=1$ differentiate the identity $\xi\eta=\zeta$ twice.
\end{proof}


\section{Hardy decompositions}\label{sec:hardy}

In this section we redo Section~\ref{sec:hardy1} in parametrized form. 
Here is where we use the local model of Section~\ref{sec:local}.

\para Throughout this section
$(\pi_A:P\to A,R_*,a_0)$ and $(\pi_B:Q\to B,S_*,b_0)$  are nodal unfoldings,
$$
f_0: P_{a_0}\to Q_{b_0}
$$
is a fiber isomorphism, and $p_1,p_2,\ldots,p_k$ are
the nodal points of the central fiber $P_{a_0}$,
so $q_i:=f_0(p_i)$ (for $i=1,\ldots,k$)
are the nodal points of the central fiber $Q_{b_0}$.
Let $m:=\dim_\C(A)$
and  $d:=\dim_\C(B)$.
Let $C_A\subset P$ and $C_B\subset Q$ denote the
critical points of   $\pi_A$ and $\pi_B$ respectively.
\arap

\dfn\label{def:hardyDecomposition}
A \jdef{Hardy decomposition} for $(\pi_A,R_*,a_0)$
is a decomposition
$$
P=M\cup N, \qquad \p M=\p N = M\cap N,
$$
into manifolds with boundary such that
$$
N=N_1\cup\cdots \cup N_k,
$$
each $N_i$ is a neighborhood of $p_i$,
the closures of the  $N_i$ are pairwise disjoint,
$N$ is disjoint from the elements of $R_*$,
and $N$  is the domain of a \jdef{nodal coordinate system}.
This consists of three sequences of holomorphic maps
$$
(x_i,y_i):N_i\to\D^2,\qquad z_i: A \to\C,\qquad t_i:A\to\C^{m-1},
$$
such that each map
$$
A\to   \D\times\C^{m-1}:a\mapsto(z_i(a),t_i(a))
$$
is a holomorphic coordinate system,
each  map
$$
N_i\to \D^2\times \C^{m-1}:
p\mapsto \bigl(x_i(p),y_i(p),t_i(\pi_A(p))\bigr)
$$
is a holomorphic coordinate system and
$$
x_i(p_i)=y_i(p_i)=0,\qquad    z_i\circ\pi_A =x_iy_i.
$$
(Note that here $N_i$ has a boundary whereas its analog
$U_i$ in~\ref{UV} was open.)
Restricting to a fiber gives a decomposition
$$
P_a=M_a\cup N_a, \qquad
M_a:=M\cap P_a, \qquad
N_a:=N\cap P_a
$$
where $M_a$ is a Riemann surface with boundary
and each component of $N_a$ is either
a closed annulus or a pair of transverse closed disks.
The nodal coordinate system determines a trivialization
\begin{equation}\label{eq:trivialize}
\iota:A\times\Gamma\to\p  N,\qquad
\Gamma:=\bigcup_{i=1}^k\{(i,1),(i,2)\}\times S^1,
\end{equation}
where $\iota^{-1}$ is  the disjoint union of the maps
$$
\pi\times x_i:\p_1N_i\to A\times S^1,\qquad \p_1 N_i:=\{|x_i|=1\},
$$
$$
\pi\times y_i:\p_2 N_i\to A\times S^1,\qquad \p_2 N_i:=\{|y_i|=1\}.
$$
The indexing is so that
$\iota((i,1)\times S^1)=x_i(\p_1N_i)$ and $\iota((i,2)\times S^1)=y_i(\p_2N_i)$.
For $a\in A$ define $\iota_a:\Gamma\to \p N$ by
$\iota_a(\lambda):=\iota(a,\lambda)$.
\nfd

\begin{lemma}\label{le:f=id}
After shrinking $A$ and $B$ if necessary,
there is a Hardy decomposition $P=M\cup N$ as
in~\ref{def:hardyDecomposition} and there are open
subsets $U=U_1\cup\cdots\cup U_k$, $V$, $W$ of $Q$
and functions $\xi_i,\eta_i,\zeta_i,\tau_i$ as described in~\ref{UV}
such that
$$
f_0(M_{a_0})\subset V_{b_0},\qquad f_0(N_{a_0})\subset U_{b_0}.
$$
and
$$
\xi_i\circ f_0\circ x_i^{-1}(x,0,0)=x,\qquad
\eta_i\circ f_0\circ y_i^{-1}(0,y,0)=y
$$
for $x,y\in\D$.
\end{lemma}

\begin{proof}
Choose any Hardy decomposition
$P=M\cup N$ as in~\ref{def:hardyDecomposition}
as well as open subsets $U=U_1\cup\cdots\cup U_k$, $V$, $W$ of $Q$
and functions $\xi_i,\eta_i,\zeta_i,\tau_i$ as described in~\ref{UV}.
Read $\xi_i\circ f_0\circ x_i^{-1}(x)$ for $\xi_0(x)$ and
$\eta_i\circ f_0\circ y_i^{-1}(y)$ for $\eta_0(y)$
in Lemma~\ref{le:nodalExtend}, let $\Phi$ and $\zeta$
be as in the conclusion of that Lemma, and replace
$(\xi_i,\eta_i)$ by $\Phi^{-1}\circ(\xi_i,\eta_i)$
and $\zeta_i$ by $\zeta^{-1}\circ\zeta_i$.
This requires shrinking $U_i$ (and $B$).
Then shrink $N$ so that $f_0(N_{a_0})\subset U_{b_0}$
and enlarge $V$ so that $f_0(M_{a_0})\subset V_{b_0}$.
\end{proof}

\para\label{cWcUcV}
We use a Hardy decomposition
to mimic the construction of~\ref{WUV0}
with $a\in A$ as a parameter.
Choose a Hardy decomposition
$P=M\cup N$ for $(\pi_A,R_*,a_0)$,
open subsets $U=U_1\cup\cdots\cup U_k$, $V$, $W$ of $Q$,
and functions $\xi_i,\eta_i,\zeta_i,\tau_i$ as described in~\ref{UV},
such that the conditions of Lemma~\ref{le:f=id} are satisfied.
Fix a constant $s>1/2$ and define an  open subset
$$
\cW(a,b)\subset H^s(\p N_a,W_b)
$$
by the condition that for $\gamma\in H^s(\p N_a,W_b)$
we have $\gamma\in\cW(a,b)$ iff
$$
\gamma\bigl(x_i^{-1}(S^1)\bigr)\subset W_{i,1},\qquad
\gamma\bigl(y_i^{-1}(S^1)\bigr)\subset W_{i,2},
$$
(see~\ref{UV} for the notation $W_{i,1}$ and $W_{i,2}$)
and the curves $\xi_i\circ \gamma\circ x_i^{-1}$ and
$\eta_i\circ \gamma\circ y_i^{-1}$ from $S^1$ to $\C\setminus0$
both have winding number one about the origin.
For $a\in A$ and $b\in B$ let
\begin{equation*}
\begin{split}
\cV(a,b)&:=\{\gamma=v|\p N_a\in\cW(a,b):
v\in \Hol^s(M_a,V_b),\,
v(R_*\cap P_a)=S_*\cap Q_b\},\\
\cU(a,b)&:=\{\gamma=u|\p N_a\in\cW(a,b):
u\in \Hol^s(N_a,U_b),\,u(C_A\cap P_a)=C_B\cap Q_b\}.
\end{split}
\end{equation*}
Here $\Hol^s(X,Y)$ denotes the set of continuous maps
from $X$ to $Y$ which are holomorphic
on the interior of $X$ and restrict to $H^s$-maps
on the boundary; holomorphicity at a nodal point
is defined as in \ref{standardNode}. Define
$$
\cW_a:=\bigsqcup_{b\in B}\cW(a,b),\qquad
\cV_a:=\bigsqcup_{b\in B}\cV(a,b),\;\qquad
\cU_a:=\bigsqcup_{b\in B}\cU(a,b),
$$
$$
\cW:=\bigsqcup_{a\in A}\cW_a,\qquad
\cV:=\bigsqcup_{a\in A}\cV_a,\qquad
\cU:=\bigsqcup_{a\in A}\cU_a.
$$
Our notation means that the three formulas
$(a,\gamma,b)\in\cW$, $(\gamma,b)\in\cW_a$,
and $\gamma\in\cW(a,b)$ have the same meaning.
\arap

We will use the implicit function  theorem on a manifold of maps.
The main difficulty in defining a suitable manifold of maps
is that the nodal family $\pi_A$ is not  locally trivial as the homotopy
type of the fiber changes. To circumvent this difficulty we
need the   following

\dfn\label{def:hardyTrivialization}
A \jdef{Hardy trivialization} for  $(\pi_A:P\to A,R_*,a_0)$ is a triple
$
(M\cup N,\iota,\rho)
$
where $P=M\cup N$ is a Hardy decomposition with
corresponding trivialization
$
\iota:A\times\Gamma\to \p N
$
as in~\ref{def:hardyDecomposition}
and
$$
\rho:M\to\Omega:=M_{a_0}
$$
is a trivialization
such that $\iota^{-1}$ agrees with  $\pi_A\times \p\rho$, or equivalently
$$
        \rho_a\circ\iota_a=\iota_{a_0}
$$
for $a\in A$.  The trivialization is a smooth map
$
\rho:M\to \Omega
$
such that $\rho_a:=\rho|M_a:M_a\to \Omega$ is a diffeomorphism
satisfying
$$
\rho(R_*)=R_*\cap \Omega=:r_*
$$
for every $a\in A$.
The map $\rho$ determines a trivialization
$$
    \p\rho:=\rho|\p M: \p M\to  \p  \Omega
$$
of the bundle $\pi|\p M:\p M\to A$.
Note that $\iota_{a_0}:\Gamma\to\p\Omega$ identifies
the boundary of $\Omega$ with a disjoint union of circles.
\nfd

\para\label{cW}
Fix a Hardy trivialization $(P=M\cup N,\iota,\rho)$;
we use it construct an auxiliary Hilbert manifold structure on $\cW$.
The domains of the maps in this space vary with $a$
so we replace them with a constant domain by using
an appropriate trivialization.   Define an open set
$$
\cW_0\subset \left\{(a,\gamma,b)\in A\times H^s(\Gamma,W)\times B:
\pi_B\circ\gamma=b\right\}
$$
by the condition that the map
$$
\cW_0\to\cW:(a,\gamma,b)\mapsto(a,\gamma\circ\rho_a|\p N_a,b)
$$
is a bijection.  In particular $\gamma((i,1)\times S^1)\subset W_{i,1}$
and $\gamma((i,2)\times S^1)\subset W_{i,2}$ for $(a,\gamma,b)\in\cW_0$.
(By a standard construction $H^s(\Gamma,W)$ is a complex
Hilbert manifold and the subset $\{(a,\gamma,b):\pi_B\circ\gamma=b\}$
is a complex Hilbert submanifold of $A\times H^s(\Gamma,W)\times B$.
This is because the map $H^s(\Gamma,W)\to H^s(\Gamma,B)$
induced by $\pi_B$ is a holomorphic submersion.
Note that $\cW_0$ is a connected component
of $\{(a,\gamma,b):\pi_B\circ\gamma=b\}$ and hence inherits
its Hilbert manifold structure.)
We emphasize that the resulting Hilbert manifold structure
on $\cW$ depends on the choice of the Hardy trivialization.
Two different Hardy trivializations give rise to a homeomorphism
which is of class $C^k$ on the dense subset $\cW\cap H^{s+k}$.
\arap

\para\label{shrinking}
The fiber isomorphism $f_0:P_{a_0}\to Q_{b_0}$ determines
a point
$$
(a_0,\gamma_0:=f_0|\p N_{a_0},b_0)\in\cW;
$$
this point lies in $\cU\cap\cV$ as
$$
\gamma_0=u_0|\p N_{a_0}=v_0|\p M_{a_0},
\qquad\mbox{where}   \qquad
u_0:=f_0|  N_{a_0},\quad v_0:=f_0|M_{a_0}.
$$
In the sequel we will denote neighborhoods of $a_0$ in $A$ and
$(a_0,\gamma_0,b_0)$ in $\cU$, $\cV$, or $\cW$ by the same letters
$A$, respectively $\cU$,  $\cV$, or $\cW$, and signal this with
the text ``shrinking $A$, $\cU$,  $\cV$, or $\cW$, if necessary''.
\arap

\begin{lemma}\label{le:fgamma}
For every $(a,\gamma,b)\in\cU\cap\cV$
there is a unique fiber isomorphism
$f:P_a\to Q_b$ with $f|\p N_a=\gamma$.
\end{lemma}

\begin{proof}
This follows immediately from Lemma~\ref{le:ugamma}.
\end{proof}

\begin{theorem}\label{thm:UV}
Assume $s>7/2$.
After shrinking $A$, $\cU$, $\cV$, $\cW$,
if necessary, the following holds.
\begin{description}
\item[(i)]
For each $a\in A$,
$\cU_a$ and $\cV_a$ are complex submanifolds of~$\cW_a$.
\item[(ii)]
$\cU$ is a  complex submanifold of $\cW$ and
$\cV$ is a smooth submanifold of~$\cW$.
\item[(iii)]
The projections $\cW\to A$ and $\cU\to A$ are holomorphic
submersions and the projection $\cV\to A$ is a smooth submersion.
\item[(iv)]
The  unfolding $(\pi_B,S_*,b_0)$ is infinitesimally
universal if and only if
$$
T_{w_0}\cW_{a_0}=T_{w_0}\cU_{a_0}\oplus
T_{w_0}\cV_{a_0},\qquad
w_0=(a_0,\gamma_0,b_0).
$$
\end{description}
\end{theorem}

\begin{proof}
In~\ref{UV} it was not assumed that $k$ was precisely the number of nodal pairs
so the arguments of Section~\ref{sec:hardy1} will apply when the central
fiber $Q_{b_0}$ is replaced by a nearby fiber $Q_b$ with possibly fewer
nodal points. Hence~(i) and~(iv) follow from Theorem~\ref{thm:infuniv}.

\bigskip

{\em We prove that $\cU$ is a complex Hilbert submanifold of $\cW$.}

\medskipnoindent The image of the nodal coordinate system
$(x_i,y_i,t_i)$ on $N_i$ in $\C\times\C\times\C^{m-1}$ has the form
$$
\{(x,y,t)\in\D^2\times\C^{m-1}: (xy,t)\in A_i,\; |x|<1,\; |y|<1\}.
$$
where $A_i\subset\C\times\C^{m-1}$ is contained in
the open set $\{|z_i|<1\}\times\C^{m-1}$.
The image of the nodal coordinate systems
$(\xi_i,\eta_i,\tau_i)$ on $U_i$ in $\C\times\C\times\C^{d-1}$
has the form
$$
     \left\{(\xi_i,\eta_i,\tau_i)\in\C\times\C\times\C^{d-1}\,:\,
     |\xi_i|<2,\,|\eta_i|<2,(\xi_i\eta_i,\tau_i)\in B_i\right\},
$$
where $B_i\subset\C\times\C^{d-1}$ is contained in
the open set $\{|\zeta_i|<4\}\times\C^{d-1}$.
By assumption (see~\ref{cWcUcV}), the fiber isomorphism
$f_0:P_0\to Q_0$ between the fibers over the origin
is the identity in these coordinates.

Consider the map
$$
\cW\to A\times (H^s)^{2k}\times B:
(a,\gamma,b)\mapsto(a,\alpha_1,\beta_1,\dots,\alpha_k,\beta_k,b)
$$
where $\gamma\in\cW(a,b)$ and
$\alpha_i=\xi_i\circ\gamma\circ x_i^{-1}$
and $\beta_i=\eta_i\circ\gamma\circ y_i^{-1}$.
This map is a diffeomorphism from $\cW$, with the
manifold structure of~\ref{cW}, onto an open subset
of the Hilbert manifold $A\times (H^s)^{2k}\times B$.
The image of the subset $\cU\subset\cW$ under this
diffomorphism consists of all tuples
$(a,\alpha_1,\beta_1,\dots,\alpha_k,\beta_k,b)$ in the
image of $\cW$ such that
$$
xy=z_i(a)\implies\alpha_i(x)\beta_i(y)=\zeta_i(b)\mbox{ for }i=1,\dots,k.
$$
That this subset is a complex submanifold of
$A\times (H^s)^{2k}\times B$ follows from
Theorem~\ref{thm:localmodel}

\bigskip

{\em We prove that $\cV$ is a smooth Hilbert submanifold of $\cW$.}

\medskipnoindent Define
$$
\cB := \left\{(a,v,b)\;\bigg|\;\begin{array}{l}
a\in A,\quad b\in B,\quad v\in H^{s+1/2}(M_a,V_b),\\
v(R_*\cap P_a)=S_*\cap Q_b,\;
v|\p N_a\in\cW(a,b)
\end{array}
\right\}
$$
and
\begin{equation}\label{eq:cZ}
\cZ:=\bigl\{(a,v,b)\in\cB\,:\,
v\in\Hol^s(M_a,V_b)
\bigr\}.
\end{equation}
We construct an auxiliary Hilbert manifold structure on $\cB$
and show that $\cZ$ is a smooth submanifold of $\cB$.
In analogy with~\ref{cWcUcV} denote
$$
\cB_0:=\left\{(a,v,b): \;\bigg|\;\begin{array}{l}
a\in A,\quad b\in B,\quad v\in H^{s+1/2}(\Omega,V_b),\\
v(r_*)=S_*\cap Q_b,\; v\circ\rho_a|\p N_a\in\cW(a,b)
\end{array}
\right\},
$$
where $\Omega:=M_{a_0}$ and $r_*:=R_*\cap\Omega=\rho(R_*)$
as in Definition~\ref{def:hardyTrivialization}.
This space is a Hilbert manifold and the Hardy trivialization
$(P=N\cup M,\iota,\rho)$ induces a bijection
$$
\cB_0\to\cB:(a,v,b)\mapsto(a,v\circ\rho_a,b).
$$
This defines the Hilbert manifold structure on $\cB$.
Note the commutative diagram
$$
\Rectangle{\cB_0}{}{\cB}{}{}{\cW_0}{}{\cW}.
$$
Here we identify $\p\Omega$ with $\Gamma$ via
the diffeomorphism $\iota_{a_0}$
(see Definition~\ref{def:hardyDecomposition}).
The bijection $\cB_0\to\cB$ identifies the subset
$\cZ\subset\cB$ with the subset
$\cZ_0\subset\cB_0$ given by
$$
\cZ_0:=\bigl\{(a,v,b)\in\cB_0,\;v\in\Hol^s((\Omega,j_a),Q_b)\bigr\},
$$
where $j_a:=(\rho_a)_*J|M_a$, $\rho_a:M_a\to\Omega$
is the Hardy trivialization,
and $J$ is the complex structure on $P$.
(Note that the map $a\mapsto j_a$ need not be holomorphic.)

\bigskip

{\em We prove that $\cZ_0$ is a smooth Hilbert submanifold of $\cB_0$.}

\medskipnoindent The tangent space of $\cB_0$ at a
triple $(a,v,b)$ is
\begin{equation*}
\begin{split}
T_{a,v,b}\cB_0
= T_aA\times \bigl\{&(\hat v,\hat b)
\in H^{s+1/2}(\Omega,v^*TQ)\times T_bB\,:\, \\
& d\pi_B(v)\hat v\equiv\hat b,\,\hat v(r_i)\in T_{v(s_i)}S_i\bigr\}.
\end{split}
\end{equation*}
Let $\cE\to\cB_0$ be the complex Hilbert space bundle
whose fiber
$$
\cE_{a,v,b} := H^{s-1/2}(\Omega,\Lambda_{j_a}^{0,1}T^*\Omega\otimes v^*TQ_b)
$$
over $(a,v,b)\in\cB_0$ is the Sobolev space of $(0,1)$-forms
on $(\Omega,j_a)$ of class $H^{s-1/2}$ with values in the
pullback tangent bundle $v^*TQ_b$. As before the
Cauchy--Riemann operator defines a smooth
section $\dbar:\cB_0\to\cE$ given by
\begin{equation}\label{eq:section}
\dbar(a,v,b) := \bar\p_{j_a,J}(v)
= \frac12\left(dv + J\circ dv\circ j_a\right).
\end{equation}
Here $J$ denotes the complex structure on $Q$.
The zero set of this section is the set $\cZ_0$ defined above.
It follows as in the proof of Theorem~\ref{thm:infuniv} that the
linearized operator $D_{a,v,b}:T_{a,v,b}\cB_0\to\cE_{a,v,b}$
is surjective and has a right inverse.  Hence the zero set $\cZ_0$
is a smooth Hilbert submanifold of $\cB_0$.

Again as in the proof of Theorem~\ref{thm:infuniv} restriction to
the boundary gives rise to a smooth embedding
$$
\cZ\to\cW:(a,v,b)\mapsto(a,\gamma,b),\qquad
\gamma:=v\circ\rho_a^{-1}|\p M_a,
$$
whose image is $\cV$.  The only difference in the proof that
the restriction map $\cZ\to\cW$ is proper is that now we
have a $\Cinf$ convergent sequence of complex structures
on $\Omega$.  The proof is otherwise word for word the same.
(Note that~\cite[Theorem~B.4.2]{MS2} allows for a sequence
of complex structures on the domain.) Hence $\cV$ is
a smooth Hilbert submanifold of $\cW$.

\bigskip

{\em We prove~(iii)}.

\medskipnoindent
That the projections $\cW\to A$ and $\cU\to A$ are holomorphic and
the projection $\cV\to A$ is smooth is obvious from the construction.
We prove that these three maps are submersions.  For the map
$\cU\to A$, and hence for $\cW\to A$, this follows immediately
from Proposition~\ref{prop:localmodel3}.   For $\cV$ this follows from the
fact that the linearized operator of the section~(\ref{eq:section}) is already
surjective when differentiating in the direction of a vector field $\hat v$
along $v$.  This proves~(iii).
\end{proof}


\section{Proofs of the main theorems}\label{sec:proof}

\dfn \label{def:regularNodal} The  set $C$ of critical points of
a nodal family $\pi:Q\to B$ is a submanifold of $Q$ and the restriction
of $\pi$ to this set is an immersion. The family is said to be
\jdef{regular nodal} at $b\in B$
if all self-intersections of $\pi(C)$ in $\pi^{-1}(b)$  are transverse, i.e.
$$
     \dim_\C\left(\im\,d\pi(q_1)\cap\cdots\cap \im\,d\pi(q_m)\right)
     = \dim_\C(B)-m
$$
whenever $q_1,\dots,q_m\in C$ are pairwise distinct
and $\pi(q_1)=\cdots=\pi(q_m)=b$;
the nodal family is called \jdef{regular nodal}
if it is regular nodal at each $b\in B$.
\nfd

\begin{lemma}\label{le:D}
Let $u$ be a desingularization of the fiber $Q_b$,
$g$ be the arithmetic genus of the fiber,  and $n$
be the number of marked points.
Then the following hold:
\begin{description}
\item[(i)]
We have $D_{u,b}(\hat u,\hat b)\in\cY_u$
for $(\hat u,\hat b)\in\cX_{u,b}$.
\item[(ii)]
The operator $D_{u,b}:\cX_{u,b}\to\cY_u$ is Fredholm.
\item[(iii)]
The Fredholm index satisfies
$$
\INDEX_\C(D_{u,b}) \ge 3-3g - n + \dim_\C(B)
$$
with equality if and only if $\pi$ is regular nodal  at $b$.
\end{description}
\end{lemma}

\bigbreak

\begin{proof}
We prove~(i). Choose $(\hat u,\hat b)\in\cX_{u,b}$ and
let
$$
\bar\p:\Omega^0(\Sigma,T_bB)\to \Omega^{0,1}(\Sigma,T_bB)
$$
denote usual Cauchy--Riemann operator. Then
$d\pi_B(u)D_u\hat u=\bar\p d\pi_B(u)\hat u=0$
since $d\pi_B(u)\hat u$ is a constant vector.
Hence $D_u\hat u\in\cY_u$.  Item~(ii) is immediate as  $D_u$
is Fredholm as a map from vertical vector fields to vertical
$(0,1)$-forms and $D_{u,b}$ is obtained from $D_u$ by a finite dimensional
modification of the domain. (A vertical vector field is an element
$\hat u\in \Omega^0(\Sigma,u^*TQ)$  such that $d\pi(u)\hat u=0$;
a vertical $(0,1)$-form is an element $\eta\in\Omega^{0,1}(\Sigma,u^*TQ)$
such that $d\pi(u)\eta=0$, i.e. an element of $\cY_u$.)

We prove~(iii).
The arithmetic  genus $g$ of the fiber is given by
\begin{equation}\label{eq:genus}
        g =\#\mbox{edges}-\#\mbox{vertices}+1+ \sum_i g_i
\end{equation}
where
$\#\mbox{vertices}=\sum_i1$ is the number of components of $\Sigma$,
$\#\mbox{edges}$ is the number of pairs of nodal points, and
$g_i$ is the genus of the $i$th component.
Now consider the subspace
$$
\cX_u:=\{(\hat u,\hat b)\in\cX_{u,b}:\hat b=0\}
$$
of all vertical vector fields along $u$
satisfying the nodal and marked point conditions.
If $(\hat u,\hat b)\in\cX_{u,b}$, then
the vector $\hat b$ belongs to the image of $d\pi(q)$ for every $q\in Q_b$.
Hence the codimension of $\cX_u$ in $\cX_{u,b}$ is
$$
        \codim_{\cX_{u,b}}(\cX_u)
       = \dim\left(\bigcap_{q\in Q_b}\im \,d\pi(q)\right)
       \ge \dim_\C(B)-\#\mbox{edges}
$$
with equality if and only if $\pi_B$ is regular nodal at $b$.
By Riemann-Roch  the restricted operator  has Fredholm
index
$$
         \INDEX_\C(D_u:\cX_u\to\cY_u)=\sum_i(3-3g_i)-2\,\#\mbox{edges}-n.
$$
Here the summand $-2\,\#\mbox{edges}$ arises from the nodal point
conditions in the definition of $\cX_u$.  To obtain the Fredholm index
of $D_{u,b}$ we must add the codimension
of $\cX_u$ in $\cX_{u,b}$ to the last identity. Hence
$$\renewcommand{\arraystretch}{1.5}
\begin{array}{lcl}
       \ds  \INDEX_\C(D_{u,b})
        &\ge &\ds \sum_i(3-3g_i)-3\,\#\mbox{edges}-n+\dim_\C(B)\\
       &= &\ds 3-3g-n+\dim_\C(B).
\end{array}
$$
The last identity follows from equation~(\ref{eq:genus}).
Again, equality holds if and only if $\pi_B$ is regular nodal at $b$.
\end{proof}

\begin{proof}[Proof of Theorem~\ref{thm:stable}.]
The proof is an easy application of the openness of trans\-versality.
Take $P=Q$, $A=B$, $\pi_A=\pi_B$, and $f_0=\id$,
so $\gamma_0$ is the inclusion of $\p N_{b_0}$ in $Q_{b_0}$.
Assume the unfolding $(\pi_B:Q\to B,S_*,b_0)$ is infinitesimally
universal. Choose $b\in B$ near $b_0$, fix a constant $s>7/2$,
and let $\cU_b$ and $\cV_b$ be the manifolds in~\ref{cWcUcV}
with $P=Q$ and $a=b$.
To show that $(\pi_B,S_*,b)$ is infinitesimally universal
we must show that $\cU_b$ and $\cV_b$ intersect transversally
at $\gamma$ where $\gamma$ is the inclusion of $\p N_b$ in $Q_b$.
Since $b$ is near $b_0$, so also is  $\gamma$ near $\gamma_0$.
Consider the three subspaces $T_\gamma \cW_b$,
$T_\gamma \cU_b$, $T_\gamma \cV_b$, of  $T_\gamma \cW$.
We have that $T_\gamma\cU_b=T_\gamma\cW_b\cap T_\gamma\cU$
and the intersection is transverse as the projection $\cU\to B$
is a submersion. Similarly,
$T_\gamma\cV_b=T_\gamma\cW_b\cap T_\gamma\cV$.
Hence the subspaces $T_\gamma\cU_b$ and $T_\gamma\cV_b$
depend contiuously on $(b,\gamma)$.
By part~(iv) of Theorem~\ref{thm:UV} the submanifolds
$\cU_{b_0}$ and $\cV_{b_0}$ intersect transversally at $\gamma_0$,
i.e. $T_{\gamma_0}\cW_{b_0}=T_{\gamma_0}\cU_{b_0}+T_{\gamma_0}\cV_{b_0}$.
Hence  $T_\gamma\cW_b=T_\gamma\cU_b+T_\gamma\cV_b$
for $(b,\gamma)$ near $(b_0,\gamma_0)$.
Hence the unfolding $(\pi_B,S_*,b)$ is infinitesimally universal
for $b$ near $b_0$ by Theorem~\ref{thm:UV} again as required.
\end{proof}

\begin{proof}[Proof of Theorem~\ref{thm:idm}.]
We proved `only if' in Section~\ref{sec:universal}.
For `if' assume that the unfolding $(\pi_B,S_*, b_0)$ is infinitesimally universal.
Let $(\pi_A,R_*,a_0)$ be another unfolding and $f_0:P_{a_0}\to Q_{b_0}$
be a fiber isomorphism. Assume the notation introduced
in Section~\ref{sec:hardy}.  In particular assume the
hypotheses of Theorem~\ref{thm:UV}.

\medskip\noindent{\bf Step~1.}
{\em We show that $\cU$ and $\cV$
intersect transversally at $(a_0,\gamma_0,b_0)$.}

\medskipnoindent
Abbreviate $w_0:=(a_0,\gamma_0,b_0)$. Choose
$\hat w\in T_{w_0}\cW$ and let $\hat a=d\pi(w_0)\hat w$.
As the restriction of $\pi$ to $\cU$ is a submersion
there is a vector $\hat u\in T_{w_0}\cU$ with
$d\pi(w_0)\hat u=\hat a$. Then $\hat w-\hat u\in \cW_{a_0}$
so by part~(iv) of Theorem~\ref{thm:UV} there are vectors
$\hat u_0\in T_{w_0}\cU_{a_0}$ and $\hat v_0\in T_{w_0}\cV_{a_0}$
with $\hat w -\hat u= \hat u_0+\hat v_0$.

\medskip\noindent{\bf Step~2.}
{\em We show that the projection $\cU\cap\cV\to A$
is a diffeomorphism.}

\medskip\noindent
By Step~1 the intersection $\cU\cap\cV$ is a smooth submanifold
of $\cW$ (after shrinking) and
$$
     T_w(\cU\cap \cV)=(T_w\cU)\cap(T_w\cV)
$$
for $w\in\cU\cap \cV$. By the inverse function theorem it is enough
to show that $d\pi(w_0):T_w(\cU\cap \cV)\to T_{a_0}A$ is bijective.
Injectivity follows from part~(iv) of Theorem~\ref{thm:UV} and the fact that
$T_{w_0}\cU_{a_0}=T_{w_0}\cU\cap\ker d\pi(w_0)$ and
$T_{w_0}\cV_{a_0}=T_{w_0}\cV\cap\ker d\pi(w_0)$.
We prove surjectivity. Choose $\hat a\in T_{a_0}A$. Since
the restrictions of $\pi$ to $\cU$ and $\cV$ are submersions,
there exist tangent vectors $\hat u\in T_{w_0}\cU$ and
$\hat v\in T_{w_0}\cV$ with $d\pi(w_0)\hat u= d\pi(w_0)\hat v=\hat a$.
The difference $\hat u-\hat v$ lies in $T_{w_0}\cW_{a_0}$
so, by part~(iv) of Theorem~\ref{thm:UV},  there are vectors
$\hat u_0\in T_{w_0}\cU_{a_0}$ and $\hat v_0\in T_{w_0}\cV_{a_0}$
with $\hat u-\hat v= \hat u_0+\hat v_0$.
Hence $\hat u -\hat u_0=\hat v+\hat v_0$
lies in $T_{w_0}(\cU\cap\cV)$ and projects to~$\hat a$.

Now define $\Phi:P\to Q$ and $\phi:A\to B$ by
$$
\phi(a):=b_a,\qquad \Phi|P_a:=f_a,
$$
where $f_a:P_a\to Q_{b_a}$ is the unique fiber isomorphism that satisfies
$f_a|\p N_a=\gamma_a$. (See Lemma~\ref{le:fgamma}.)

\medskip\noindent{\bf Step~3.}
{\em The maps $\phi$ and  $\Phi$
and smooth.}

\medskipnoindent
The map $\phi$ is the composition of the inverse of the projection
$\cU\cap\cV\to A$ with the projection $\cU\cap\cV\to B$
and is therefore smooth by Step~2.

By  Lemma~\ref{le:smooth} the restriction of $\Phi$
to $\INT(N)$ is holomorphic  and hence smooth.
To prove that the restriction of $\Phi$ to $\INT(M)$ is smooth
write it as the composition
$$
\INT(M)\to
\left\{(a,v,b,p)\,:\,(a,v,b)\in\cZ,\,p\in\INT(M_a)\right\}
\to Q
$$
where $\cZ$ is given by~(\ref{eq:cZ}) in
the proof of Theorem~\ref{thm:UV}; the
first map sends $p\in\INT(M)$ to the quadruple
$(a,v_a,\phi(a),p)$, where $(a,v_a,\phi(a))\in\cZ$ is the unique
point with $(a,v_A|\p N_a,\phi(a))\in\cU_a\cap\cV_a$,
and the second map is the evaluation map
$(a,v,b,p)\mapsto v(p)$.  Use the Hardy trivialization
$\rho$  as in the proof of
Theorem~\ref{thm:UV} to define an auxiliary Hilbert manifold structure
on the set $\left\{(a,v,b,p)\,:\,(a,v,b)\in\cZ,\,p\in\INT(M_a)\right\}$.
Then both maps are smooth and hence so is their composition
$\Phi$ on $\INT(M)$.

Thus we have proved that $\Phi$ is smooth on $P\setminus\p N$.
Repeat the same argument with a different Hardy trivialization
$(P=N'\cup M',\iota',\rho')$ such that $N'\subset\INT(N)$.
Then the resulting morphisms $\Phi$ and $\Phi'$ must agree
because each $\Phi_a'$ for $a$ close to $a_0$ is a fiber isomorphism
near $f_0$ and so determines an intersection point
$(\gamma_a',b_a')\in \cU_a\cap\cV_a$.  By uniqueness, this
point agrees with $(\Phi_a|\p N_a,\phi(a))$ and hence $\Phi_a=\Phi'_a$
as claimed.  It follows that $\Phi$ is smooth.

\medskip\noindent{\bf Step~4.}
{\em The maps $\phi$ and  $\Phi$
and holomorphic.}

\medskipnoindent By Theorems~\ref{thm:stable} and~\ref{thm:teich-mark}
$\phi$ is holomorphic on the open set $A\setminus\pi_A(C_A)$ and
$\Phi$ is holomorphic on the open set $P\setminus\pi_A^{-1}(\pi_A(C_A))$.
Hence by Step~3 they are holomorphic everywhere.
This proves the theorem.
\end{proof}

\begin{proof}[Proof of Theorem~\ref{thm:pseudo}.]
Assume that the unfolding $(\pi_B,S_*, b_0)$ is infinitesimally universal
and let $(\phi,\Phi)$ is a pseudomorphism
from $(\pi_A,R_*, a_0)$ to $(\pi_B,S_*, b_0)$.
Then, in the notation of the proof of Theorem~\ref{thm:idm},
we have that $(\gamma_a,b_a):=(\Phi|\p N_a,\phi(a))$ is the unique
intersection point of $\cU_a$ and $\cV_a$.  Hence $(\phi,\Phi)$
agrees with the unique (holomorphic) morphism constructed in the
proof of Theorem~\ref{thm:idm}.
\end{proof}

\begin{proof}[Proof of Theorem~\ref{thm:existence}]
We proved `only if' in Section~\ref{sec:universal}.
For `if' assume that $(\Sigma,s_*,\nu,j)$ is stable.
We first consider the case where
$\Sigma$ is disconnected and there are no
nodal points.  Let $\Sigma_1,\dots,\Sigma_k$ be the components
of $\Sigma$, $g_j$ be the genus of $\Sigma_j$, and $n_j$ be the number
of marked points on $\Sigma_j$. Let $I_j\subset\{1,\dots,n\}$
be the index set associated to the marked points in $\Sigma_j$.
Then $\{1,\dots,n\}$ is the disjoint union of the sets
$I_1,\dots,I_k$ and
$
n_j = |I_j| > 2-2g_j.
$
By Theorem~\ref{thm:teich-mark}
there exists, for each $j$, a universal unfolding
$(\pi_j:Q_j\to B_j,\{S_{ji}\}_{i\in I_j},b_{0j},v_{0j})$
of $\Sigma_j$.
Note that
$
\dim_\C(B_j)=3g_j-3+n_j.
$
Define
$$
B_0 := B_1\times\cdots\times B_k,
$$
$$
Q_0 := \bigsqcup_{j=1}^k B_1\times\cdots\times B_{j-1}\times Q_j
\times B_{j+1}\times\cdots\times B_k,
$$
$$
\pi_0(b_1,\dots,b_{j-1},q_j,b_{j+1},\dots,b_k)
:= (b_1,\dots,b_{j-1},\pi_j(q_j),b_{j+1},\dots,b_k),
$$
$$
S_{0i}:=\left\{(b_1,\dots,b_{j-1},q_j,b_{j+1},\dots,b_k)\,:\,
q_j\in S_{ij}\right\},\qquad i\in I_j,
$$
$$
b_0:=(b_{01},\dots,b_{0k}),
$$
$$
v_0(z) := (b_{01},\dots,b_{0,j-1},
v_{0j}(z),b_{0,j+1},\dots,b_{0k}),\qquad
z\in\Sigma_j.
$$
Then the quadruple $(\pi_0,S_{0*},b_0,v_0)$ is a universal
unfolding of $\Sigma$.

Next consider the general case.  Denote the nodal points on
$\Sigma$ by $\nu=\{\{r_1,s_1\},\dots,\{r_m,s_m\}\}$
and the marked points by $t_1,\dots,t_n$.
Assume, without loss of generality, that the signature
of $(\Sigma,t_*,\nu,j)$ is a connected graph
(see Definition~\ref{def:signature}).
Replace all the nodal points by marked points.
Then, by what we have just proved, there exists a universal
unfolding $(\pi_0,R_{0*},S_{0*},T_{0*},b_0,v_0)$
of $(\Sigma,r_*,s_*,t_*,j)$. Choose disjoint open sets
$U_1,\dots,U_m,V_1,\dots,V_m\subset Q_0$
such that
$$
R_{0i}\subset U_i,\qquad
S_{0i}\subset V_i,\qquad
U_i\cap T_{0j}=V_i\cap T_{0j}=\emptyset
$$
for $i=1,\dots,m$ and $j=1,\ldots,n$.  Choose holomorphic functions
$x_i:U_i\to\C$ and $y_i:V_i\to\C$ such that
$$
x_i(R_{0i}) = 0,\qquad
y_i(S_{0i}) = 0
$$
and $(\pi_0,x_i)$ and $(\pi_0,y_i)$ are coordinates
on $Q_0$.  Shrink $B$ and the open sets $U_i,V_i$, if necessary.
Assume without loss of generality that $x_i(U_i)=y_i(V_i)=\D$.
Define
$$
B := B_0\times\D^m,\qquad
Q := Q_0\times\D^m/\sim.
$$
Two points $(q,z)$ and $(q',z)$ with
$q\in U_i$ and $q'\in V_i$ are identified if and only
if $\pi_0(q)=\pi_0(q')$ and either
$$
x_i(q)y_i(q') = z_i\ne 0\qquad\mbox{or}\qquad
x_i(q) = y_i(q') = z_i = 0.
$$
The equivalence relation on $Q_0\times\D^m$ is generated by these
identifications.  (Two points $(q,z)$ and $(q',z)$ with
$\pi_0(q)=\pi_0(q')$, $q\in U_i$, $q'\in V_i$, $z_i=0$
are {\em not} identified in the case
$x_i(q)=0$ and $y_i(q')\ne 0$ nor in the case
$x_i(q)\ne 0$ and $y_i(q')=0$.)
The projection $\pi:Q\to B$ and the sections
$T_j:B\to Q$ are defined by
$$
\pi([q,z]) := (\pi_0(q),z),\qquad
T_j:=\left\{[q,z]\,:\,q\in T_{0j}\right\}
$$
for $j=1,\dots,n$.

We have thus defined $Q$ as a set.  The manifold structure
is defined as follows.  For $i\in\{1,\dots,m\}$
denote by $C_i\subset Q$ the set of all equivalence classes
$[q,z]\in Q$ that satisfy $z_i=0$ and $q\in R_{0i}$.
Note that any such point is equivalent to the pair $[q',z]$ with
$q'\in S_{0i}$ and $\pi_0(q')=\pi_0(q)$.  Let
$$
     C := \bigcup_{i=1}^mC_i.
$$
The manifold structure on $Q\setminus C$ is induced by the
product manifold structure on $Q_0\times\D^m$.
We now explain the manifold structure near $C_i$.
Fix a constant $0<\eps<1$ and define an open neighborhood
$N_i\subset Q$ of $C_i$ by
\begin{eqnarray*}
     N_i
&:= &
     C_i \cup\left\{[q,z]\in Q:
     q\in U_i,\;\frac{|z_i|}{\eps}<|x_i(q)|<\eps\right\} \\
&&
     \cup \left\{[q,z]\in Q:
     q\in V_i,\;\frac{|z_i|}{\eps}<|y_i(q)|<\eps\right\}.
\end{eqnarray*}
A coordinate chart on $N_i$ is the map
$$
     [q,z]\mapsto
     (b_0,z_1,\dots,z_{i-1},x_i,y_i,z_{i+1},\dots,z_m),
$$
where $b_0:=\pi_0(q)\in B_0$,
$$
     x_i := \left\{\begin{array}{rl}
     x_i(q),&\mbox{if }q\in U_i,\\
     z_i/y_i(q),&\mbox{if }
     q\in V_i,\,z_i\ne 0, \\
     0,&\mbox{if }
     q\in V_i,\,z_i=0,
     \end{array}\right.
$$
(if $[q,z]\in N_i$ and $q\in V_i$ then $z_i\ne 0$ implies
$y_i(q)\ne 0$), and
$$
     y_i := \left\{\begin{array}{rl}
     y_i(q),&\mbox{if }q\in V_i,\\
     z_i/x_i(q),&\mbox{if }
     q\in U_i,\,z_i\ne 0, \\
     0,&\mbox{if }
     q\in U_i,\,z_i=0.
     \end{array}\right.
$$
With this construction the transition maps are
holomorphic and so $Q$ is a complex manifold.
In the coordinate chart on $N_i$ the projection $\pi$
has the form
$$
     (b_0,z_1,\dots,z_{i-1},x_i,y_i,z_{i+1},\dots,z_m)
     \mapsto (b_0,z),\qquad
     z_i := x_iy_i.
$$
It follows that $\pi$ is holomorphic, the critical set
of $\pi$ is $C$, and each critical point is nodal.
Moreover, $\pi$ restricts to a diffeomorphism from $C_i$
onto the submanifold $\{z_i=0\}\subset B$.
Hence $\pi$ is a regular nodal  family
(see Definition~\ref{def:regularNodal}).

Denote $b:=(b_0,0)\in B$, let $\iota:Q_0\to Q$ be the
holomorphic map defined by
$
\iota(q):=[q,0],
$
and define $v:\Sigma\to Q$ by
$
v:=\iota\circ v_0.
$
Then $v$ is a desingularization of the fiber $Q_b=\pi^{-1}(b)$.

We prove that the triple $(\pi,T_*,b)$ is a universal
unfolding. Since the signature of the marked nodal Riemann surface
$\Sigma$ is a connected graph, the first Betti number of this
graph is $1-k+m$ (since $m$ is the number of edges,
i.e. of equivalent pairs of nodal points,
and $k$ is the number of vertices, i.e. of components
of $\Sigma$).  Hence the arithmetic genus $g$
(see Definition~\ref{def:arithGenus})
of the central fiber $Q_b$ is given by
$$
g - 1 = m + \sum_{j=1}^k(g_j-1).
$$
Now recall that $n_j$ is the number of special points
on $\Sigma_j$ and
$$
\sum_{j=1}^kn_j=n+2m.
$$
Since $\dim_\C(B_j)=3g_j-3+n_j$ this implies
$$
\dim_\C(B)
= \dim_\C(B_0) + m
= \sum_{j=1}^k(3g_j-3+n_j) + m
= 3g-3+n.
$$
Since the Riemann family $\pi:Q\to B$ is  regular nodal
it follows from lemma~\ref{le:D} that
the operator $D_{v,b}$ (see Definition~\ref{def:D})
has Fredholm index zero.  Hence $D_{v,b}$ is bijective
if and only if it is injective.

We prove in three steps that $D_{v,b}$
is injective. First, every vector
$(\hat v,\hat b)\in\cX_{v,b}$ with
$
\hat b=:(\hat b_0,\hat z)
$
satisfies $\hat z=0$.  To see this note that
$d\pi(v)\hat v\equiv\hat b$.   For every $i$ there is
a unique pair of equivalent nodal points in $\Sigma$
that are mapped to $C_i$ under $v$. Since the image
of $d\pi$ at each point in $C_i$ is contained in
the subspace $\{\hat z_i=0\}$ it follows that $\hat z=0$.
Second, we define a linear operator
$$
\cX_{v_0,b_0}\to X_{v,b}:(\hat v_0,\hat b_0)\mapsto(\hat v,\hat b)
$$
by
$$
\hat b:=(\hat b_0,0),\qquad
\hat v(s) := (\hat v_0(s),0)\in T_{v(s)}Q,
$$
for $s\in\Sigma\setminus\{r_1,\dots,r_m,s_1,\dots,s_m\}$.
Then $\hat v$ extends uniquely to a smooth vector field along $v$.
In the above coordinates on $N_i$ the tangent vector
$\hat v(r_i)=\hat v(s_i)\in T_{v(r_i)}Q=T_{v(s_i)}Q$
has the form $(\hat x_i,\hat y_i,\hat b_0,0)$, where
$
\hat x_i := dx_i(v_0(r_i))\hat v_0(r_i)
$
and
$
\hat y_i := dy_i(v_0(s_i))\hat v_0(s_i).
$
It is easy to see that this operator is bijective.
Third, since the map $\iota:Q_0\to Q$ is holomorphic, it follows
that
$$
D_{v,b}(\hat v,\hat b) = d\iota(v)D_{v_0,b_0}(\hat v_0,\hat b_0).
$$
Hence the operator $\cX_{v_0,b_0}\to X_{v,b}$ restricts to a
vector space isomorphism from the kernel of $D_{v_0,b_0}$
to the kernel of $D_{v,b}$.   By construction, the operator
$D_{v_0,b_0}$ is injective and hence, so is $D_{v,b}$.
Thus we have proved that $D_{v,b}$ is bijective.
Hence, by Theorem~\ref{thm:idm}, the quadruple
$(\pi,T_*,b,v)$ is a universal unfolding of
$(\Sigma,s_*,\nu,j)$.
\end{proof}


\section{Topology}\label{sec:topology}

The orbit space of a groupoid  inherits a topology from
an orbifold structure on the groupoid. This topology
is independent of the choice of the structure in the sense
that equivalent orbifold structures determine the same topology
(see Section~\ref{sec:orbifold}). In the case of the
Deligne--Mumford orbifold $\bar\cM_{g,n}$, the topology
has as a basis for the open sets the collection of all sets
$\{[\Sigma_b]_\cB:b\in U\}$ where $(\pi_B:Q\to B,S_*)$ is a
universal family as in Definition~\ref{def:universal},
the functor
$$
B\to\bar\cB_{g,n}:b\mapsto\Sigma_b
$$
is the corresponding orbifold structure as in Definition~\ref{B-Gamma}, 
and $U$ runs over all open subsets of $B$. 
In section~\ref{sec:compact} we show that $\bar\cM_{g,n}$ is
compact and Hausdorff. (See Example~\ref{ex:nonHausdorff}
for an example which shows why it is not obvious that the
moduli space is Hausdorff.) For this purpose we introduce
in this section
a notion of convergence of sequences of marked nodal Riemann
surfaces which we call DM-convergence.

\para\label{suture}
Let $\Sigma$ be a compact oriented surface
and $\gamma\subset\Sigma$ be a disjoint union of
embedded circles. We denote by $\Sigma_\gamma$
the compact surface with boundary which results
by \jdef{cutting open} $\Sigma$ along $\gamma$.
This implies that there  is  a local embedding
$$
\sigma:\Sigma_\gamma\to\Sigma
$$
which
maps $\INT(\Sigma_\gamma)$ one to one onto $\Sigma\setminus\gamma$
and
maps $\p \Sigma_\gamma$ two to one onto $\gamma$.
One might call $\sigma$ the {\em suture map} and $\gamma$ the {\em incision}.
\arap

\dfn \label{deformation}
Let   $(\Sigma',\nu')$ and $(\Sigma,\nu)$ be nodal surfaces.
A smooth map $\phi:\Sigma'\setminus\gamma'\to\Sigma$
is called a $(\nu',\nu)$-\jdef{deformation}
iff  $\gamma'\subset\Sigma'\setminus\bigcup\nu'$ is a disjoint union
of embedded circles such that (where $\sigma:\Sigma'_{\gamma'}\to\Sigma'$
is the suture map just defined) we have
\begin{itemize}
\item
$
\phi_*\nu':=\bigl\{  \{\phi(y'_1),\phi(y'_2)\}:
\{y'_1,y'_2\}\in\nu'\bigr\} \subset\nu.
$
\item
$\phi$ is a diffeomorphism from $\Sigma'\setminus \gamma'$
onto $\Sigma\setminus \gamma$, where
$\gamma:=\bigcup(\nu\setminus\phi_*\nu')$.
\item
$\phi\circ\sigma|\INT(\Sigma'_{\gamma'})$
extends to a continuous surjective map
$\Sigma'_{\gamma'}\to\Sigma$ such that
the preimage of each nodal point in $\gamma$
is a component of $\p\Sigma'_{\gamma'}$ and two boundary components
which map under $\sigma$ to the same component of $\gamma'$
map to a nodal pair $\{x,y\}\in\gamma$.
\end{itemize}
A sequence $\phi_k:(\Sigma_k\setminus\gamma_k,\nu_k)\to(\Sigma,\nu)$
of $(\nu_k,\nu)$-deformations
is called \jdef{monotypic} iff $(\phi_k)_*\nu_k$
is independent of $k$.
\nfd

\dfn\label{def:convergence}
A sequence $(\Sigma_k,s_{k,*},\nu_k,j_k)$ of marked
nodal Riemann surfaces of type $(g,n)$ is said to
\jdef{converge monotypically} to a  marked nodal
Riemann surface $(\Sigma,s_*,\nu,j)$ of type $(g,n)$ iff
there is a monotypic sequence
$\phi_k:\Sigma_k\setminus\gamma_k \to\Sigma$
of $(\nu_k,\nu)$-deformations such that
for $i=1,\ldots,n$ the sequence $\phi_k(s_{k,i})$
converges to $s_i$ in $\Sigma$, and
the sequence $(\phi_k)_*j_k$ of complex structures on $\Sigma\setminus \gamma$
converges to $j|(\Sigma\setminus \gamma)$ in the $\Cinf$ topology.
The sequence $(\Sigma_k,s_{k,*},\nu_k,j_k)$ is said to
\jdef{DM-converge} to $(\Sigma,j,s,\nu)$ iff,
after discarding finitely many terms, it is the disjoint union
of finitely many sequences which converge monotypically to  $(\Sigma,s,\nu,j)$.
\nfd

\rmk
Assume
that $(\Sigma_k,s_{k,*},\nu_k,j_k)$ DM-converges to $(\Sigma,s_*,\nu,j)$,
that $(\Sigma_k,s_{k,*},\nu_k,j_k)$ is isomorphic
to  $(\Sigma'_k,s'_{k,*},\nu'_k,j'_k)$, and
that $(\Sigma,s_*,\nu,j)$ is isomorphic to $(\Sigma',s'_*,\nu',j')$.
Then $(\Sigma'_k,s'_{k,*},\nu'_k,j'_k)$
DM-converges to $(\Sigma',s'_*,\nu',j')$.
\kmr

\rmk
Our definition of deformation agrees with~\cite[page 79]{HUMMEL}.
Our definition of monotypic convergence is Hummel's
definition of weak convergence to cusp curves
in~\cite[page~80]{HUMMEL}
(with the target manifold $M$ a point)
except that he does not allow marked points.
However, the conclusion of Proposition~5.1
in~\cite[page 71]{HUMMEL}
allows marked points in the guise of what
Hummel calls degenerate boundary components.
We will apply Proposition~5.1 of~\cite{HUMMEL} in the
proof of Theorem~\ref{thm:compact} below
after some preliminary constructions to fit its hypotheses.
\kmr

\begin{theorem}\label{thm:convergence}
Let $(\pi:Q\to B,S_*,b_0)$ be a universal unfolding
of a marked nodal Riemann surface $(\Sigma_0,s_{0,*},\nu_0,j_0)$
of type $(g,n)$ and $(\Sigma_k,s_{k,*},\nu_k,j_k)$ be a
sequence of marked nodal Riemann surfaces of type $(g,n)$.
Then the following are equivalent.
\begin{description}
\item[(i)]
The sequence $(\Sigma_k,s_{k,*},\nu_k,j_k)$ DM-converges
to $(\Sigma_0,s_{0,*},\nu_0,j_0)$.
\item[(ii)]
After discarding finitely many terms
there is a sequence $b_k\in B$ such that $b_k$
converges to $b_0$ and $(\Sigma_k,s_{k,*},\nu_k,j_k)$
arises from a desingularization $u_k:\Sigma_k\to Q_{b_k}$.
\end{description}
\end{theorem}

We postpone the proof of Theorem~\ref{thm:convergence} till
after we treat the analogous theorem for fiber isomorphisms.

\dfn\label{def:proper}
Let $(\pi_A:P\to A,R_*,a_0)$ and $(\pi_B:Q\to B,S_*b_0)$
be two universal unfoldings of type $(g,n)$ and
$a_k\to a_0$ and $b_k\to b_0$ be convergent sequences.
A sequence of fiber isomorphisms $f_k:P_{a_k}\to Q_{b_k}$ is said to
\jdef{DM-converge} to a fiber isomorphism $f_0:P_{a_0}\to Q_{b_0}$
iff for every Hardy decomposition $P=M\cup N$
as in Definition~\ref{def:hardyDecomposition}
the sequence $f_k\circ\iota_{a_k}$
converges to $f_0\circ\iota_{a_0}$ in  the $\Cinf$ topology.
\nfd

\begin{theorem}\label{thm:convergence-f}
Let $(\Phi,\phi):(P,A)\to(Q,B)$ be the germ of a morphism satisfying
$\phi(a_0)=b_0$ and $\Phi_{a_0}=f_0$.  Then the following are equivalent.
\begin{description}
\item[(i)]
The sequence $(a_k,f_k,b_k)$ DM-converges to $(a_0,f_0,b_0)$.
\item[(ii)]
For $k$ sufficiently large we have $\phi(a_k)=b_k$ and $\Phi_{a_k}=f_k$.
\end{description}
\end{theorem}

\begin{proof}
That~(ii) implies~(i) is obvious. We prove that~(i) implies~(ii).
Recall the Hardy decomposition  as in the definition of
the spaces $\cU$, $\cV$, $\cW$ in~\ref{cWcUcV}.
The proof of Theorem~\ref{thm:idm} in Section~\ref{sec:proof}
shows that
$$
    \bigl(a,\Phi_a|\p N\cap P_a,\phi(a)\bigr) = \cU_a\cap \cV_a.
$$
But $(a_k,f_k|\p N\cap P_{a_k},b_k)
\in\cU_{a_k}\cap \cV_{a_k}\subset\cW$
for $k$ sufficiently large by DM-convergence.
Both sequences $(a_k,\Phi_{a_k}|\p N\cap P_{a_k},\phi(a_k))$
and $(a_k,f_k|\p N\cap P_{a_k},b_k)$ converge
to the same point $(a_0,f_0|\p N\cap P_{a_0},b_0)$.
Hence by transversality in Therorem~\ref{thm:UV}
they are equal for large $k$.  Now use Lemma~\ref{le:fgamma}.
\end{proof}

\begin{proof}[Proof of Theorem~\ref{thm:convergence}.]
We prove that~(ii) implies~(i).
Let $u_0:\Sigma_0\to Q_{b_0}$ be a desingularization.
Assume that $b_k$ converges to $b$ and that
$u_k:\Sigma_k\to Q_{b_k}$ is a sequence of desingularizations.
Choose a Hardy trivialization
$(Q=M\cup N,\iota,\rho)$ for $(\pi,S_*,b_0)$ as in
Definition~\ref{def:hardyTrivialization}.
For each $b\in B$ choose a smooth map
$$
\psi_b:Q_b\to Q_{b_0}
$$
as follows.
The restriction of $\psi_b$ to $M_b$ agrees with $\rho_b$.
Next,  using the nodal coordinates of Definition~\ref{def:hardyDecomposition},
extend $\psi_b$ to a neighborhood of the common boundary of $M$ and $N$ via
$(x_i,0,t_i)\mapsto(x_i,0,0)$ for $2\sqrt{|z_i(b)|}\le|x_i|\le1$
and
$(0,y_i,t_i)\mapsto(0,y_i,0)$ for $2\sqrt{|z_i(b)|}\le|y_i|\le1$.
Finally,  when $z_i(b)\ne 0$, extend  to a smooth map
$Q_b\cap N_i\to Q_{b_0}\cap N_i$ that maps the circle
$|x_i|=|y_i|=\sqrt{|z_i(b)|}$ onto the nodal point $q_i$ and is a
diffeomorphism from the complement of this circle in
$Q_b\cap N_i$ onto the complement of $q_i$
in $Q_{b_0}\cap N_i$.
Denote
$$
\gamma_k:=\bigcup_{z_i(b_k)\ne 0}
u_k^{-1}\bigl(\bigl\{q\in Q_b\cap N_i\,:\,
|x_i(q)|=|y_i(q)|=\sqrt{|z_i(b_k)|}\bigr\}\bigr)
\subset\Sigma_k.
$$
Then, for every $k$, there is a unique smooth map
$\phi_k:\Sigma_k\setminus \gamma_k\to\Sigma_0$
such that
$$
u_0\circ\phi_k = \psi_{b_k}\circ u_k:
\Sigma_k\setminus \gamma_k\to Q_{b_0}.
$$
It follows that $\phi_k$ is a sequence of deformations as in
Definition~\ref{deformation} and that this sequence satisfies the
requirements of Definition~\ref{def:convergence}.
(The sequence $\phi_k$ is monotypic whenever there is
an index set $I$ such that, for every $k$, we have
$z_i(b_k)=0$ for $i\in I$ and $z_k(b_i)\ne 0$ for $i\notin I$;
after discarding finitely many terms, we can write $\phi_k$
as a finite union of monotypic sequences.)
Hence the sequence $(\Sigma_k,s_{k,*},\nu_k,j_k)$
DM-converges to $(\Sigma,j,s_*,\nu)$.
Thus we have proved that~(ii) implies~(i).

We prove that~(i) implies~(ii).
Let $(\Sigma_k,s_{k,*},\nu_k,j_k)$
be a sequence of marked nodal Riemann surfaces
of type $(g,n)$ that DM-converges to $(\Sigma,j,s_*,\nu)$.
If $\Sigma$ has no nodes then
$\Sigma_k$ has no nodes and the maps $\phi_k:\Sigma_k\to\Sigma$
in the definition of DM-convergence are diffeomorphisms.
Since $(\phi_k)_*j_k$ converges to $j$, assertion~(ii)
follows from the fact that a slice in $\cJ(\Sigma)$ determines
a universal unfolding. The same reasoning works
when $(\Sigma_k,s_{k,*},\nu_k)$ has the
same signature as $(\Sigma,s_*,\nu)$. To avoid excessive
notation we consider  the case where $(\Sigma,\nu)$
has precisely one node and $(\Sigma_k,\nu_k)$ has no nodes, i.e.
$$
\nu_k=\emptyset,\qquad \nu=:\{\{z_0,z_\infty\}\}.
$$
Choose holomorphic diffeomorphisms
$$
x: (\Delta_0,z_0)\to(\D,0), \qquad y: (\Delta_\infty,z_\infty)\to(\D,0),
$$
where  $\Delta_0,\Delta_\infty\subset\Sigma$ are
disjoint closed disks
centered at $z_0,z_\infty$ respectively  and
$$
\Delta:=\Delta_0\cup \Delta_\infty
$$
does not contain  any marked points. For $\delta\in(0,1)$ let
$\Delta(\delta):=\Delta_0(\delta)\cup\Delta_\infty(\delta)$ where
$$
\Delta_0(\delta):= \{p\in \Delta_0:|x(p)|\le\delta\},\qquad
\Delta_\infty(\delta):=\{q\in \Delta_\infty:|y(q)|\le\delta\}.
$$
A  decreasing sequence  $\delta_k\in(0,1)$ converging  to zero determines
a sequence of decompositions
$$
  \Sigma = \Omega_k\cup\Delta(\delta_k), \qquad
\p\Omega_k=\p\Delta(\delta_k)=\Omega_k\cap\Delta(\delta_k).
$$
Thus $\Omega_k$ is obtained from $\Sigma$ by removing a
nested sequence of pairs of open disks centered at the nodal
points so
$
      \bigcup_k \Omega_k=\Sigma\setminus\{z_0,z_\infty\},
$
$\Omega_k\subset \Omega_{k+1}$, and
$\Omega_k\cap\Delta=(\Omega_k\cap\Delta_0)\cup(\Omega_k\cap\Delta_\infty)$
is a disjoint union of two closed annuli.

\medskipnoindent{\bf Claim.}
{\em There are sequences of real numbers
$\delta_k$, $r_k$, $\theta_{0 k}$, $\theta_{\infty k}$,
smooth embeddings
$$
f_k:\Omega_k\to\Sigma_k,\qquad
\xi_k:\A(\delta_k,1)\to\A(r_k,1),\qquad \eta_k:\A(\delta_k,1)\to\A(r_k,1),
$$
and holomorphic diffeomorphisms
$$
h_k:\A(r_k,1)\to A_k:=\Sigma_k\setminus f_k(\Sigma\setminus\Delta),
$$
satisfying the following conditions.
\begin{description}
\item[1)] $f_k^*j_k$ converges to $j$ in the $\Cinf$ topology
on $\Sigma\setminus\{z_0,z_\infty\}$.
\item[2)] $f_k^*j_k$ is equal to $j$ on $\Omega_k\cap \Delta(\half)$.
\item[3)] $\xi_k(S^1)=\eta_k(S^1)=S^1$.
\item[4)] $\ds h_k(\xi_k(x(p))) = f_k(p)$
for $p\in \Omega_k\cap\Delta_0$.
\item[5)] $\ds h_k\left(\frac{r_k}{\eta_k(y(q))}\right) = f_k(q)$
for $q\in\Omega_k\cap\Delta_\infty$.
\item[6)] $\xi_k(x(p))=e^{i\theta_{0k}}x(p)$
for $p\in\Omega_k\cap \Delta_0(\half)$.
\item[7)] $\eta_k(y(q))=e^{i\theta_{\infty k}}y(q)$
for $q\in\Omega_k\cap \Delta_\infty(\half)$.
\end{description}
}

\medskipnoindent{\bf Proof of the Claim.}
Let $\delta_k\in(0,1]$ be any sequence decreasing to zero,
for example $\delta_k:=1/k$, and denote
$\Omega_k:=\Sigma\setminus\INT(\Delta(\delta_k))$.
Define $f_k:\Omega_k\to\Sigma_k$ by
$
f_k:=\left(\phi_k|\phi_k^{-1}(\Omega_k)\right)^{-1}
$
where
$\phi_k:\Sigma_k\setminus \gamma_k\to\Sigma$
is as in Definition~\ref{def:convergence}.
Then $f_k$ satisfies~1). We will modify $\delta_k$ and
$f_k$ to satisfy the other conditions.

By the path lifting arguments used in
in Section~\ref{sec:S2} (see also Appendix~C.5 of~\cite{MS2})
there is a sequence of holomorphic embeddings
$$
g_k:(\Omega_k\cap \Delta,j)\to(\Omega_k,f_k^*j_k)
$$
that converges to the identity in the $\Cinf$ topology
and preserves the boundary of $\Omega_k$.
Extend $g_k$ to a diffeomorphism, still denoted by $g_k:\Omega_k\to \Omega_k$,
so that the extensions converge to the identity in the $\Cinf$ topology
and  replace $f_k$ by $f_k\circ g_k$.  This new sequence  satisfies~1) and~2);
in fact,  $f_k$ is now holomorphic on $\Omega_k\cap\Delta$.
(Below we modify $f_k$ again.)

The set $A_k\subset\Sigma_k$ is an annulus with boundary
$f_k(\p \Delta_0)\cup f_k(\p \Delta_\infty)$  so
there is a unique number $r_k>0$ and a holomorphic diffeomorphism
$
h_k:\A(r_k,1)\to A_k,
$
unique up to composition with a rotation, such that
$$
h_k(S^1)=f_k(\p \Delta_0),\qquad
h_k(r_kS^1)=f_k(\p \Delta_\infty).
$$
The embeddings $\xi_k:\A(\delta_k,1)\to \A(r_k,1)$ and
$\eta_k:\A(\delta_k,1)\to \A(r_k,1)$  defined by~4)  and~5)
satisfy~3);  they are holomorphic because $f_k$ is holomorphic on
$\Omega_k\cap\Delta$. Hence by Lemma~\ref{le:1annulus}
below $r_k<\delta_k$ and so $r_k$ converges to zero.

By Lemma~\ref{le:annulus} below, there are
sequences $\eps_k>0$, $\rho_k>\delta_k$, and
$\theta_{0k},\theta_{\infty k}\in[0,2\pi]$
such that $\eps_k$ and $\rho_k$ converge to zero
and
$$
\left|x^{-1} \xi_k(x)-e^{i\theta_{0k}}\right|\le\eps_k,\qquad
\left|y^{-1} \eta_k(y)-e^{i\theta_{\infty k}}\right|\le\eps_k,
$$
for $x,y\in\A(\rho_k,1)$.  To see this let $\delta(m)>0$
be the constant of Lemma~\ref{le:annulus} with $\eps=\rho=1/m$,
choose an increasing sequence of integers $k_m$
such that $\delta_k\le\delta(m)$ for $k\ge k_m$,
and define $\eps_k:=\rho_k:=1/m$ for $k_m\le k<k_{m+1}$.
We call this kind of argument {\em proof by patience}.
It follows that the maps $x\mapsto e^{-i\theta_{0k}}\xi_k(x)$
and $y\mapsto e^{-i\theta_{\infty k}}\eta_k(y)$
converge to the identity uniformly with all derivatives
on every compact subset of $\INT(\D)\setminus0$.

Next we construct two sequences of diffeomorphisms
$\alpha_k,\beta_k:\D\to\D$, converging
to the identity in the $\Cinf$ topology,
and an exhausting sequence of closed annuli
$B_k\subset\INT(\D)\setminus 0$, such that
$\alpha_k$ and $\beta_k$ are equal to the
identity in a neighborhood of $S^1=\p\D$  and
$$
\xi_k(\alpha_k(x))=e^{i\theta_{0k}}x,\qquad
\eta_k(\beta_k(y))=e^{i\theta_{\infty k}}y
$$
for $x,y\in B_k$.  The assertion is obvious by an interpolation
argument when the sequence $B_k$ is replaced by a single
closed annulus $B\subset\INT(\D)\setminus 0$.
Now argue by patience as above.

Increasing $\delta_k$ if necessary we may assume that
$A(\delta_k,\half)\subset B_k$ for all $k$.
Denote $\Omega_k:=\Sigma\setminus\INT(\Delta(\delta_k))$
as above.  Now replace $\xi_k$ by $\xi_k\circ\alpha_k$
and $\eta_k$ by $\eta_k\circ\beta_k$.
These functions satisfy~6) and~7).  Redefine
$f_k$ so that~4) and~5) hold with the new definitions
of $\xi_k$ and $\eta_k$.  Thus we have  proved the claim.

\medskip

In the following we assume w.l.o.g. that $(\pi:Q\to B,S_*,b_0)$
is the universal unfolding of $(\Sigma,s_*,\nu,j)$
constructed in the proof of Theorem~\ref{thm:existence}.
Now define the marked Riemann surface $(\Sigma_k',s_{k,*}',j_k')$  by
$$
\Sigma_k':=\Sigma\setminus\INT(\Delta(\sqrt{r_k}))
$$
where
$j_k':=f_k^*j_k$,
$s_{k,i}':=f_k^{-1}(s_{k,i})$, and
where the equivalence relation is defined by
$$
p\sim q\qquad\iff\qquad x(p)y(q)=z_k, \qquad
z_k := r_ke^{-i(\theta_{0k}+\theta_{\infty k})}
$$
for $p\in \Delta_0$ and $q\in \Delta_\infty$ with $|x(p)|=|y(q)|=\sqrt{r_k}$.
Then, after removing finitely many terms,
there is a sequence of regular values $b_k\in B$ of $\pi:Q\to B$
and a sequence of desingularizations $u_k':\Sigma_k'\to Q_{b_k}$
such that $u_k'(s_{k,i}')=S_i\cap Q_{b_k}$ and $j_k'$
is the pullback of the complex structure on $Q_{b_k}$ under $u_k'$.
This  follows from Theorem~\ref{thm:teich-mark}
and the construction of a universal unfolding in
the proof of Theorem~\ref{thm:existence}.
Moreover, $(\Sigma_k',s_{k,*}',j_k')$ is isomorphic to
$(\Sigma_k,s_{k,*},j_k)$.  An explicit isomorphism is the
map $\psi_k:\Sigma_k'\to\Sigma_k$ defined by
$${
\renewcommand{\arraystretch}{1.9}
\psi_k(p):=
\left\{\begin{array}{ll}
f_k(p) &
 \mbox{for }p\in\Sigma_k'\setminus \Delta,\\
h_k(\xi_k(x(p))) &
 \mbox{for }p\in \Delta_0\mbox{ with } \delta_k\le|x(p)|\le 1,\\
h_k(e^{i\theta_{0k}}x(p)) &
  \mbox{for }p\in \Delta_0\mbox{ with } \sqrt{r_k}\le|x(p)|\le\delta_k,\\
\ds h_k\left(\frac{r_k}{\eta_k(y(p))}\right) &
   \mbox{for } p\in \Delta_\infty\mbox{ with }\delta_k\le|y(p)|\le 1,\\
\ds h_k\left(\frac{r_k}{e^{i\theta_{\infty k}}y(p)}\right) &
  \mbox{for }p\in \Delta_\infty\mbox{ with } \sqrt{r_k}\le|y(p)|\le\delta_k.
\end{array}\right.
}$$
That~(i) implies~(ii) in the case of a single node follows
immediately with
$$
u_k:=u_k'\circ\psi_k^{-1}:\Sigma_k\to Q_{b_k}.
$$
The case of several nodes is analogous.
This proves Theorem~\ref{thm:convergence}.
\end{proof}

\rmk
The sequence $u_k$ just constructed  is such that
$u_k(\gamma_k)$ converges to the nodal set in $Q_{b_0}$.
To prove this, note that
$$
\psi_k^{-1}(\gamma_k)\subset
\{[p]=[q]\in\Sigma_k'\,:\,\sqrt{r_k}\le |x(p)|,|y(q)|\le\delta_k\}.
$$
Hence, by Step~1, $u_k(\gamma_k)$ converges to the nodal point
in $Q_{b_0}$.
\kmr

\begin{lemma}\label{le:1annulus}
If there is a holomorphic map $f:\A(r_1,R_1)\to\A(r_2,R_2)$
inducing an isomorphism of fundamental groups,
then $R_1/r_1\le R_2/r_2$.
\end{lemma}

\begin{proof}
The result is due to Huber~\cite{HUBER}; an exposition appears
in~\cite[Theorem~6.1, page~14]{KOBAYASHI}.
The proof uses the Schwarz Pick Ahlfors Lemma
(a holomorphic map from the unit disk to itself is a contraction in the
Poincar\'e metric). The circle of radius $\sqrt{r_1R_1}$ is a geodesic
in the hyperbolic metric of length $2\pi^2/\log(R_1/r_1)$;
its image under $f$ is shorter and hence so is the central
geodesic in $\A(r_2,R_2)$.
\end{proof}

\begin{lemma}\label{le:annulus}
For every $\eps>0$ and every $\rho>0$ there is a constant
$\delta\in(0,\rho)$ such that the following holds.
If $u:\A(\delta,1)\to\D\setminus 0$
is a holomorphic embedding such that $u(S^1)=S^1$
then there is a real number $\theta$ such that
$$
x\in\A(\rho,1)\qquad\implies\qquad
\left|x^{-1}u(x)-e^{i\theta}\right|<\eps.
$$
\end{lemma}

\begin{proof}
It suffices to assume $u(1)=1$ and then prove the claim with
$\theta=0$.  Suppose, by contradiction that there exist
constants $\eps>0$ and $\rho>0$ such that the assertion
is wrong. Then there exists a sequence $\delta_i>0$ converging
to zero and a sequence of holomorphic embeddings
$u_i:\A(\delta_i,1)\to\D\setminus 0$ such that
$$
u_i(S^1)=S^1,\qquad u_i(1)=1,\qquad
\sup_{\rho\le|x|\le1}\left|u_i(x)-x\right|\ge\eps\rho.
$$
We claim that $u_i$ converges to the identity, uniformly
on every compact subset of $\D\setminus 0$.  To see this extend
$u_i$ to the annulus $\A(\delta_i,1/\delta_i)$ by the formula
$$
u_i(z) := \frac{1}{\overline{u_i(1/\bar z)}}
$$
for $1\le|z|\le 1/\delta_i$. Think of the extended map
as a holomorphic embedding
$u_i:\A(\delta_i,1/\delta_i)\to S^2\setminus\{0,\infty\}$.
Next we claim that
\begin{equation}\label{eq:bounded}
\sup_i\sup_{z\in K}|du_i(z)|<\infty
\end{equation}
for every compact subset $K\subset\C\setminus0$.
Namely, the energy of the holomorphic curve $u_i$ is bounded by the area
of the target manifold $S^2$.  So if $|du_i(z_i)|\to\infty$
for some sequence $z_i\to z_0\in\C\setminus0$, then a holomorphic
sphere bubbles off near $z_0$ and it follows that a subsequence of
$u_i$ converges to a constant, uniformly on every compact
subset of $\C\setminus\{0,z_0\}$.  But this contradicts the fact
that $u_i(S^1)=S^1$.  Thus we have proved~(\ref{eq:bounded}).
Now it follows from the standard elliptic bootstrapping techniques
(or alternatively from Cauchy's integral formula and the
Arzela--Ascoli theorem) that there is a subsequence,
still denoted by $u_i$, that converges in the
$\Cinf$ topology to a holomorphic curve
$u_0:\C\setminus 0\to S^2\setminus\{0,\infty\}$.
By the removable singularity theorem, $u_0$ extends to a
holomorphic curve $u_0:S^2\to S^2$.  Since $u_i$ is an embedding
for every $i$, it follows that $u_0$ is an embedding and hence
a M\"obius transformation.  Since $u_i(S^1)=S^1$, $u_i(1)=1$
and $0\notin u_i(\A(\delta_i,1))$, it follows that
$$
u_0(S^1)=S^1,\qquad u_0(1)=1,\qquad u_0(0)=0.
$$
This implies that $u_0=\id$.  Thus we have proved that $u_i$
converges to the identity, uniformly on every compact subset
of $\D\setminus 0$.  This contradicts the inequality
$\sup_{\rho\le|x|\le1}\left|u_i(x)-x\right|\ge\eps\rho$
and this contradiction proves the lemma.
\end{proof}


\section{Compactness}\label{sec:compact}

In this section we prove that every sequence of stable marked nodal
Riemann surfaces of type $(g,n)$ has a DM-convergent subsequence.
Our strategy is to perform some preliminary constructions to reduce
our compactness theorem to Proposition~5.1 of~\cite[page 71]{HUMMEL}.
We begin by rephrasing Hummel's result in a weaker form that we will
apply directly (see Proposition~\ref{prop:hummel} below).

\para\label{para:finiteType}
Let $W$ be a smooth oriented surface,
possibly with boundary and not necessarily compact or connected.
A \jdef{finite extension} of $W$ is a smooth orientation preserving
embedding $\iota:W\to S$ into a compact oriented surface $S$
such that $S\setminus\iota(W)$ is finite. If $\iota_1:W\to S_1$
and $\iota_2:W\to S_2$ are two such extensions,
the map $\iota_2\circ\iota_1^{-1}$ extends to  a homeomorphism,
but not necessarily to a diffeomorphism.
Let $W_1,\ldots,W_\ell$ be the components of $W$,
$S_1,\ldots,S_\ell$ be the corresponding components of
a finite extension $S$, $g_i$ be the genus of $S_i$,
$m_i$ be the number of boundary components
of $W_i$, and $n_i$ be the number of points in $S_i\setminus\iota(W_i)$.
The (unordered) list $(g_i,m_i,n_i)$ is called the \jdef{signature} of $W$.
Two surfaces of finite type are diffeomorphic if and only if they have the
same signature. (Compare~\ref{def:signature}  and~\ref{rmk:signature}.)
We say that $W$ is of \jdef{stable type} if $n_i>\chi(S_i)$
(at least one puncture point on an annulus or torus,
at least two on a disk, and at least three on a sphere).
\arap

\para\label{hyperbolic}
A \jdef{hyperbolic metric}  on $W$ is a complete Riemannian metric
$h$ of constant curvature $-1$ such that each boundary component
is a closed geodesic. A \jdef{finite extension} of a complex structure $j$
on $W$ is a finite extension $\iota:W\to S$ such that $\iota_*j$
extends to a complex structure on $S$; we say $j$ has \jdef{finite type}
if it admits a finite extension.
\arap

\begin{proposition}\label{prop:hyperbolic}
Let $W$ be a surface of stable type.
Then the operation which assigns to each hyperbolic metric
on $W$ its corresponding complex structure (rotation by $90^\circ$)
is bijective. It restricts to  a bijection between hyperbolic metrics
of finite area and complex structures of finite type.
\end{proposition}

\begin{proof}
The operation $h\mapsto j$ is injective by the removable
singularities theorem and surjective by applying the
uniformization theorem to the holomorphic double.
If $j$ is of finite type, then
the area is finite by~\cite[Proposition~3.9 page~68]{HUMMEL}.
If $h$ is of finite area, then
$j$ is of finite type by~\cite[Proposition~3.6 page~65]{HUMMEL}.
\end{proof}

\begin{proposition}[Mumford--Hummel]\label{prop:hummel}
Let $S$ be a compact connected surface with boundary and
$x_1,\dots,x_n$ be a sequence of marked points in the interior
of $S$ such that $W:=S\setminus\{x_1,\dots,x_n\}$ is of stable type.
Denote
$$
\p S =: \p_1S\cup\cdots\cup\p_mS,
$$
where each $\p_iS$ is a circle.
Let $j_k$ be a sequence of complex structures on $S$
and $h_k$ be the corresponding sequence of hyperbolic
metrics on $W$. Assume:
\begin{description}
\item[(a)]
The lengths of the closed geodesics in $W\setminus\p W$
are bounded away from zero.
\item[(b)]
The lengths of the boundary geodesics converge to zero.
\end{description}
Then there exists a subsequence, still denoted by $(j_k,h_k)$,
a closed Riemann surface $(\Sigma,j)$ with distinct
marked points $\xi_1,\dots,\xi_n,\eta_1,\dots,\eta_m$,
a hyperbolic metric $h$ of finite area
on $\Sigma\setminus\{\xi_1,\dots,\xi_n,\eta_1,\dots,\eta_m\}$,
and a sequence of continuous maps $\phi_k:S\to\Sigma$
satisfying the following conditions.
\begin{description}
\item[(i)]
$\phi_k(x_i)=\xi_i$ for $i=1,\dots,n$ and $\phi_k(\p_iS)=\eta_i$
for $i=1,\dots,m$.
\item[(ii)]
The restriction of $\phi_k$ to $S\setminus\p S$ is a diffeomorphism
onto $\Sigma\setminus\{\eta_1,\dots,\eta_m\}$.
\item[(iii)]
$(\phi_k)_*j_k$ converges to $j$ on
$\Sigma\setminus\{\eta_1,\dots,\eta_m\}$.
\item[(iv)]
$(\phi_k)_*h_k$ converges to $h$ on
$\Sigma\setminus\{\xi_1,\dots,\xi_n,\eta_1,\dots,\eta_m\}$.
\end{description}
\end{proposition}

\begin{proof}
This follows from Proposition~5.1 in~\cite[page~71]{HUMMEL}.
The discussion preceding Proposition~5.1  in~\cite{HUMMEL}
explains how to extract the subsequence and how to construct
the Riemann surface $(\Sigma,j)$ and the hyperbolic metric $h$.
\end{proof}

\begin{theorem}\label{thm:compact}
Every sequence of stable marked nodal Riemann surfaces of
type $(g,n)$ has a DM-convergent subsequence.
\end{theorem}

\begin{proof}
Let $(\Sigma_k,s_{k,*},\nu_k,j_k)$ be a sequence of marked nodal
Riemann surfaces of type $(g,n)$. Passing to a subsequence, if necessary,
we may assume that all marked nodal surfaces in our sequence
have the same signature (see Definition~\ref{def:signature})
and hence are diffeomorphic.  Thus we assume that
$$
(\Sigma_k,s_{k,*},\nu_k)=(\Sigma,s_*,\nu)
$$
is independent of $k$. Denote by $\Sigma^*$
the possibly disconnected and noncompact surface obtained from
$\Sigma$ by removing the special points. Let $h_k$ be the hyperbolic
metric on $\Sigma^*$ determined by $j_k$ (see Proposition~\ref{prop:hyperbolic}).

Let $\ell_k^1$ be the length of the shortest geodesic in $\Sigma^*$
with respect to $h_k$.  If a subsequence of the $\ell_k^1$ is bounded away
from zero we can apply Proposition~\ref{prop:hummel} to each component
of $\Sigma$ and the assertion follows. Namely, the maps
$\phi_k$ in Proposition~\ref{prop:hummel} are deformations
as in Definition~\ref{deformation}.

Hence assume $\ell_k^1$ converges to zero as $k$ tends to infinity
and, for each $k$,  choose a geodesic $\gamma_k^1$ with length
$\ell_k^1$.  Passing to a further subsequence and, if necessary,
modifying $h_k$ by a diffeomorphism that fixes the marked and
nodal points we may assume that the geodesics $\gamma_k^1$
are all homotopic and indeed equal.  Thus
$$
\gamma_k^1 = \gamma^1
$$
for every $k$. Now let $\ell_k^2$ be the length of the shortest
geodesic in $\Sigma\setminus\gamma^1$ with respect to $h_k$.
If a subsequence of $\ell_k^2$ is bounded away from zero
we cut open $\Sigma$ along $\gamma^1$.  Again the assertion
follows by applying Proposition~\ref{prop:hummel} to each component
of the resulting surface with boundary.

Continue by induction. That the induction terminates follows from
the fact that the geodesics in $(\Sigma^*,h_k)$ of lengths at most
$2\mathrm{arcsinh}(1)$ are pairwise disjoint and their number is
bounded above by $3g-3+N$, where $N$ is the number of
special points (see~\cite[Lemma~4.1 page~68]{HUMMEL}).
This proves the theorem.
\end{proof}

\begin{lemma}\label{le:vanishingcycle}
Let $(\pi:P\to A,R_*,a_0)$ be a nodal unfolding and $C\subset P$ be
the set of critical points of $\pi$.  Then, after shrinking $A$ if necessary,
there exists a closed subset $V\subset P$ and a smooth map
$$
\rho:P\setminus V\to P_{a_0}\setminus V
$$
satisfying the following conditions.
\begin{description}
\item[(i)]
For every $a\in A$ we have $C\cap P_a\subset V\cap P_a=:V_a$;
moreover $C\cap P_{a_0}=V_{a_0}$.
\item[(ii)]
Each component of $V$ intersects $P_a$
either in a simple closed curve or in a nodal point.
\item[(iii)]
For every $a\in A$ the restriction
$\rho_a:=\rho|P_a\setminus V_a:
P_a\setminus V_a\to P_{a_0}\setminus V_{a_0}$
is a diffeomorphism; moreover $\rho_{a_0}=\id$.
\end{description}
\end{lemma}

\begin{proof}
Choose a Hardy trivialization $(P=M\cup N,\iota,\rho)$
as in~\ref{def:hardyTrivialization} and write
$$
N = N_1\cup\cdots\cup N_k.
$$
as in Definition~\ref{def:hardyDecomposition}.
Let $(z_i,t_i):A\to U_i\subset\D\times\C^{m-1}$ and
$(x_i,y_i,t_i):N_i\to\D^2\times\C^{m-1}$ be the holomorphic
coordinates of Definition~\ref{def:hardyDecomposition}
so that $z_i(\pi(p))=x_i(p)y_i(p)$ and the critical
set $C\subset P$ has components
$$
C_i:=\left\{p\in N_i\,:\,x_i(p)=y_i(p)=0\right\}.
$$
Define
$$
V := V_1\cup\cdots\cup V_k,\qquad
V_i:=\left\{p\in N_i\,:\,|x_i(p)|=|y_i(p)|=\sqrt{|z_i(\pi(p))|}\right\}.
$$
This set satisfies~(i) and~(ii). The restriction of the
trivialization $\rho:M\to M_{a_0}$ to
$\p N_i\subset \p N=\p M$ is, in the above coordinates,
given by $\rho(x_i,y_i,t_i)=(x_i,0,t_i)$ for $|x_i|=1$
and by $\rho(x_i,y_i,t_i)=(0,y_i,t_i)$ for $|y_i|=1$.
We extend this map by an explicit formula.
Choose a smooth cutoff function $\beta:[1,\infty)\to[0,1]$
such that $\beta'(r)\ge 0$ for every $r$ and
$$
\beta(r) := \left\{\begin{array}{rl}
r-1,&\mbox{for } 1\le r\le 3/2,\\
1,&\mbox{for } r\ge 2.
\end{array}\right.
$$
Then define the extension $\rho:N_i\setminus V_i\to P_{a_0}$
in local coordinates by
$$
\rho(x_i,y_i,t):=\left\{\begin{array}{rl}
\left(\beta\left(\sqrt{|x_i|/|y_i|}\right)x_i,0,t_i\right),&\mbox{if }|x_i|>|y_i|,\\
\left(0,\beta\left(\sqrt{|y_i|/|x_i|}\right)y_i,t_i\right),&\mbox{if }|y_i|>|x_i|.
\end{array}\right.
$$
The resulting map $\rho:P\setminus V\to P_{a_0}$ is smooth
and satisfies~(iii).  This proves the lemma.
\end{proof}

\begin{proof}[Proof of Theorem~\ref{thm:proper}.]
Let $(\pi:Q\to B,S_*)$ be a universal family and denote by $(B,\Gamma)$
the associated etale groupoid of Definition~\ref{B-Gamma}
(see Theorem~\ref{thm:Gamma}).   We prove that this groupoid is proper.
Thus let $(a_k,f_k,b_k)$ be a sequence in $\Gamma$ such that
$a_k$ converges to $a_0$ and $b_k$ converges to $b_0$.
We must show that there is a fiber isomorphism $f_0:Q_{a_0}\to Q_{b_0}$
such that a suitable subsequence of $f_k$ DM-converges to $f_0$
(see Definition~\ref{def:proper}). To see this choose desingularizations
$$
\iota:\Sigma\to Q_{a_0},\qquad \iota':\Sigma'\to Q_{b_0}.
$$
Denote by $(\Sigma,s_*,\nu,j)$ and $(\Sigma',s'_*,\nu',j')$
the induced marked nodal Riemann surfaces. Consider the
following diagram
$$
\xymatrix{
&{Q_{a_k}\setminus V_{a_k}}\ar[r]^{f_k}\ar[d]^{\rho_{a_k}}
&{Q_{b_k}\setminus V_{b_k}}\ar[d]^{\rho_{b_k}}\\
{\Sigma\setminus\nu}\ar[r]^{\iota}\ar[ru]^{\iota_k}
&{Q_{a_0}\setminus V_{a_0}}
&{Q_{b_0}\setminus V_{b_0}}
&{\Sigma'\setminus\nu'}\ar[l]_{\iota'}\ar[lu]_{\iota_k'}
}
$$
Here the sets $V_a:=V\cap Q_a$ and the diffeomorphisms
$\rho_a:Q_a\setminus V_a\to Q_{a_0}\setminus V_{a_0}$
are as in Lemma~\ref{le:vanishingcycle} for $a$ near $a_0$;
similarly for $b$ near $b_0$. Moreover,
$$
\iota_k := \rho_{a_k}^{-1}\circ\iota,\qquad
\iota'_k := \rho_{b_k}^{-1}\circ\iota'.
$$
By definition, the pullback complex structures
$$
j_k:=\iota_k^*J|Q_{a_k},\qquad
j_k':={\iota_k'}^*J|Q_{b_k}
$$
converge to $j$, respectively $j'$, in the $\Cinf$ topology
on every compact subset $\Sigma\setminus\nu$,
respectively $\Sigma'\setminus\nu'$.
By Lemma~\ref{le:vanishingcycle}, there exist exhausting
sequences of open sets
$$
U_k\subset\Sigma\setminus\nu,\qquad
U_k'\subset\Sigma'\setminus\nu',\qquad
f_k(U_k)\subset U_k',
$$
such that $j_k$ can be modified outside $U_k$ so as to
converge in the $\Cinf$ topology on all of $\Sigma$ to $j$,
and similarly for $j_k'$. Then
$$
u_k := (\iota_k')^{-1}\circ f_k\circ\iota_k:U_k\to\Sigma'
$$
is a sequence of $(j_k,j_k')$-holomorphic embeddings such that
$u_k(s_*)=s'_*$.  The argument in Remark~\ref{rmk:proper}
shows that, if the first derivatives of $u_k$ are uniformly bounded,
then $u_k$ has a $\Cinf$ convergent subsequence.
It also shows that a nonconstant holomorphic
sphere in $Q$ bubbles off whenever the first
derivatives of $u_k$ are not bounded.  But bubbling cannot
occur (in $\Sigma\setminus\nu$).  To see this argue as follows.
Suppose $z_k$ converges to $z_0\in\Sigma\setminus\nu$ and
the derivatives of $u_k$ at $z_k$ blow up.  Then the standard
bubbling argument (see~\cite[Chapter~4]{MS2}) applies.
It shows that, after passing to a subsequence and modifying
$z_k$ (without changing the limit), there are $(i,j_k)$-holomorphic
embeddings $\eps_k$ from the disk $\D_k\subset\C$, centered
at zero with radius $k$, to $\Sigma$ such that $\eps_k(0)=z_k$,
the family of disks $\eps_k(\D_k)$ converges to $z_0$,
and $u_k\circ\eps_k$ converges to a nonconstant $J$-holomorphic
sphere $v_0:S^2=\C\cup\infty\to Q_{b_0}$. (The convergence is
uniform with all derivatives on every compact subset of $\C$.)
Hence the image of $v_0$ contains at least three special points.
It follows that the image of $u_k\circ\eps_k$ contains at least two
special points for $k$ sufficiently large.
But the image of $u_k$ contains no nodal points,
the image of $\eps_k$ contains at most one marked point,
and $u_k$ maps the marked points of $\Sigma$
bijectively onto the marked points of $\Sigma'$.
Hence the image of $u_k\circ\eps_k$ contains
at most one special point, a contradiction.

This shows that bubbling cannot occur, as claimed,
and hence a suitable subsequence of $u_k$ converges
in the $\Cinf$-topology to a $(j,j')$-holomorphic curve
$u_0:\Sigma\setminus\nu\to\Sigma'\setminus\nu'$.
Now the removable singularity theorem shows that
$u_0$ extends to a $(j,j')$-holomorphic curve on all of $\Sigma$
and maps $\nu$ to $\nu'$. That $u_0$ is bijective
follows by applying the same argument to $f_k^{-1}$.
Hence there exists a unique fiber isomorphism $f_0:Q_{a_0}\to Q_{b_0}$
such that $\iota'\circ u_0=f_0\circ\iota$. By construction,
the subsequence of $f_k$ DM-converges to $f_0$

Thus we have proved that the map $s\times t:\Gamma\to B\times B$
is proper. Hence, by Corollary~\ref{cor:properHausdorff},
the quotient space $B/\Gamma$ is Hausdorff.
Moreover, by Theorems~\ref{thm:convergence}
and~\ref{thm:compact}, it is sequentially compact.
Since $B$ is second countable it follows that $B/\Gamma$
is compact. This completes the proof of Theorem~\ref{thm:proper}.
\end{proof}

\begin{corollary}\label{cor:hausdorff}
Suppose that a sequence of marked nodal Riemann surfaces
of type $(g,n)$ DM-converges to both $(\Sigma,s_*,\nu,j)$ and $(\Sigma',s'_*,\nu',j')$.
Then $(\Sigma,s_*,\nu,j)$ and $(\Sigma',s'_*,\nu',j')$ are isomorphic.
\end{corollary}

\begin{proof}
Let $(\pi_B:Q\to B,S_*)$ be a universal family.
By Theorem~\ref{thm:convergence} there exist points $a_0,b_0\in B$
and sequences $a_k\to a_0$, $b_k\to b_0$ such that
$(\Sigma,s_*,\nu,j)$  arises from a desingularization of $Q_{a_0}$,
$(\Sigma',s'_*,\nu',j')$ arises from a desingularization of $Q_{b_0}$,
and the fibers $Q_{a_k}$ and $Q_{b_k}$ are isomorphic.
Hence, by Theorem~\ref{thm:proper}, there exists
a fiber isomorphism from $Q_{a_0}$ to $Q_{b_0}$,
and so $(\Sigma,s_*,\nu,j)$ and $(\Sigma',s'_*,\nu',j')$
are isomorphic.
\end{proof}



\begin{thebibliography}{99}
\small

\bibitem{CGS}
K.~Cieliebak, A.R.~Gaio, \& D.A.~Salamon:
J-holomorphic curves, moment maps,
and invariants of Hamiltonian group actions,
{\em IMRN} {\bf 10} (2000), 831-882.

\bibitem{EARLE}
C.~Earle:
On holomorphic cross sections in Teichm\"uller spaces,
Duke Math. Journal, {\bf 36} (1969) 409--415.

\bibitem{EE}
C.~Earle \& J.~Eells:
A fibre bundle approach to Teichm\"uller theory,
{\it J. Differential Geometry}, {\bf 3} (1969) 19--43.

\bibitem{DISCRETE}
W.J.~Harvey (ed.):
{\em Discrete Groups and Automorphic Functions},
Academic Press, 1977.

\bibitem{GN}
A.~Friedman:
{\it Partial Differential Equations},
Holt, Rinehart, and Winston, 1969.

\bibitem{GROTHENDIECK}
A.~Grothendieck:
Familles d'espaces complexes et fondements de la g\'eom\'etrie analytique,
S\'eminaire Henri Cartan, 13i\`eme ann\'ee: 1960/61.  Fasc. 1, Exp. 7, 9--13,
deuzime dition, Secrtariat mathmatique, Paris, 1962; Fasc. 2, Exp. 14--17,
deuzime dition, Secrtariat mathmatique, Paris, 1962.

\bibitem{GROT}
A.~Grothendieck:
Sur la classification des fibr\'es holomorphes sur la sph\`ere
de Riemann, {\it Amer. J. Math.} {\bf 76}, (1957), 121--138.

\bibitem{HUBER} H.~Huber:
\"Uber analytische Abbildungen von Ringgebieten
in Ringgebiete, {\it Compos. Math.}  {\bf 9} (1951), 161--168.

\bibitem{HUMMEL}
C.~Hummel:
{\em Gromov's Compactness Theorem for Pseudo-holomorphic Curves},
Birkh\"auser Progress in Mathematics {\bf 151}, 1997.

\bibitem{KEEL}
S.~Keel:
Intersection theory of moduli spaces of $n$-pointed curves of genus zero.
{\it Trans. Amer. Math. Soc} {\bf 330} (1992), 545--574.

\bibitem{KELLEY}
J.L.~Kelley:
{\em General Topology}, Van Nostrand, 1955.

\bibitem{KNUDSEN}
F.F.~Knudsen:
The projectivity of moduli spaces of stable curves II,
{\it Funct. Anal. Appl.} {\bf 19} (1983), 161--199.

\bibitem{KOBAYASHI} S.~Kobayashi:
{\it Hyperbolic Manifolds and Holomorphic Mappings},
Marcell Dekker, 1970.

\bibitem{MS1}
D.~McDuff \& D.A.~Salamon:
{\it Introduction to Symplectic Topology},
2nd edition, Oxford University Press, 1998.

\bibitem{MS2}
D.~McDuff \& D.A.~Salamon:
{\it J-holomorphic Curves and Symplectic Topology},
AMS Colloquium Publications, 2004.

\bibitem{SMALE} S.~Smale:
Diffeomorphisms of the $2$-sphere,
{\em Proc. Amer. Math. Society}, {\bf 10} (1959) 621--626.

\bibitem{PS} A.~Pressley \& G.~Segal:
{\em Loop Groups}, Oxford Mathematical Monographs,
The Clarendon Press Oxford University Press, 1986.

\bibitem{TROMBA}
A.~Tromba:
{\em Teichm\"uller Theory in Riemannian Geometry},
Lectures in Mathematics ETH Z\"urich, Birkh\"auser, 1992.

\end{thebibliography}
\end{document}